\documentclass{amsart}
\usepackage[dvips]{epsfig}
\usepackage{graphicx}
\usepackage{latexsym}
\usepackage{amsmath}
\usepackage{amsthm}
\usepackage{amssymb}
\usepackage{mathrsfs}
\usepackage[shortlabels]{enumitem}
\newlist{propenum}{enumerate}{1} 
\setlist[propenum]{label=(\roman*)}

\usepackage{caption}
\usepackage{subcaption}

\usepackage{eucal}

\usepackage{xcolor,soul}

\usepackage{mathtools}

 \usepackage[applemac]{inputenc}

\newtheorem{thm}{Theorem}[section]

\newtheorem{lem}[thm]{Lemma}

\newtheorem{defi}[thm]{Definition}
\newtheorem{hyp}[thm]{Assumption}

\theoremstyle{remark}

\newtheorem{rem}[thm]{Remark}

\newcommand{\inL} {\underset{n \rightarrow \infty}{\overset{\rm (d)}{\xrightarrow{\hspace*{0.75cm}}}} }

\newcommand{\vip}{\vskip.2cm}
\newcommand{\COMMENTAIRE}[1]{}

\newcommand{\field}[1]{\mathbb{#1}}
\newcommand{\EE}{\field{E}}

\newcommand{\GG}{\field{G}}

\newcommand{\NN}{\field{N}}

\newcommand{\RR}{\field{R}}
\newcommand{\TT}{\field{T}}

\newcommand{\XX}{\field{X}}

\newcommand{\Cc}{{\mathcal C}}

\newcommand{\Nn}{{\mathcal N}}
\newcommand{\Pp}{{\mathcal P}}

\newcommand{\Qq}{{\mathcal Q}}

\def \ep {\varepsilon}

\newcommand{\rd}{{\rm d}}

\newcommand{\bF}{{\mathfrak f}}
\newcommand{\bC}{{\mathfrak c}}
\newcommand{\bH}{{\mathfrak h}}
\newcommand{\bQ}{{\mathfrak q}}

\newcommand{\cb}{{\mathcal B}}
\newcommand{\cc}{{\mathcal C}}

\newcommand{\ce}{{\mathcal E}}

\newcommand{\cf}{{\mathcal F}}

\newcommand{\cn}{{\mathcal N}}

\newcommand{\cp}{{\mathcal P}}
\newcommand{\cq}{{\mathcal Q}}

\newcommand{\cs}{{\mathscr S}}

\newcommand{\A}{{\mathbb A}}

\newcommand{\E}{{\mathbb E}}

\newcommand{\G}{\mathbb{G}}
\newcommand{\N}{{\mathbb N}}

\newcommand{\R}{{\mathbb R}}

\newcommand{\T}{\mathbb{T}}

\newcommand{\ind}{{\bf 1}}

\newcommand{\Var}{{\mathrm{Var}}}

\newcommand{\sot}{\otimes_{\rm sym}}

\newcommand{\norm}[1]{\mathop{\parallel\! #1 \! \parallel}\nolimits}
\newcommand{\normm}[1]{\mathop{\parallel\! #1 \! \parallel}\nolimits_{L^{2}(\mu)}}

\newcommand{\inv}[1]{\mathop{\frac{1}{ #1}}\nolimits}
\newcommand{\expp}[1]{\mathop {\mathrm{e}^{ #1}}}

\newcommand{\reff}[1]{(\ref{#1})}

\begin{document}

\title[clt for kernel density estimator of BMC]{Central limit theorem for kernel estimator of invariant density in bifurcating Markov chains models.}

\author{S. Val\`ere Bitseki Penda and Jean-Fran\c cois Delmas}

\address{S. Val\`ere Bitseki Penda, IMB, CNRS-UMR 5584, Universit\'e Bourgogne Franche-Comt\'e, 9 avenue Alain Savary, 21078 Dijon Cedex, France.}

\email{simeon-valere.bitseki-penda@u-bourgogne.fr}

\address{Jean-Fran\c cois Delmas, CERMICS, Ecole des Ponts,  France.}

\email{jean-francois.delmas@enpc.fr}
\date{\today}

\begin{abstract}
  Bifurcating Markov  chains (BMC) are  Markov chains indexed by  a full
  binary tree representing  the evolution of a trait  along a population
  where each  individual has two  children. Motivated by  the functional
  estimation of the  density of the invariant  probability measure which
  appears  as the  asymptotic distribution  of the  trait, we  prove the
  consistence and  the Gaussian fluctuations  for a kernel  estimator of
  this  density based  on  late  generations.  In  this  setting, it  is
  interesting to note  that the distinction of the three  regimes on the
  ergodic  rate  identified in  a  previous  work (for  fluctuations  of
  average over large generations) disappears. This result is a first step
  to go  beyond the  threshold condition  on the  ergodic rate  given in
  previous statistical papers on  functional estimation.
\end{abstract}

\thanks{We warmly thank Ad\'ela\"ide Olivier for numerous discussions in the preliminary phases of this project. }

\maketitle

\textbf{Keywords}: Bifurcating Markov chains, bifurcating
auto-regressive process, binary trees, fluctuations for tree indexed
Markov chain, density estimation.\\

\textbf{Mathematics Subject Classification (2020)}: 62G05, 62F12, 60J05, 60F05, 60J80.




\section{Introduction}

Bifurcating  Markov chains  (BMC) are  a class  of stochastic  processes
indexed by  regular binary tree  and which satisfy the  branching Markov
property (see below for a precise definition). This model represents the
evolution of  a trait along a  population where each individual  has two
children.   The  recent  study  of  BMC  models  was  motivated  by  the
understanding  of the  cell division  mechanism (where  the trait  of an
individual is given by its growth  rate).  The first model of BMC, named
``symmetric''  bifurcating auto-regressive  process  (BAR), see  Section
\ref{sec:sym-bar}  for  more  details  in  a  Gaussian  framework,  were
introduced  by Cowan  \& Staudte  \cite{CS86} in  order to  analyze cell
lineage data.   In \cite{Guyon}, Guyon  has studied more general
asymmetric  BMC to
prove statistical  evidence of  aging in Escherichia  Coli. We  refer to
\cite{BD}  for  more detailed  references  on  this subject.   Recently,
several statistical  works have been  devoted to the estimation  of cell
division rates, see Doumic, Hoffmann, Krell \& Roberts \cite{dhkr:segf},
Bitseki, Hoffmann  \& Olivier  \cite{bho:aebmc} and Hoffmann  \& Marguet
\cite{MR4025704}. Moreover,  another studies, such  as Doumic, Escobedo
\& Tournus \cite{MR3906858},
can be  generalized using  the BMC  theory (we  refer to  the conclusion
therein).  \medskip


In this  paper, our objective is  to study the functional  estimation of
the density of the invariant probability measure $\mu$ associated to the
BMC.  For this purpose, we develop a kernel estimation in the $L^2(\mu)$
framework under reasonable hypothesis (which are in particular satisfied
by    the     Gaussian    symmetric     BAR    model     from    Section
\ref{sec:sym-bar}). This  approach is  in the  spirit of  the $L^2(\mu)$
approach developed \cite{BD2}.  In BMC model, the evolution of the trait
along the genealogy of an individual  taken at random is Markovian.  Let
us assume  it is  geometrically ergodic with  rate $\alpha\in  (-1, 1)$,
with $\mu$ is its invariant measure.  In \cite{BD2}, three regimes where
identified  for  the   rate  of  convergence  of   averages  over  large
generations according to  the ergodic rate of  convergence $\alpha$ with
respect  to   the  threshold  $1/\sqrt{2}$.   It   is  interesting,  and
surprising as well, to note that  the distinction of those three regimes
disappears  for the  rate  of convergence  when  considering the  kernel
density   estimation   of   the    density   of   $\mu$,   see   Theorem
\ref{thm:MK}. However, let  us mention that some  further restriction on
the admissible bandwidths  of the kernel estimator are to  be taken into
account  in   the   super-critical   regime  (\emph{i.e.}
$\alpha>    1/\sqrt{2}$),    to     be    precise    see    Condition
\eqref{eq:bandwidth}  which  is  in   force  for  Theorem  \ref{thm:MK}.
Furthermore, we get that estimations using different generations provide
asymptotically independent fluctuations, see Remark \ref{rem:indep} (see
also the  form of the  asymptotic variance in Theorem  \ref{thm:flx} and
Remark  \ref{rem:s=+}  in a  more  general  framework); this  phenomenon
already appear in \cite{DM10}.  The  convergence of the kernel estimator
in Theorem \ref{thm:MK} relies on different type of assumptions:
\begin{itemize}
   \item Geometric ergodic rate $\alpha\in (0, 1)$ of convergence  for the evolution
of the  trait along the genealogy  of an individual taken  at random, see
Assumption \ref{hyp:Q1}.
\item Regularity  (density and integrability conditions) for the
  evolution kernel $\cp$ and the initial distribution   of the
  BMC, see Assumptions \ref{hyp:densite-cp},  and
  \ref{hyp:densite-mu2}. The former is in  the spirit
  of \cite{BD2} (see Assumption \ref{hyp:DenMu} which is a consequence
  of Assumption \ref{hyp:densite-cp} as proven in Section \ref{sec:hyp=hyp}). 
\item
  Regularity (isotropic H\"older regularity) of the density of $\mu$ with respect to the Lebesgue
  measure on $S=\R^d$, see Assumption
  \ref{hyp:estim-tcl}~\ref{item:den-Holder}.
\item Regularity of  the kernel function $K$ and on  the bandwidth given
  in               Assumption              \ref{hyp:K}               and
  Assumption~\ref{hyp:estim-tcl}~\ref{item:K}-\ref{item:bandew}.
  
\item A  condition on the bandwidth given in Equation
  \eqref{eq:bandwidth} which add a further restriction  only in
  the  super-critical regime
  $\alpha>1/\sqrt{2}$. 

\end{itemize}

Eventually, we present some simulations  on the kernel estimation of the
density of $\mu$.   We note that in statistical studies  which have been
done in \cite{dhkr:segf, bho:aebmc,  bitsekiroche2020}, the ergodic rate
of convergence  is assumed to be  less than 1/2, which  is strictly less
than  the  threshold $1/\sqrt{2}$  for  criticality.   Moreover, in  the
latter works, the authors are  interested in the non-asymptotic analysis
of the  estimators.  Now, with the  new perspective given  by the
present results, 
see in particular Remark \ref{rem:bandwidth}, we think that the works in
\cite{dhkr:segf,  bho:aebmc, bitsekiroche2020}  can be  extended to  the
case where the ergodic rate of convergence belongs to $(1/2,1)$.

\medskip
 
The paper is organized as follows. We introduce the BMC model in Section
\ref{sec:BMC0} as  well as the  $L^2$ ergodic assumption. We  define the
kernel estimator  and state the  main results  on the estimation  of the
density  of $\mu$,  see Lemma  \ref{lem:cge-ps-muhat} (consistency)  and
Theorem \ref{thm:MK}  (asymptotic normality), in  Section \ref{sec:Noy}.
The proofs of those result rely  on a general central limit theorem, see
Theorem  \ref{thm:flx}  in  Section \ref{sec:general-res}.   In  Section
\ref{sec:sym-bar}, we  illustrate our results by  studying the symmetric
BAR, and we  provide a numerical study in  Section \ref{sec:simul}.  The
Sections  \ref{sec:p-l-th-MK}-\ref{sec:proof-Scrit-L2} are  dedicated to
the proofs of the main results.

\section{Bifurcating Markov chain (BMC)}
\label{sec:BMC0}
We denote by   $\N$     the   set   of  non-negative   integers   and
$\N^*=   \N   \setminus   \{0\}$. If $(E,  \ce)$ is a  measurable space, then $\cb(E)$  (resp. $\cb_b(E)$,
resp.    $\cb_+(E)$)  denotes   the  set   of  (resp.    bounded,  resp.
non-negative)  $\R$-valued measurable  functions  defined  on $E$.   For
$f\in \cb(E)$, we set $\norm{f}_\infty =\sup\{|f(x)|, \, x\in E\}$.  For
a finite measure $\lambda$ on $(E,\ce)$ and $f\in \cb(E)$ we shall write
$\langle \lambda,f  \rangle$ for  $\int f(x) \,  \rd\lambda(x)$ whenever
this  integral is  well defined, and $\norm{f}_{L^2(\lambda)}=\langle
\lambda , f^2 \rangle^{1/2}$. For $n\in \N^*$, the product space $E^n$ is
endowed with the product $\sigma$-field  $\ce^{\otimes n}$. 
If $(E,d)$  is a  metric space,  then
$\ce $ will denote its Borel $\sigma$-field and the set $\cc_b(E)$ (resp.
$\cc_+(E)$) denotes the set of bounded (resp.  non-negative) $\R$-valued
continuous functions defined on $E$.

Let  $(S, \cs)$   be  a  measurable  space. 
Let $Q$ be a   
probability kernel   on $S \times \cs$, that is:
$Q(\cdot  , A)$  is measurable  for all  $A\in \cs$,  and $Q(x,
\cdot)$ is  a probability measure on $(S,\cs)$ for all $x \in
S$. For any $f\in \cb_b(S)$,   we set for $x\in S$:
\begin{equation}
   \label{eq:Qf}
(Qf)(x)=\int_{S} f(y)\; Q(x,\rd y).
\end{equation}
We define $(Qf)$, or simply $Qf$, for $f\in \cb(S)$ as soon as the
integral \reff{eq:Qf} is well defined, and we have $\cq f\in \cb(S)$. For $n\in \N$, we denote by $Q^n$  the
$n$-th iterate of $Q$ defined by $Q^0=I$, the identity map on
$\cb(S)$, and $Q^{n+1}f=Q^n(Qf)$ for $f\in \cb_b(S)$.  

Let $P$ be a   
probability kernel   on $S \times \cs^{\otimes 2}$, that is:
$P(\cdot  , A)$  is measurable  for all  $A\in \cs^{\otimes 2}$,  and $P(x,
\cdot)$ is  a probability measure on $(S^2,\cs^{\otimes 2})$ for all $x \in
S$. For any $g\in \cb_b(S^3)$ and $h\in \cb_b(S^2)$,   we set for $x\in S$:
\begin{equation}
   \label{eq:Pg}
(Pg)(x)=\int_{S^2} g(x,y,z)\; P(x,\rd y,\rd z)
\quad\text{and}\quad
(Ph)(x)=\int_{S^2} h(y,z)\; P(x,\rd y,\rd z).
\end{equation}
We define $(Pg)$ (resp. $(Ph)$), or simply $Pg$ for $g\in \cb(S^3)$
(resp. $Ph$ for $h\in \cb(S^2)$), as soon as the corresponding 
integral \reff{eq:Pg} is well defined, and we have  that $Pg$ and
$Ph$ belong to $\cb(S)$.
\medskip 

We  now introduce  some notations  related to  the regular  binary tree.
 We   set   $\T_0=\G_0=\{\emptyset\}$,
$\G_k=\{0,1\}^k$  and $\T_k  =  \bigcup _{0  \leq r  \leq  k} \G_r$  for
$k\in  \N^*$, and  $\T  =  \bigcup _{r\in  \N}  \G_r$.   The set  $\G_k$
corresponds to the  $k$-th generation, $\T_k$ to the tree  up to the $k$-th
generation, and $\T$ the complete binary  tree. For $i\in \T$, we denote
by $|i|$ the generation of $i$ ($|i|=k$  if and only if $i\in \G_k$) and
$iA=\{ij; j\in A\}$  for $A\subset \T$, where $ij$  is the concatenation
of   the  two   sequences  $i,j\in   \T$,  with   the  convention   that
$\emptyset i=i\emptyset=i$.

We recall the definition of bifurcating Markov chain  from
\cite{Guyon}. 
\begin{defi}
  We say  a stochastic process indexed  by $\T$, $X=(X_i,  i\in \T)$, is
  a bifurcating Markov chain (BMC) on a measurable space $(S, \cs)$ with
  initial probability distribution  $\nu$ on $(S, \cs)$ and probability
  kernel $\cp$ on $S\times \cs^{\otimes 2}$ if:
\begin{itemize}
\item[-] (Initial  distribution.) The  random variable  $X_\emptyset$ is
  distributed as $\nu$.
   \item[-] (Branching Markov property.) For  any sequence   $(g_i, i\in
     \T)$ of functions belonging to $\cb_b(S^3)$, we have for all $k\geq 0$,
\[
\E\Big[\prod_{i\in \G_k} g_i(X_i,X_{i0},X_{i1}) |\sigma(X_j; j\in \T_k)\Big] 
=\prod_{i\in \G_k} \cp g_i(X_{i}).
\]
\end{itemize}
\end{defi}

Let $X=(X_i,  i\in \T)$ be a BMC  on a measurable space $(S, \cs)$ with
  initial probability distribution  $\nu$ and probability
  kernel $\cp$. 
We    define     three
probability kernels $P_0, P_1$ and $\cq$ on $S\times \cs$ by:
\[
P_0(x,A)=\cp(x, A\times S), \quad
P_1(x,A)=\cp(x, S\times A) \quad\text{for
$(x,A)\in S\times  \cs$, and}\quad
\cq=\inv{2}(P_0+P_1).
\] 
Notice  that  $P_0$ (resp.   $P_1$)  is  the  restriction of  the  first
(resp. second) marginal of $\cp$ to $S$.  Following \cite{Guyon}, we
introduce an  auxiliary Markov  chain $Y=(Y_n, n\in  \N) $  on $(S,\cs)$
with  $Y_0$ distributed  as $X_\emptyset$  and transition  kernel $\cq$.
The  distribution of  $Y_n$ corresponds  to the  distribution of  $X_I$,
where $I$  is chosen independently from  $X$ and uniformly at  random in
generation  $\G_n$.    We  shall   write  $\E_x$   when  $X_\emptyset=x$
(\textit{i.e.}  the initial  distribution  $\nu$ is  the  Dirac mass  at
$x\in S$).

\begin{rem}\label{rem:ineq-simple}
By convention, for  $f,g\in  \cb(S)$,   we  define  the  function $f\otimes  g\in \cb(S^2)$  by $(f\otimes g)(x,y)=f(x)g(y)$ for  $x,y\in S$ and introduce the notations:
\[
f\sot g= \inv{2}(f\otimes g + g\otimes f) \quad\text{and}\quad f\otimes ^2= f\otimes f.
\]
Notice that  $\cp(g\sot \ind)=\cq(g)$ for  $g\in \cb_+(S)$. For $f \in \cb_+(S)$, as $ f\otimes f \leq f^2 \sot \ind $, we get:
\begin{equation}\label{eq:majo-pfxf}
\cp(f\otimes^2)=\cp(f\otimes f) \leq  \cp(f^2\sot \ind) = \cq\left(f^2\right). 
\end{equation}
\end{rem}

\begin{rem}
  \label{rem:jg}
If the Markov chain $Y$ is ergodic and if $\mu$ denotes its unique
invariant probability measure, then Guyon proves in \cite{Guyon} that,
when $S$ is a metric space,  for all $f \in \Cc_{b}(S)$,
\begin{equation*}
|\A_{n}|^{-1} \sum_{u \in \A_{n}} f(X_{u}) \underset{n \rightarrow \infty}{\xrightarrow{\hspace*{0.75cm}}} \langle \mu,f \rangle \quad \text{in probability,} \quad \text{where $\A_{n} \in \{\GG_{n}, \TT_{n}\}.$}
\end{equation*} 
One can then see that the study of BMC is strongly related to the knowledge of $\mu$. However, when it exists, the invariant probability $\mu$ is  generally not known. The aim of this article is then to estimate $\mu$ and study, under appropriate hypotheses, the fluctuations of the estimators of $\mu$. 
\end{rem} 

\noindent  We consider the following ergodic properties of $\cq$, which in particular  implies that $\mu$ is indeed the unique invariant probability measure for $\cq$.   We refer to \cite{dmps:mc} Section 22 for a detailed account on $L^2(\mu)$-ergodicity (and in particular Definition 22.2.2 on exponentially convergent Markov kernel). 

\begin{hyp}[Geometric ergodicity]
  \label{hyp:Q1}
  The Markov kernel $\cq$ has  an (unique) invariant probability measure
  $\mu$, and $\cq$ is $L^2(\mu)$ exponentially convergent, that is there
  exists  $\alpha  \in   (0,1)$  and  $M$  finite  such   that  for  all
  $f \in L^{2}(\mu)$:
\begin{equation}\label{eq:L2-erg}
\normm{\cq^{n}f - \langle \mu, f\rangle}\leq M \alpha^{n} \normm{f} \quad \text{for all $n \in \NN$}.
\end{equation}
\end{hyp}

\begin{rem}\label{rem:LmupfgR}
By Cauchy-Schwartz we have for $f,g\in L^2(\mu)$:
\begin{align}
\label{eq:Pfog2} |\cp ( f\otimes g)| ^2&\leq  \cp (f^2 \otimes 1)
                   \,\cp(1\otimes g^2) \leq 4 \cq(f^2)\, \cq(g^2), \\ 
\label{eq:mp} \langle \mu,  \cp ( f\otimes g)\rangle& \leq 2\normm{f}\, \normm{g}.
\end{align}
\end{rem}

\section{Main result}
\label{sec:main}
\subsection{Kernel estimator of the density $\mu$}
\label{sec:Noy}
The purpose of  this Section is to study asymptotic  normality of kernel
estimators for the  density of the stationary measure of  a BMC.  Assume
that  $S  =  \RR^{d}$, with $d\geq 1$,   and  that the  invariant  measure  $\mu$  of  the
transition  kernel $\cq$ exists is unique and has a  density, still  denoted by  $\mu$, with
respect to  the Lebesgue measure.   Our aim  is to estimate  the density
$\mu$ from the observation of  the population over the $n$-th generation
$\G_n$ of over  $\T_n$, that is up to generation  $n$. For that purpose,
assume   we   observe   $\XX^{n}   =  (X_u)_{u   \in   \A_{n}}$,   where
$\A_{n} \in  \{\GG_{n}, \TT_{n}\}$  {\it i.e.}  we have  $2^{n+1}-1$ (or
$2^{n}$) random variables with value  in $S$.  We consider an integrable
kernel function $K\in  \cb(S)$ such that $\int_{S} K(x)\, dx  = 1$ and a
sequence of positive bandwidths  $(h_{n},n \in \NN)$ which converges
to $0$ as  $n$ goes to infinity.   Then, we can define  the estimation of
the   density    of   $\mu$    at   $x    \in   S$    over   individuals
$\A_{n}  \in   \{\TT_{n},  \GG_{n}\}$  with  kernel   $K$  and  bandwidth
$(h_n, n\in \N)$ as:
\begin{equation}\label{eq:estim-K}
\widehat{\mu}_{\A_{n}}(x) = |\A_{n}|^{-1}h_n^{-d/2} \sum_{u \in \A_{n}} K_{h_{n}}(x - X_{u}), 
\end{equation}
where for $h>0$ the rescaled kernel function $K_h$ is given for $y\in S$ by:
\[
  K_{h}(y) = h^{-d/2}K(h^{-1}\, y).
\]
Those  statistics  are  strongly  inspired  from
\cite{masry1986, Roussas1969a, Roussas1991}.
For $h>0$ and $u\in \T$, we set:
\[
  K_{h}\star\mu(x) = \EE_{\mu}[K_{h}(x-X_{u})] = \int_{S}
K_{h}(x-y) \mu(y)\, dy.
\]
We have the following bias-variance type decomposition of the estimator
$\widehat{\mu}_{\A_{n}}(x)$: 
\begin{equation}\label{eq:DeBiVa}
  \widehat{\mu}_{\A_{n}}(x) - \mu(x) =
B_{h_n}(x) + V_{\A_n, h_n}(x),
\end{equation}
where for $h>0$ and $\A\subset \T$ finite:
\[
  B_h(x)=h^{-d/2}K_{h}\star\mu(x) - \mu(x)
  \quad\text{and}\quad
  V_{\A, h}(x)=  |\A|^{-1}h^{-d/2} \sum_{u \in \A} \Big(K_{h}(x -
  X_{u}) - K_{h}\star\mu(x)\Big).
\]
Our aim is to study the convergence and the asymptotic normality 
of the estimator $ \widehat{\mu}_{\A_{n}}(x)$ of $\mu(x)$. This relies
on a series of assumption on the model, that is on $\cp$, $\cq$ and
$\mu$, and on the kernel function $K$ as well as the bandwidth $(h_n, n\in
\N)$. \medskip

We first state a series of assumption of the density of the kernel $\cp$
and the initial distribution $\nu$ with respect to the invariant
measure.

\begin{hyp}[Regularity of  $\cp$ and $\nu_0$]
  \label{hyp:densite-cp}
 We assume that:
  \begin{propenum}
  \item\label{item:densite-cp}
    There exists  an invariant probability measure $\mu$  of $\cq$ and
  the  transition kernel  $\cp$  has  a density,  denoted  by $p$,  with
  respect    to    the    measure   $\mu^{\otimes    2}$,    that    is, for all $x\in S$:
  \[
    \cp(x, dy, dz) = p(x,y,z) \, \mu(dy)\mu(dz).
  \]

\item\label{item:int-q2-mu} The following function $\bH$ defined on $S$ belongs
  to $L^2(\mu)$, where:
 \begin{equation}\label{eq:def-bH0}
  \bH(x) = \left(\int_{S} q(x,y)^{2}\,  \mu(dy)\right)^{1/2},
\end{equation}
with $ q(x,y)=2^{-1}   \int_S   (p(x,y,z)+p(x,z,y))
\,\mu(dz)$, 
the density of $\cq$ with respect to $\mu$. 

\item\label{item:hk-L6}  There exists $k_{1} \geq 1$ such that $\bH_{k_1}\in
   L^6(\mu)$, where for $k\in \N^*$:
   \[
       \bH_{k} = \Qq^{k-1}\bH.
  \]
 \item\label{item:k00}  There exists $k_0\in \N$, such that the  probability measure $\nu \cq^{k_0}$ has a bounded density, say $\nu_0$, with respect to $\mu$:
\begin{equation*}
\nu \cq^{k_0}(dy) = \nu_0(y) \mu(d y) \quad\text{and} \quad \norm{\nu_0}_\infty <+\infty .
\end{equation*}    
  \end{propenum}
\end{hyp}
On one hand,  Conditions \ref{item:densite-cp}, \ref{item:int-q2-mu} and
\ref{item:k00}  can be  seen  as standard  $L^2$  condition for  ergodic
Markov chains.   On the other  hand, even  in the simpler  symmetric BAR
model presented in Section \ref{sec:sym-bar},  it may happens that $\bH$
has  no finite  higher  moments (which  are  used in  the  proof of  the
asymptotic  normality  to check  Lindeberg's  condition  using a  fourth
moment condition, see also  Assumption \ref{hyp:DenMu}).  This motivated
the introduction of Condition \ref{item:hk-L6}.

Then, we consider the real valued case, and assume further integrability
condition on the density of $\cp$ and $\cq$, and the existence of the
density of $\mu$ with respect to the Lebesgue measure. 

\begin{hyp}[Regularity of  $\mu$ and integrability conditions]
  \label{hyp:densite-mu2}
  Let $S=\R^d$ with $d\geq 1$. Assume that Assumption
  \ref{hyp:densite-cp}~\ref{item:densite-cp}
  holds. 
  \begin{propenum}
     
   \item\label{item:densite-mu} The invariant measure $\mu$ of
   the transition kernel $\cq$ has a density, still denoted by $\mu$,
   with respect to the Lebesgue measure.

    \item\label{item:C0-2} The  following constants  are finite:
\begin{align}\label{eq:ex-C0kernel}
C_0&=\sup_{x, y \in \R^d} \big(\mu(x)+ q(x,y) \mu(y)\big),\\
  \label{eq:ex-C1kernel}
  C_{1}&= \sup_{y,z \in \RR^{d}} \int_{\RR^{d}} dx\, \mu(x)\mu(y)\mu(z)p(x,y,z), \\
  \label{eq:ex-C2kernel}
  C_{2}&= \int_{\RR^{d}} dx\, \mu(x) \, \sup_{z\in \R^d}
         \left(\int _{\R^d} dy\, \mu(y) \bH(y) \, \mu(z) \big(p(x,
         y,z)+p(x, z,y)\big)\right)^2. 
\end{align}
    \end{propenum}
\end{hyp}

Following    \cite[Theorem~1A]{Parzen1962} (which we consider in
dimension $d$, see Lemma \ref{lem:bochner} below), we shall consider the 
following assumptions. For $g\in \cb(\R^d)$, we set $\norm{g}_p =
(\int_{S} |g(y)|^p\, dy )^{1/p} $. Then, we consider condition of the
kernel function. 

\begin{hyp}[Regularity of the kernel function and the bandwidths]
  \label{hyp:K}
  Let $S=\R^d$ with $d\geq 1$.
  \begin{propenum}
\item\label{item:cond-f}  The kernel function $K\in \cb(S)$  satisfies:
    \begin{equation}\label{eq:cond-f}
  \norm{K}_{\infty}<+\infty, \,\, 
  \norm{K}_1      < + \infty, \,\, 
\norm{K}_{2} <+\infty, \,\, 
\int_{S}\!  K(x)\,   dx  =  1
\quad\text{and}\quad
\lim_{|x|\rightarrow  +\infty}       |x|K(x)=0.
\end{equation} 
   \item\label{item:bandwidth} There exists $\gamma\in (0, 1/d)$ such
     that the    bandwidths  $(h_n,n \in \NN)$  
are  defined by $h_n= 2^{-n\gamma}$. 
  \end{propenum}
\end{hyp}

The following regularity  assumptions on $\mu$, the  kernel function $K$
and the bandwidth sequence $(h_n, n\in \N)$ will be useful to control de
biais    term   in    \eqref{eq:DeBiVa}. We follow Tsybakov \cite{tsybakov2008introduction}, chapter 1.    For    $s\in   \R_+$,    let
$\lfloor s  \rfloor$ denote its integer  part, that is the  only integer
$n\in \N$ such that $n \leq  s<n+1$ and set $\{s\}=s- \lfloor s \rfloor$
its fractional part.

\begin{hyp}[Further regularity on the density $\mu$, the kernel function and the bandwidths]
  \label{hyp:estim-tcl}
  Suppose  that there exists  an invariant probability measure $\mu$  of
  $\cq$ and that Assumptions \ref{hyp:densite-mu2}~\ref{item:densite-mu}
  and \ref{hyp:K} hold. 
We assume there exists  $s > 0$ such that the following holds:
\begin{propenum}
\item\label{item:den-Holder}  \textbf{The density  $\mu$ belongs  to the
    (isotropic) H\"older class  of order $(s, \ldots,  s) \in \RR^{d}$:}
  The  density $\mu$ admits partial derivatives with
  respect to  $x_{j}$, for all  $j\in \{1,\ldots  d\}$, up to  the order
  $\lfloor s  \rfloor$ and there exists  a finite constant $L  > 0$ such
  that  for all  $x=(x_1,  \ldots,  x_d), \in  \RR^{d}$,  $t\in \R$  and
  $ j \in \{1, \ldots, d\}$:
\begin{equation*}
\left|\frac{\partial^{\lfloor s \rfloor}\mu}{\partial x_{j}^{\lfloor s
      \rfloor}}(x_{-j},t)-\frac{\partial^{\lfloor s
      \rfloor}\mu}{\partial x_{j}^{\lfloor s \rfloor}}(x)\right| \leq
L|x_{j} - t|^{\{s\}},  
\end{equation*}
where $(x_{-j},t)$ denotes the vector $x$ where we have replaced the
$j^{th}$ coordinate $x_{j}$ by $t$, with the convention
${\partial^{0}\mu}/{\partial x_{j}^{0}} = \mu$. 

\item\label{item:K}\textbf{The kernel $K$ is  of order $(\lfloor s
    \rfloor, \ldots, \lfloor s \rfloor) \in \NN^{d}$:} We have
  $\int_{\RR^{d}} |x|^{s}K(x)\, dx < \infty$  and $\int_{\RR}
  x^{k}_{j}\, K(x)\, dx_{j} = 0$ for all $k \in \{1,\ldots,\lfloor s
  \rfloor\}$ and $j \in \{1,\ldots,d\}$. 

\item\label{item:bandew}\textbf{Bandwidth control:}
The bandwidths $(h_n,n \in \NN)$ satisfy $\lim_{n \rightarrow \infty}
|\GG_{n}|\,h_{n}^{2s + d} = 0$, that is 
$\gamma> 1/(2s+d)$.
   \end{propenum}
\end{hyp}
Notice that Assumption \ref{hyp:estim-tcl}-(i) implies that $\mu$ is at
least H\"older continuous as $s>0$.  
\medskip

First, we have the  following result which provides the consistency of the estimator $\widehat{\mu}_{\A_{n}}(x)$ for $x$ in the set of continuity of $\mu$. Its proof is given in Section \ref{sec:lem-Cv-mu}. 
\begin{lem}[Convergence of the kernel density estimator]\label{lem:cge-ps-muhat}
  Let $X$ be  a BMC with kernel $\cp$ and  initial distribution $\nu$,
  $K$ a kernel  function and $(h_n, n\in \N)$ a  bandwidth sequence such
  that   Assumptions  \ref{hyp:Q1}   (on   the  geometric   ergodicity),
  \ref{hyp:densite-cp}  (on  the regularity  of  $\cp$  and of  $\nu$),
  Assumptions \ref{hyp:densite-mu2} (on the density of $\mu$ and $\cp$),
  Assumptions \ref{hyp:K} (on the kernel function $K$ and the bandwidths
  $(h_n, n\in \N)$),
  and Assumptions 
  \ref{hyp:estim-tcl} (on the  density $\mu$, $K$ and  $(h_n, n\in \N)$)
  are in force.

  Furthermore, if  the ergodic rate of convergence 
  $\alpha$    (given   in    Assumption~\ref{hyp:Q1})   is    such   that
  $\alpha>1/\sqrt{2}$,  then assume  that  the  bandwidth rate  $\gamma$
  (given in Assumption~\ref{hyp:K}~\ref{item:bandwidth}) is such
  that:
  \begin{equation}
    \label{eq:bandwidth}
    2^{d\gamma} > 2\alpha^2. 
  \end{equation}
    Then,  for $x$ in the set of
  continuity of $\mu$ and $\A_{n} \in \{\GG_{n}, \TT_{n}\}$, we have the
  following convergence in probability:
\[
\lim_{n \rightarrow \infty} \widehat{\mu}_{\A_{n}}(x) = \mu(x).
\]
\end{lem}

We now study the asymptotic normality of the density kernel estimator. The
proof of the next theorem is given in Section \ref{sec:p-thm-MK}.

\begin{thm}[Asymptotic normality  of the kernel density estimator]
  \label{thm:MK}
  Under  the hypothesis  of  Lemma \ref{lem:cge-ps-muhat},  we have  the
  following convergence in distribution for $x$ in the set of continuity
  of $\mu$ and $\A_{n} \in \{\GG_{n}, \TT_{n}\}$:
  \begin{equation}
    \label{eq:MKsub2}
|\A_{n}|^{1/2} h_{n}^{d/2} (\widehat{\mu}_{\A_{n}}(x) - \mu(x)) \inL G,
\end{equation}
where $G$ is a centered Gaussian real-valued random variable with 
variance $\norm{K}_{2}^{2} \, \mu(x)$.  
\end{thm}

\begin{rem}\label{rem:bandwidth}
  The bandwidth  must be  a function  of the  geometric ergodic  rate of
  convergence via  the relation  $2^{d\gamma} >  2 \alpha^{2}$  given in
  Equation \eqref{eq:bandwidth}. Notice this condition is automatically
  satisfied     in     the     critical    and     sub-critical     case
  ($\alpha\leq 1/\sqrt{2}$)  as $\gamma>0$.  In the  super-critical case
  ($\alpha>  1/\sqrt{2}$), the  geometric rate  of convergence  $\alpha$
  could  be interpreted  as  a regularity  parameter  for the  bandwidth
  selection  problems  of the  estimation  of  $\mu(x)$, just  like  the
  regularity of the unknown function  $\mu$.  With this new perspective,
  we think that the results  in \cite{bitsekiroche2020} could  be extended to
   $\alpha \in (1/2,1)$ by studying an adaptive procedure
  with respect to the unknown geometric rate of convergence $\alpha$.
\end{rem}

\begin{rem}
We stress that the asymptotic variance is the same for $\A_{n} =
\GG_{n}$ and $\A_{n} = \TT_{n}.$ This is a consequence of the structure
of the asymptotic variance $\sigma^2$ in \eqref{eq:limf-ln} and
\eqref{eq:val-s2}, 
and the fact
that $\lim_{n \rightarrow \infty} |\TT_{n}|/|\GG_{n}| = 2.$  
\end{rem}

\begin{rem}
  \label{rem:indep}
  Using the structure of the asymptotic variance $\sigma^2$ in
  \eqref{eq:limf-ln} (see also Remark \ref{rem:s=+} or consider also the
  functions  $f_{\ell,n}=f_{\ell,n}^{\text{shift}}$ given by
  \eqref{eq:def-f-kernel} in the proofs of Lemma
  \ref{lem:cge-ps-muhat} and Theorem~\ref{sec:p-thm-MK}), it is easy to
  deduce that the estimators 
  $|\G_{n-\ell}|^{1/2} h_{n-\ell}^{d/2} (\widehat{\mu}_{\G_{n-\ell}}(x) - \mu(x)) $
  are asymptotically independent for $\ell\in \{0, \ldots, k\}$ for any
  $k\in \N$. 
\end{rem}

\subsection{Application to the study of symmetric BAR}\label{sec:sym-bar}
\subsubsection{The model}
We consider a particular case from \cite{CS86} of the real-valued
bifurcating autoregressive process (BAR), see also
\cite[Section~4]{BD2}. More precisely, let $a \in (-1,1).$ We consider
the process $X = (X_{u}, u \in \TT)$ on $S=\R$ where for all $u \in \TT$: 
\begin{equation*}
\begin{cases}
X_{u0} = a X_{u} + \ep_{u0},  \\
 X_{u1} = a X_{u}  + \ep_{u1},
\end{cases}
\end{equation*}
with $((\ep_{u0},\ep_{u1}),\, u \in \T)$  an independent  sequence of bivariate Gaussian  $\Nn(0,\Gamma)$ random vectors independent of $X_{\emptyset}$ with covariance matrix, with $\sigma>0$:
\begin{equation*}
\Gamma = \begin{pmatrix} \sigma^{2} \quad 0 \\ 0 \quad \sigma^{2} \end{pmatrix}.
\end{equation*}
Then the process $X=(X_{u},u\in\TT)$ is a BMC with transition probability $\Pp$ given by:
\begin{equation*}
\Pp(x,dy,dz) = \frac{1}{2\pi \sigma^{2}} \, \exp\left( -\frac{(y - a x)^{2} + (z - a x )^{2}}{2\sigma^{2}} \right)\, dydz = \Qq(x,dy) \Qq(x,dz),
\end{equation*}
where the transition kernel $\Qq$ of the auxiliary Markov  chain is defined by:
\begin{equation*}
\Qq(x,dy) = \frac{1}{\sqrt{2\pi\sigma^2}}\exp\left(-\frac{(y - a x )^{2}}{2\sigma^{2}}\right)\, dy.
\end{equation*} 
We have $\cq f(x)=\E[f(ax +\sigma G)]$ and more generally:
\begin{equation}
   \label{eq:Qn}
   \cq^n f(x)
   =\E\left[f\left(a^n x + \sqrt{1- a^{2n}}\,  \sigma_a G\right)\right],
\end{equation}
where $G$ is a standard $\cn(0, 1)$ Gaussian random variable and
$\sigma_a=\sigma (1- a^2)^{-1/2}$. 
The kernel 
$\Qq$ admits a unique invariant probability measure $\mu$, which is
$\cn(0, \sigma_a^2)$ and whose  density, still denoted by $\mu$, with respect
to the Lebesgue measure  is given by:
\begin{equation}
  \label{eq:mu-SBAR}
\mu(x) = \frac{\sqrt{1 - a^{2}}}{\sqrt{2\pi\sigma^2}} \exp\left(-\frac{(1 -
  a^{2})x^{2}}{2 \sigma^{2}}\right). 
\end{equation}
The density $p$ (resp. $q$) of the kernel $\cp$ (resp. $\cq$) with
respect to $\mu^{\otimes 2}$ (resp. $\mu$) are given by:
\begin{equation}
  \label{eq:cp-SBAR}
p(x,y,z) = q(x,y)q(x,z)
\end{equation}
and
\[
q(x,y) = \frac{1}{\sqrt{1 - a^{2}}} \exp\left(-\frac{(y-ax)^{2}}{2\sigma^{2}}
+ \frac{(1-a^{2})y^{2}}{2\sigma^{2}}\right) 
 = \frac{1}{\sqrt{1 - a^{2}}} \expp{- (a^2 y^2 + a^2
   x^2 -2axy )/ 2\sigma^ 2}.
\]
In particular, we have:
\[
\mu(x) \, q(x,y)= \frac{1}{\sqrt{2 \pi \sigma^2}} \exp\left(-\frac{(x-ay)^{2}}{2\sigma^{2}}
\right) .
\]


\subsubsection{Regularity of the model, and verification of the Assumptions}
We first check that Assumption \ref{hyp:Q1} on the geometric ergodicity holds. Since  $q$ is
symmetric, the    operator $\Qq$ (in
$L^2(\mu)$) is a symmetric integral Hilbert-Schmidt
operator. Furthermore its  eigenvalues are given
by $\sigma_{p}(\Qq) = (a^{n}, n \in \NN)$, with their algebraic
multiplicity being one. So Assumption \ref{hyp:Q1} holds with
$\alpha=|a|$ as $a\in (-1, 1)$.
\medskip 

We check  Assumption \ref{hyp:densite-cp}  on the regularity of  $\cp$ and $\nu_0$.  Condition \ref{item:densite-cp}  therein holds  thanks to
\eqref{eq:cp-SBAR}. Recall  $\bH$ defined  in \eqref{eq:def-bH0}.  It is
not difficult to check that for $x\in \R$:
\begin{equation}
   \label{eq:bH-a}
   \bH(x)=(1-a^4)^{-1/4}\exp\left(\frac{a^2(1-a^2)}{1+a^2}\,
     \frac{x^2}{2\sigma^2}\right) , 
\end{equation}
and     thus    $\bH         \in          L^2(\mu)$         (that         is
$\int_{\R^2}  q(x,y)^2  \,  \mu(x)  \mu(y)\, dx  dy<+\infty  $).
Thus  Condition
\ref{item:int-q2-mu} holds. 

We now consider Condition \ref{item:hk-L6}, that is
 $\bH_{k}= \cq^{k -1}
\bH $ belongs to $ L^6(\mu)$ for some
$k\geq 1$. 
We deduce from \eqref{eq:Qn} and \eqref{eq:bH-a} that there
exists a finite constant $C_{k}$ such that: 
\[                                               
\bH_{k}(x)= \cq^{k-1} \bH(x) = C_{k}  \exp\left(\frac{ a^{2k} x^{2}}{2 \sigma_{a
}^{2}(1 + a^{2k })}\right).
\]
So we deduce that  $\bH_{k}$ belongs to $ L^6(\mu)$ if and only if  $a^{2k}<1/5$,
which is satisfied for $k$ large enough as $a\in (-1 , 1)$. Thus,
Condition \ref{item:hk-L6} holds. 
 
\begin{rem}
   \label{rem:technical}
   As  we  shall see,  Assumption  \ref{hyp:densite-cp}~\ref{item:hk-L6}
   (the 6th moment  of $h_k$ being finite for some  $k\in \N^*$) is used
   to   check  \eqref{eq:Hil-SchP-1}   and  \eqref{eq:Hil-SchP-2}   from
   Assumption \ref{hyp:DenMu},  see Section \ref{sec:general-res}.  So one
   could  ask if  those two  inequalities could  hold without  Condition
   \ref{item:hk-L6}.   In fact,  using  elementary  computations, it  is
   possible   to    check   the    following.    For   $k_{1}    =   1$,
   \eqref{eq:Hil-SchP-1}     holds    for     $|a|<    3^{-1/4}$     and
   \eqref{eq:Hil-SchP-2}   also   holds   for   $|a|\leq   0.724$   (but
   \eqref{eq:Hil-SchP-2}  fails  for  $|a|\geq  0.725$).   (Notice  that
   $2^{-1/2}< 0.724< 3^{-1/4}$.)  For $k_{1} = 2$, \eqref{eq:Hil-SchP-1}
   holds for  $|a|< 3^{-1/6}$  and \eqref{eq:Hil-SchP-2} also  holds for
   $|a|\leq    0.794$     (but    \eqref{eq:Hil-SchP-2}     fails    for
   $|a|\geq 0.795$).  So we  see that checking \eqref{eq:Hil-SchP-1} and
   \eqref{eq:Hil-SchP-2}   is  rather   tricky.    This  motivated   the
   introduction   of  the   stronger  Condition   \ref{item:hk-L6}  from
   Assumption \ref{hyp:densite-cp}.
\end{rem}

We   now   comment   on   Condition   \ref{item:k00}   from   Assumption
\ref{hyp:densite-cp}.   Notice  that  $\nu  \cq^k$  is  the  probability
distribution of $a^k X_\emptyset +  \sigma_a \sqrt{1- a^{2k}}\, G$, with
$G$  a $\cn(0,  1)$ random  variable independent  of $X_\emptyset$.   So
Condition  \ref{item:k00}  holds  in  particular if  $\nu$  has  compact
support (with  $k_0=1$) or if  $\nu$ has a  density with respect  to the
Lebesgue  measure,   which  we   still  denote   by  $\nu$,   such  that
$\norm{\rd \nu/\rd \mu}_\infty  $ is finite (with  $k_0=0$). Notice that
if    $\nu$    is    the   Gaussian    probability    distribution    of
$\cn(m_0, \rho_0^2)$, then Condition \ref{item:k00} holds if and only if
$\rho_0< \sigma_a$ and $m_0\in \R$, or $\rho_0= \sigma_a$ and $m_0=0$.

\medskip

We   now  check   Assumptions  \ref{hyp:densite-mu2}
on the regularity of  $\mu$ and on the  integrability conditions on the density
of $\cp$ and $\cq$. Condition \ref{item:densite-mu} holds, see
\eqref{eq:mu-SBAR} for the density of $\mu$ with respect to the Lebesgue
measure. We now check that Condition \ref{item:C0-2} holds, that is the
constants $C_0$, $C_1$ and $ 
C_2$ defined in \reff{eq:ex-C0kernel}, \reff{eq:ex-C1kernel} and \reff{eq:ex-C2kernel} 
 are finite.
The fact that $C_0$ is finite is
clear. Notice that:
\[
C_1=\sup_{y,z \in \RR^{d}} \int_{\RR^{d}} dx\, 
       \mu(x)\mu(y)\mu(z)p(x,y,z)
=  \sup_{y,z \in \RR^{d}} \int_{\RR^{d}} dx\, 
       \mu(x)\mu(y)\mu(z)q(x,y)q(x,z) \leq  C_0^2.
\]
We also have, using Jensen for the second inequality (and the probability measure $\mu(y) q(x, y) \, dy$):
\begin{align*}
  C_{2}
  &= 4\int_{\RR^{d}} dx\, \mu(x) \, \sup_{z\in \R^d}
    \left(\int _{\R^d} dy\, \mu(y) \bH(y)  \, \mu(z)  q(x,y) q(x,z) \right)^2 \\
  &\leq  4 C_0^2 \int_{\RR^{d}} dx\, \mu(x) \, \left(\int _{\R^d} dy\,
    \mu(y) \bH(y) \,  q(x,y) \right)^2 \\
&\leq  4 C_0^2 \normm{ \bH}^2.
\end{align*}
So, we get that the constants $C_0$, $C_1$ and $
C_2$  are finite, and thus Condition \ref{item:C0-2} holds.
\medskip

Since  the  function  $\mu$  given in  \eqref{eq:mu-SBAR}  is  of  class
$\cc^\infty  $  with  all  its  derivative bounded,  we  get  that  the
H\"older type 
Assumption  \ref{hyp:estim-tcl}~\ref{item:den-Holder}   holds  (for  any
$s>0$).

\medskip

Many  choices of the  kernel  function, $K$,  and  of the bandwidths  parameter
$\gamma$     satisfy    Assumption     \ref{hyp:K}    and     Assumption
\ref{hyp:estim-tcl} \ref{item:K} and \ref{item:bandew}.
Eventually, as $d=1$ and $\alpha=|a|$, we get that  Equation
\eqref{eq:bandwidth} becomes  $2 ^\gamma > 2 a^2$, which  holds
\textit{a fortiori} if $2a^2\leq 1$.






\subsection{Numerical studies}\label{sec:simul}

\begin{figure}[!ht]
	\centering
	\begin{subfigure}{0.45\textwidth} 
		\includegraphics[width=\textwidth]{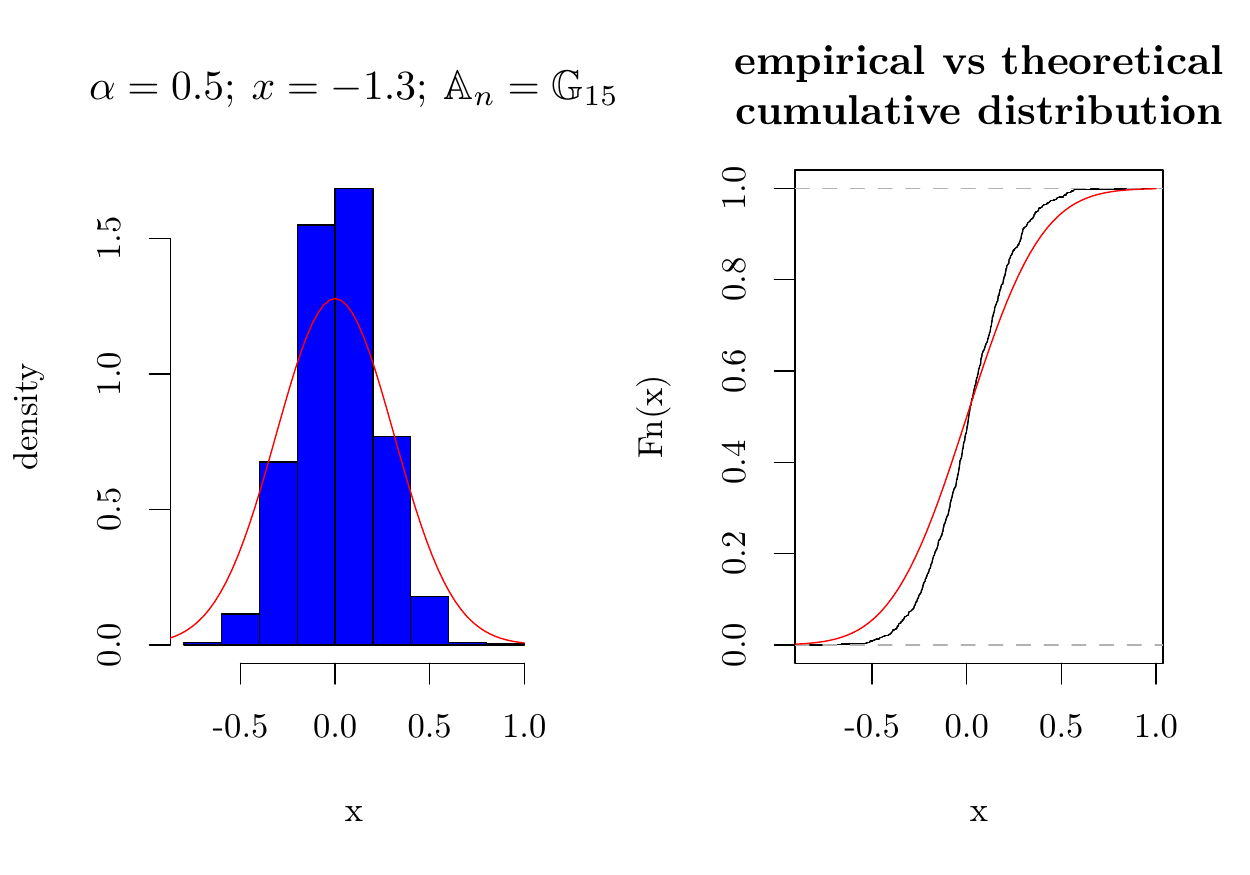}
		\caption{$\alpha=0.5$}  
	\end{subfigure}
	\begin{subfigure}{0.45\textwidth} 
		\includegraphics[width=\textwidth]{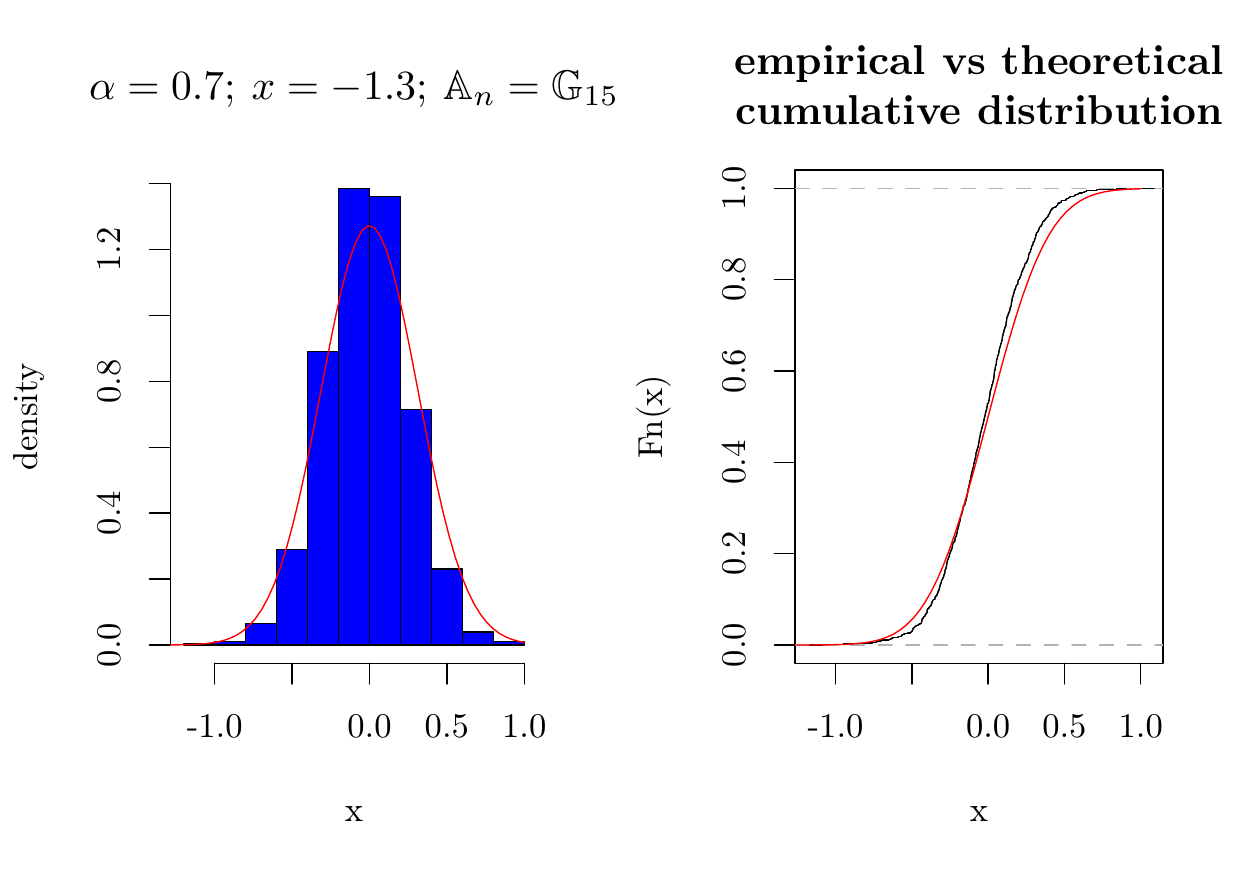}
		\caption{$\alpha=0.7$} 
	\end{subfigure}
        \caption{Histogram  and  empirical  cumulative  distribution  of
          $\zeta_{n}$  given in  \eqref{eq:def-zeta} with  $x =  - 1.3$,
          $n = 15$, $\A_{n} = \GG_{n}$  and $\gamma = 1/5 + 10^{-3}$. We
          consider  the (sub-critical) ergodic rate of convergence: $\alpha  = 0.5$
           and  $\alpha =  0.7$.} \label{fig:I-clt-noy-sub}
\end{figure}
\begin{figure}[!ht]	
	\centering
	\begin{subfigure}{0.45\textwidth} 
		\includegraphics[width=\textwidth]{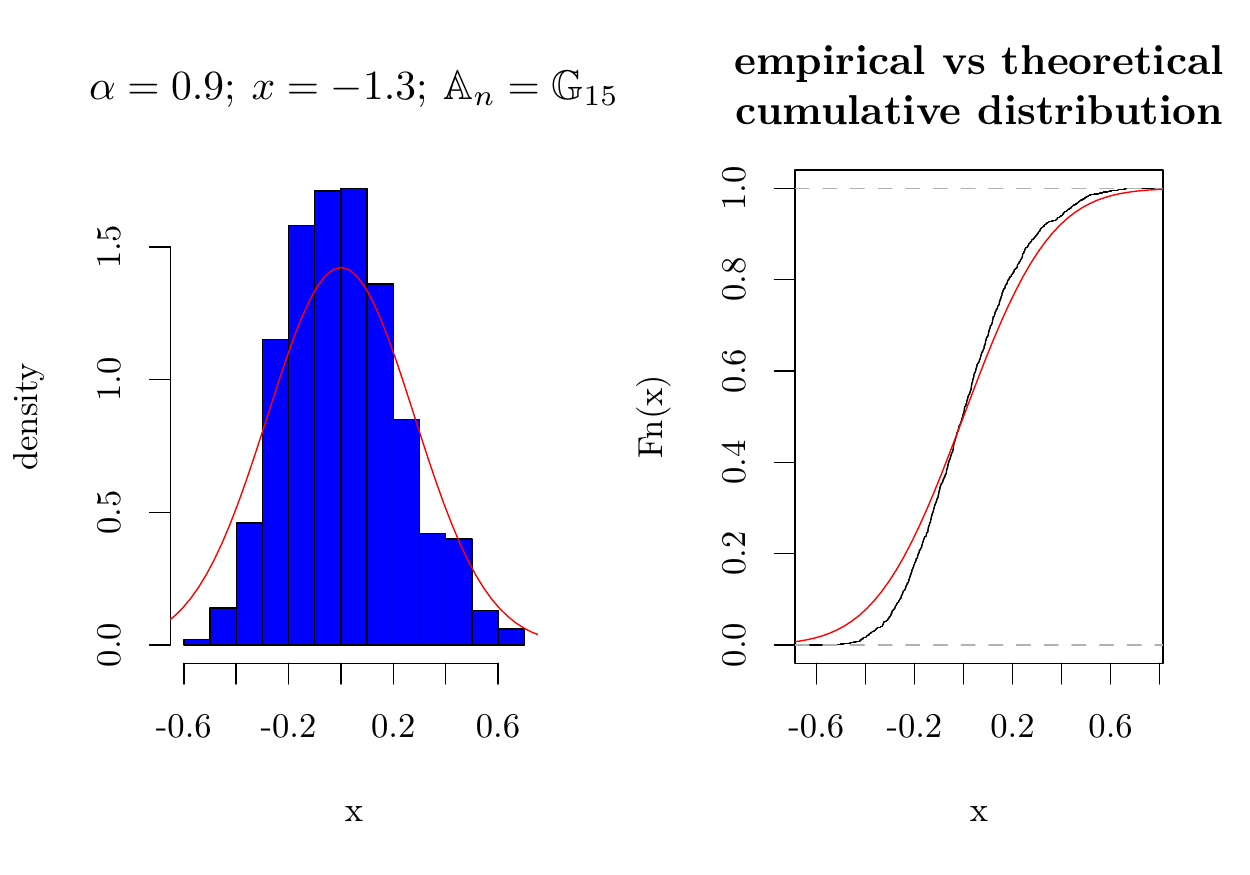}
		\caption{$\gamma = 0.696$} 
	\end{subfigure}
	\begin{subfigure}{0.45\textwidth} 
		\includegraphics[width=\textwidth]{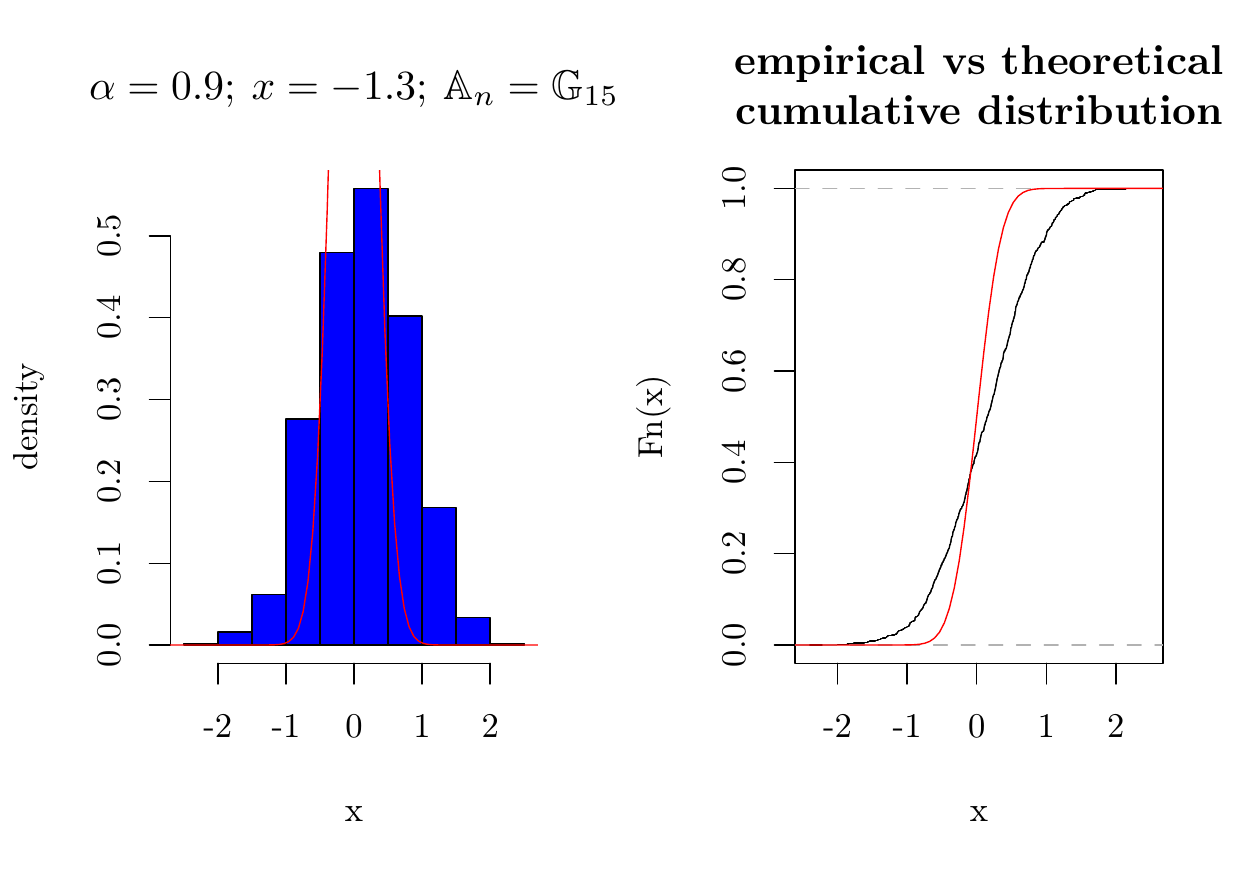}
		\caption{$\gamma = 1/5 + 10^{-3}$} 
	\end{subfigure}
	\caption{Histogram  and  empirical  cumulative  distribution  of
          $\zeta_{n}$  given in  \eqref{eq:def-zeta} with  $x =  - 1.3$,
          $n = 15$, $\A_{n} = \GG_{n}$ and the ergodic rate $\alpha=0.9$
          (super-critical  case).  We  consider  the  bandwidth  exponent
          $\gamma     =    0.696$     (which    satisfies     \eqref{eq:bandwidth}) for the two left graphics and 
          $\gamma = 1/5 + 10^{-3}$ (which does not satisfy \eqref{eq:bandwidth}) for the two right.} \label{fig:I-clt-noy-super09} 
\end{figure}

In order to illustrate the central limit theorem for the estimator of the invariant density $\mu$, we simulate $n_0=500$ samples of a symmetric BAR $X=(X_{u}^{(a)}, u \in \TT_{n})$ with different values of the autoregressive coefficient $\alpha=a\in (-1, 1)$. For each sample, we compute the estimator $\widehat{\mu}_{\A_{n}}(x)$ given in \eqref{eq:estim-K} and its fluctuation given by
\begin{equation}\label{eq:def-zeta}
\zeta_n=|\A_{n}|^{1/2} h_n^{d/2} \, (\widehat{\mu}_{\A_{n}}(x) - \mu(x))
\end{equation}
for $x\in \R$, the average over  $\A_{n} \in \{ \G_{n},  \T_{n}\}$, the Gaussian kernel
\begin{equation*}
K(x) = \frac{1}{\sqrt{2\pi}} \expp{-x^2/2}
\end{equation*}
and the bandwidth $h_n=2^{-n\gamma}$ with  $\gamma\in (0, 1)$.  Next, in order to compare theoretical and empirical  results, we plot in the same graphic, see Figures \ref{fig:I-clt-noy-sub}  and \ref{fig:I-clt-noy-super09}:
\begin{itemize}
\item The histogram of $\zeta_n$ and the density of the centered Gaussian distribution with variance $\mu(x) \norm{K}_2^2=  \mu(x)(2\sqrt{\pi})^{-1}$ (see Theorem \ref{thm:MK}). 
\item The empirical cumulative distribution  of $\zeta_n$ and the cumulative distribution of the centered Gaussian distribution with variance $\mu(x) \norm{K}_2^2=  \mu(x)(2\sqrt{\pi})^{-1}$. 
\end{itemize}
Since the Gaussian kernel is of  order $s=2$ and the dimension is $d=1$, the   bandwidth   exponent   $\gamma$   must   satisfy   the   condition $\gamma  >  1/5$,  so that  Assumption  \ref{hyp:estim-tcl}-(iii)  holds. Moreover,  in  the  super-critical   case,  $\gamma$  must  satisfy  the supplementary     condition    $2^\gamma>     2\alpha^2$, that is $\gamma   >   1   +  \log(\alpha^{2})/\log(2)$,   so  that \eqref{eq:bandwidth} holds.  In  Figure \ref{fig:I-clt-noy-sub}, we take $\alpha=0.5$  and   $\alpha=0.7$  (both  of  them   corresponds  to  the sub-critical case as  $2\alpha^2<1$) and $\gamma = 1/5  + 10^{-3}$.  The simulations  agree with  results from  Theorem \ref{thm:MK}.   In Figure \ref{fig:I-clt-noy-super09},  we  take  $\alpha =  0.9$  (super-critical case) and consider  $\gamma = 0.696$ and $\gamma = 1/5 + 10^{-3}$. In the  former   case \eqref{eq:bandwidth}  is  satisfied  as $\gamma = 0.696  > 1 + \log((0.9)^{2})/\log(2)$, and in  the latter case \eqref{eq:bandwidth}  fails. As one  can see in  the graphics Figure  \ref{fig:I-clt-noy-super09},  the   estimates   agree with the theory in the former case ($\gamma = 0.696$), whereas they are poor in the latter case.

\subsection{A general CLT for additive functionals of BMC}
\label{sec:general-res}
The proof of Lemma \ref{lem:cge-ps-muhat} and Theorem \ref{thm:MK} rely
on a general  central limit result for additive functionals of BMC. 
In the spirit of 
\cite{BD2}, we introduce  the following series of
assumptions in a general $L^2(\mu)$ framework, with increasing
 conditions as the geometric ergodic rate $\alpha$ exceed the critical
threshold of $1/\sqrt{2}$. In fact,
  we believe that the general framework presented in this section may be
  used also for others nonparametric smoothing methods for BMC than the
  one presented in Section \ref{sec:Noy}.
  \medskip

  Let $X=(X_u, u\in \T)$ be a BMC on $(S, \cs)$ with initial probability
  distribution $\nu$, and probability kernel $\cp$.  Recall $\cq$ is the
  induced Markov kernel. In the spirit  of Assumption 2.4 and Remark 2.5
  in \cite{BD2}, we consider the  following hypothesis on asymptotic and
  non-asymptotic distribution of the process.
\begin{hyp}[$L^2(\mu)$ regularity for the probability kernel $\cp$ and
  density of the initial distribution]
  \label{hyp:DenMu}
There exists an  invariant probability measure $\mu$ of $\cq$ and:
\begin{propenum}
\item\label{item:P-L2}
There exists $k_1\in \N$ and a finite constant $M$ such that for all $f, g\in
L^{2}(\mu)$:
\begin{equation}
 \label{eq:Hil-SchP-1} \normm{\cp (\cq^{k_1} f \otimes \cq^{k_1} g)} \leq
  M  \normm{f}  \normm{g},
 \end{equation}
 and for all $h\in L^{2}(\mu)$, and all $m\in \{0, \ldots, k_1\}$:
\begin{equation}
 \label{eq:Hil-SchP-2} 
\normm{\cp \left(\cq^{m} \cp(\cq^{k_{1}} f \sot \cq^{k_{1}} g)\sot
    \cq^{k_{1}} h\right)} \leq M  \normm{f}  \normm{g} \normm{h}.
 \end{equation}

\item\label{item:k0}  There exists $k_0\in \N$, such that the  probability measure $\nu \cq^{k_0}$ has a bounded density, say $\nu_0$, with respect to $\mu$:
\begin{equation*}
\nu \cq^{k_0}(dy) = \nu_0(y) \mu(d y) \quad\text{and} \quad \norm{\nu_0}_\infty <+\infty .
\end{equation*}
\end{propenum}
\end{hyp}

The next family of three assumptions are related to the sequence of
functions which will be considered.

\begin{hyp}[Regularity of the approximation functions in the
  sub-critical regime]
  \label{hyp:f}
Let  $ (f_{ \ell,n}, \, n\geq  \ell\geq 0)$  be  a sequence of
real-valued measurable functions defined on $S$  such that: 
\begin{itemize} 
\item[(i)] There exists $\rho \in (0,1/2)$ such that 
$\sup_{n\geq \ell\geq 0}2^{-n\rho}\,\|f_{\ell,n}\|_{\infty} $ is finite.
\item[(ii)] The constants $\bC_2=\sup_{n\geq \ell\geq 0}\|f_{\ell,n}\|_{L^{2}(\mu)}$ and $\bQ_2=\sup_{n\geq \ell\geq 0}
\norm{\cq (f_{\ell,n}^2)}^{1/2}_{\infty}$ are finite.
\item[(iii)] There exists a sequence $(\delta_{\ell,n}, n\geq \ell\geq 0)$ of
  positive numbers such that $\Delta=\sup_{n\geq \ell\geq 0} \delta_{\ell,n}$
  is finite, $\lim_{n\rightarrow \infty } \delta_{\ell,n}=0$ for all $\ell\in
  \N$, and for all $n\geq \ell\geq 0$:
\[
\langle \mu, |f_{\ell, n}| \rangle + 
|\langle \mu, \cp(f_{\ell, n} \otimes^2) \rangle|\leq
 \delta_{\ell, n} ;
\]
and  for all  $g\in \cb_+(S)$:
\begin{equation}\label{eq:pfqg-d}
\normm{\cp( |f_{\ell,n}|\sot  Q g)}\leq  \delta_{\ell,n} \normm{g}.
\end{equation}
\item[(iv)] The following limit exists and is finite:
\begin{equation}\label{eq:limf-ln}
\sigma^2= \lim_{n\rightarrow\infty } \sum_{\ell=0}^n 2^{-\ell} \norm{f_{\ell,n} }_{L^2(\mu)} ^2<+\infty. 
\end{equation}
\end{itemize}
\end{hyp}


\begin{rem}\label{rem:fln}
We stress that (i) and (ii) of Assumption \ref{hyp:f} imply the
existence of finite constant $C$ such that for all $n\geq \ell\geq 0$:
\begin{equation*}
\langle \mu,f_{\ell,n}^{4}\rangle \leq \|f_{\ell,n}\|^{2}_{\infty}
\langle \mu,f_{\ell,n}^{2}\rangle \leq  C\, \bC_2^2  \, 2^{2n\rho}
\quad\text{and}\quad
\langle \mu,f_{\ell,n}^{6}\rangle \leq  C\, \bC_2^2  \, 2^{4n\rho}. 
\end{equation*}
\end{rem}
We will use the following notations: for $n\in \N$, set
$\bF_{n}=(f_{\ell, n}, \ell\in \N)$ with the convention that
$f_{\ell,n}=0$ if $\ell>n$; and for $k\in \N^*$: 
\begin{equation}\label{eq:def-ck}
  c_{k}(\bF_{n})=\sup_{ \ell \geq 0}\norm{f_{\ell, n}}_{L^{k}(\mu)}
  \quad \text{and}\quad
  q_k(\bF_{n})=\sup_{ \ell \geq 0}\norm{\cq( f_{\ell,n}^k)}_{\infty }^{1/k}.
\end{equation}
In particular, we have $\bC_2=\sup_{n\geq 0} c_{2}(\bF_{n})$ and
$\bQ_2=\sup_{n\geq 0} q_{2}(\bF_{n})$.
\medskip

For the critical case, $2\alpha^{2} = 1$, we shall assume Assumption
\ref{hyp:f} as well as the following. 
\begin{hyp}[Regularity of the approximation functions in the critical
  regime]
  \label{hyp:QRjfln-crit}
Keeping the same notations as in  Assumption \ref{hyp:f}, we further 
assume that:
\begin{itemize}
   \item[(v)] 
\begin{equation}
   \label{eq:cv-nldln}
\lim_{n\rightarrow \infty } n \sum_{\ell=0}^ n 2^{-\ell/2}\,
  \delta_{\ell, n} =0. 
 \end{equation}
\item[(vi)]
 For all $n \geq \ell\geq 0$:
\begin{equation}\label{eq:bQfln-crit}
\|\Qq(|f_{\ell,n}|)\|_{\infty} \leq  \delta_{\ell,n}.
\end{equation}
\end{itemize}
\end{hyp}
For  the  super-critical  case,  $2\alpha^{2}  >  1$,  we  shall  assume
Assumptions \ref{hyp:f}, \ref{hyp:QRjfln-crit} as well as the following.

\begin{hyp}[Regularity of the approximation functions in the
  super-critical regime]
  \label{hyp:fln-Scrit}
Keeping the same notations as in  Assumption \ref{hyp:QRjfln-crit}, we further 
assume that Assumption \ref{hyp:Q1} holds with   $2\alpha^2>1$ and that:
\begin{equation}\label{eq:Cfln-Scrit}
\sup_{0\leq \ell\leq  n} (2\alpha^{2})^{n - \ell}
\delta_{\ell,n}^{2}<+\infty 
\quad\text{and, for all $\ell\in \N$,}\quad
\lim_{n \rightarrow \infty} (2\alpha^{2})^{n - \ell} \delta_{\ell,n}^{2}
= 0  .
\end{equation}
\end{hyp}
Notice that condition \reff{eq:Cfln-Scrit} implies \reff{eq:cv-nldln} as
well as $\Delta<+\infty $ and $\lim_{n\rightarrow \infty } \delta_{\ell,n}=0$ for all $\ell\in
  \N$ (see Assumption \ref{hyp:f} (iii)) when $2\alpha^2>1$. 
  \medskip
  
  Following  \cite{BD2}, for a finite set $\A\subset \T$ and a function $f\in \cb(S)$, we set:
\begin{equation}
   \label{eq:def-MA}
M_\A(f)=\sum_{i\in \A} f(X_i).
\end{equation}
We shall be interested in the cases $\A=\G_n$ (the  $n$-th generation)
and $\A=\T_n$ (the tree up to the $n$-th generation). We shall assume
that $\mu$ is an invariant probability measure of $\cq$. In view of
Remark \ref{rem:jg}, one is interested in the fluctuations of
$|\G_n|^{-1} \, M_{\G_n}(f)$ around $ \langle \mu, f \rangle$. So, we will use frequently the following notation: 
\begin{equation}
   \label{eq:def-tf}
\boxed{\tilde f= f - \langle \mu, f \rangle\quad \text{for $f\in L^1(\mu)$.}}
\end{equation}
Let $\bF=(f_\ell, \ell\in \N)$ be a sequence of elements of
$L^1(\mu)$. We set for $n\in \N$:
\begin{equation}
  \label{eq:def-NOf}
\boxed{N_{n, \emptyset}(\bF)= |\G_n|^{-1/2 }\sum_{\ell=0}^{n}M_{\G_{n-\ell}}(\tilde f_{\ell}).} 
\end{equation}
The notation $N_{n, \emptyset}$ means  that we consider the average from
the root $\emptyset$ up to the $n$-th generation.

\begin{rem}\label{rem:simpleN0n}
The following two simple cases are frequently used in the literature. Let  $f\in  L^1(\mu)$  and consider the sequence $\bF=(f_\ell,  \ell\in \N)$.   If  $f_0=f$ and $f_\ell=0$ for  $\ell \in \N^*$,  then we get:
\begin{equation*}
N_{n, \emptyset}(\bF)= |\G_n|^{-1/2} M_{\G_n}(\tilde f).
\end{equation*}
If $f_\ell=f$ for $\ell \in \N$, then we get, as  $|\T_n|=2^{n+1} - 1 $ and $|\G_n|= 2^n$:
\begin{equation*}
N_{n, \emptyset}(\bF)= |\G_n|^{-1/2} M_{\T_n}(\tilde f)
=  \sqrt{2 - 2^{-n}}\,\, |\T_n|^{-1/2} M_{\T_n}(\tilde f).
\end{equation*}
Thus, we will easily deduce the fluctuations of $M_{\T_n}(f)$ and $M_{\G_n}(f)$ from the asymptotics of $N_{n, \emptyset}(\bF)$. 
\end{rem}

The    main  result  of   this  section  is  motivated  by  the
decomposition given in \eqref{eq:DeBiVa}. It  will allow us to treat the
variance term  of kernel  estimators defined in  \eqref{eq:estim-K}. The
proof is  given in  Section \ref{sec:proof-sub-L2} for  the sub-critical
case ($\alpha\in  (0, 1/\sqrt{2})$), in  Section \ref{sec:proof-crit-L2}
for   the   critical   case   ($\alpha=1/\sqrt{2}$)   and   in   Section
\ref{sec:proof-Scrit-L2}      for       the      supercritical      case
($\alpha\in  (  1/\sqrt{2}, 1)$),  with  $\alpha$  the rate  defined  in
Assumption \ref{hyp:Q1}.
  Recall $N_{n, \emptyset}(\bF)$ defined in \eqref{eq:def-NOf}. 

  \begin{thm}\label{thm:flx}
  Let $X$ be  a BMC with kernel $\cp$ and  initial distribution $\nu$,
 such
  that   Assumption  \ref{hyp:Q1}   (on   the  geometric   ergodic rate
  $\alpha\in (0, 1)$), Assumption 
    \ref{hyp:DenMu}  (on  the regularity  of  $\cp$  and of  $\nu$)
and Assumption \ref{hyp:f} (on the approximation functions $(f_{\ell,n}, n\geq
\ell\geq 0)$)   are in force. 

Furthermore,    if    $\alpha=     1/\sqrt{2}$    then    assume    that
Assumption~\ref{hyp:QRjfln-crit} holds; and if $\alpha> 1/\sqrt{2}$ then
assume  that  Assumption~\ref{hyp:QRjfln-crit}  and Assumption  \ref{hyp:fln-Scrit}
hold.
Then, we  have the  following convergence
in distribution:  
\begin{equation*}
N_{n, \emptyset}(\bF_n) \; \xrightarrow[n\rightarrow \infty ]{\text{(d)}} \; G,
\end{equation*}
where  $\bF_n=(f_{\ell, n},  \,  \ell\in \N)$  and  the convention  that
$f_{\ell,n}=0 $ for $\ell>n$, and with $G$ a centered Gaussian random
variable with finite variance $\sigma^2$ defined in \reff{eq:limf-ln}.
\end{thm}

\begin{rem}
  \label{rem:s=+}
  Assume
  $\sigma^2_\ell=\lim_{n\rightarrow  \infty   }  \normm{f_{\ell,  n}}^2$
  exists  for  all  $\ell  \in   \N$;  so  that  $\sigma^2$  defined  in
  \reff{eq:limf-ln}          is          also          equal          to
  $\sum_{\ell\in  \N}2^{-\ell}  \sigma^2_\ell$.  According  to  additive
  form  of the  variance $\sigma^2$,  we deduce    that for
  fixed        $k\in       \N$,        the       random        variables
  $\left( |\G_n|^{-1/2  }M_{\G_{n-\ell}}(\tilde f_{\ell,n}),  \, \ell\in
    \{0, \ldots, k\} \right)$ converges  in distribution, as $n$ goes to
  infinity towards  $(G_\ell, \, \ell\in  \{0, \ldots, k\} )$  which are
  independent   real-valued   Gaussian     centered    random
  variables     with variance 
  $\Var(G_\ell)=2^{-\ell} \sigma^2_\ell$.
\end{rem}

\section{Proof of Lemma \ref{lem:cge-ps-muhat} and Theorem \ref{thm:MK}}\label{sec:p-l-th-MK}

\subsection{Checking Assumptions \ref{hyp:DenMu}, \ref{hyp:f},
  \ref{hyp:QRjfln-crit} and \ref{hyp:fln-Scrit}}
\label{sec:hyp=hyp}
We shall check that Assumptions \ref{hyp:densite-cp},
\ref{hyp:densite-mu2} and \ref{hyp:K}, and Equation
\eqref{eq:bandwidth}, for the density estimation, 
implies the more general Assumptions \ref{hyp:DenMu}, \ref{hyp:f},
  \ref{hyp:QRjfln-crit} and \ref{hyp:fln-Scrit}. 

  \medskip
We   check  that  Assumption \ref{hyp:densite-cp} implies Assumption \ref{hyp:DenMu}  on  the  $L^2(\mu)$
regularity for the  probability kernel $\cp$ and density  of the initial
distribution. Notice that  Assumption
\ref{hyp:densite-cp}~\ref{item:k00} and 
Assumption \ref{hyp:DenMu}~\ref{item:k0} coincide. 
So, it is enough  to check that  Assumption
\ref{hyp:densite-cp}~\ref{item:densite-cp}-\ref{item:hk-L6}  implies
Assumption \ref{hyp:DenMu}~\ref{item:P-L2}. 
Since $|\cq f|\leq \normm{f} \bH $, we deduce that 
$|\cq^{k_{1}}  f| \leq  \normm{f}  \bH_{k_{1}}$. 
We deduce that:
\[
 \normm{\cp (\cq^{k_1} f \otimes \cq^{k_1} g)} \leq
 \normm{f}  \normm{g} \normm{\cp (\bH_{k_{1}}\otimes^2)}.
\]
Then use \eqref{eq:majo-pfxf} to get that $\normm{\cp
  (\bH_{k_{1}}\otimes^2)}\leq  \normm{\cq\left(\bH_{k_{1}}^2\right)}
\leq   \normm{\bH_{k_{1}}^2}\leq  \norm{\bH_{k_{1}}}^2_{L^6(\mu)}<+\infty
$. This gives
\eqref{eq:Hil-SchP-1}. 
Similarly, we have:
\begin{multline*}
\normm{\cp \left(\cq^{m} \cp(\cq^{k_{1}} f \sot \cq^{k_{1}} g)\sot
    \cq^{k_{1}} h\right)} \\
\leq   \normm{f}  \normm{g} \normm{h}
   \normm{\cp \left(\cq^{m} \cp(\bH_{k_{1}}  \otimes^2 )\sot
    \bH_{k_{1}} \right)}.
\end{multline*}
On the other hand, using  \eqref{eq:Pfog2}, the  H\"older inequality and
\eqref{eq:majo-pfxf}, we also have: 
\begin{align*}
   \normm{\cp \left(\cq^{m} \cp(f_0  \otimes^2 )\sot
       g_0\right)}^2
&\leq 4 \langle \mu, \cq (  (\cq^{m} \cp(f_0\otimes^2))^2) \cq(g_0^2)\rangle\\
&\leq 4 \langle \mu, \cp(f_0\otimes^2)^{3}\rangle^{2/3}
    \langle \mu, \cq (  g_0^2)^{3}\rangle^{1/3}\\
& \leq 4 \langle \mu, f_0^6\rangle^{2/3}
    \langle \mu, g_0^6\rangle^{1/3}. 
\end{align*} 
Taking $f_0=g_0=\bH_{k_1}$ gives that $  \normm{\cp \left(\cq^{m}
    \cp(\bH_{k_{1}}  \otimes^2 )\sot 
      \bH_{k_{1}} \right)} \leq  2 \norm{\bH_{k_1}}_{L^6}^3<+\infty$. 
This gives \eqref{eq:Hil-SchP-2}. Thus,
Assumption~\ref{hyp:DenMu}~\ref{item:P-L2} holds.
\medskip

We suppose that $S=\R^d$ and that Assumptions \ref{hyp:densite-cp},
\ref{hyp:densite-mu2} hold. 
Let $K$ be a kernel function satisfying Assumption
\ref{hyp:K}~\ref{item:cond-f} and  bandwidths $(h_n, n\in \N)$ 
satisfying Assumption
\ref{hyp:K}~\ref{item:bandwidth}.
For  $x \in S$, we define the sequences of functions  $(f^x_{\ell},\, \ell\in \N)$ given by:
\begin{equation*}\label{eq:flxsub}
 f^x_\ell (y)=K_{h_\ell}(x-y)={h_\ell^{-d/2}}
 K\left(\frac{x-y}{h_\ell}\right) \quad\text{for $y\in \R^d$}.  
\end{equation*}
Then,      we     consider      the      sequences     of      functions
$(f_{\ell,n}^\text{shift},\,         n\geq         \ell\geq         0)$,
$(f_{\ell,n}^{\text{id}},\,    n     \geq    \ell    \geq     0)$    and
$(f^0_{\ell,n},\, n\geq \ell\geq 0)$ defined by:
\begin{equation}\label{eq:def-f-kernel}
  f_{\ell, n}^\text{shift} = f^x_{n-\ell},
  \quad f_{\ell,n}^{\text{id}} = f^{x}_{n}
  \quad \text{and} \quad
  f^{0}_{\ell, n} = f^x_n \ind_{\{\ell=0\}} . 
\end{equation}
\medskip

Under  those hypothesis,  we shall  check that  Assumptions \ref{hyp:f},
\ref{hyp:QRjfln-crit} and \ref{hyp:fln-Scrit} hold for those three sequences
of functions.  We first check that
(i-iii)  from   Assumption  \ref{hyp:f}   and  (v-vi)   from  Assumption
\ref{hyp:QRjfln-crit}.      We     consider    only     the     sequence
$(f_{\ell,n},\,  n\geq \ell\geq  0)$ with $f_{\ell, n}=f_{\ell,
  n}^\text{shift}$, the  arguments for  the other  two 
being                 similar.                   We                 have
$\norm{f_{\ell,n}}_{\infty}  =   h_{n-\ell}^{-d/2}  \norm{K}_{\infty}  =
2^{(n-\ell)d\gamma/2}   \norm{K}_{\infty}$.    Thus  property   (i)   of
Assumption \ref{hyp:f} holds with $\rho = d\gamma/2$.  We have:
\[
  \norm{f_{\ell, n}}_1= h_{n-\ell}^{d/2} \norm{K}_1=2^{-(n-\ell)
    d\gamma/2} \norm{K}_1
  \quad\text{and}\quad
  \norm{f_{\ell, n}}_2=\norm{K}_2.
\]
This gives $
\normm{f_{\ell,n}} ^2 \leq  \norm{f_{\ell, n}}^2_2 \, \norm{\mu}_\infty
\leq  C_0 \norm{K}_2^2$ and
\[
\norm{\cq (f_{\ell, n}^2)}_\infty  \leq  \norm{f_{\ell, n}}^2_2 \,
\sup_{x,y\in \R^d} q(x, y) \mu(y) \leq  C_0 \norm{K}_2^2.
\]
We conclude that (ii) of Assumption \ref{hyp:f} holds with $\bC_{2} =
\bQ_{2} = C_{0}^{1/2}\norm{K}_2$.
We have  $\langle \mu,|f_{\ell,n}| \rangle \leq
C_{0} \norm{K}_{1} h_{n-\ell}^{d/2}$ and $|\langle \mu, \cp(f_{\ell, n}
\otimes^2) \rangle| \leq C_{1} \norm{K}_{1}^{2}
h_{n-\ell}^d$. Furthermore, for all  $g\in \cb_+(\RR^{d})$, we have
$\normm{\cp( |f_{\ell,n}|\sot  Q g)}\leq  C_{2} h_{n-\ell}^{d/2}
\norm{K}_{1} \normm{g}$. We also have  $\norm{\cq (|f_{\ell,
    n}|)}_\infty  \leq  C_0 \norm{K}_{1} h_{n-\ell}^{d/2}$. This implies
that (iii) of Assumption \ref{hyp:f} and (vi) of Assumption
\ref{hyp:QRjfln-crit} hold  with $\delta_{\ell,n} = c \,
h_{n-\ell}^{d/2}=c 2^{- (n-\ell) d\gamma  /2}$ for some finite constant
$c$ depending only on $C_0, C_1, C_2$ and $\norm{K}_1$.
With this choice
of $\delta_{\ell, n}$, notice that (v) of Assumption
\ref{hyp:QRjfln-crit} also holds as $d\gamma<1$.  
\medskip

Recall    that    $d\gamma<1$.    Moreover,   if Equation
\eqref{eq:bandwidth} holds,, that is 
$2^{d\gamma}  >  2  \alpha^2$  where  $\alpha$ is  the  rate  given  in
Assumption \ref{hyp:Q1}  (this is  restrictive on  $\gamma$ only  in the
super-critical      regime      $2\alpha^2>1$),     then      Assumption
\ref{hyp:fln-Scrit}   also    holds   with   the   latter    choice   of
$\delta_{\ell,n}.$ \medskip

Eventually we prove (iv) of Assumption \ref{hyp:f}. 
 We recall  the following result due to Bochner
(see \cite[Theorem 1A]{Parzen1962}  which can be easily
extended to any dimension $d\geq 1$). 
\begin{lem}
  \label{lem:bochner}
Let $(h_{n},n\in\NN)$ be a sequence of positive numbers converging to $0$ as $n$ goes to infinity. Let $g: \RR^{d} \rightarrow \RR$ be a measurable function such that $\int_{\RR^{d}} |g(x)|dx < +\infty$. Let $f: \RR^{d} \rightarrow \RR$  be a measurable function such that  $\norm{f}_{\infty}<+\infty $, $\int_{\R^d} |f(y)|\, dy < + \infty$  and $\lim_{|x|\rightarrow +\infty} |x|f(x)=0$. Define
\begin{equation*}
g_{n}(x) = h_{n}^{-d}\int_{\RR^{d}} f(h_{n}^{-1}(x-y))g(y)dy.
\end{equation*}
Then, we have at every point $x$ of continuity of $g$,
\begin{equation*}
\lim_{n\rightarrow +\infty} g_{n}(x) = g(x)\int_{\RR} f(y)dy.
\end{equation*} 
\end{lem}

Let $x$ be in the set of continuity of $\mu$.  Thanks to Lemma
\ref{lem:bochner}, we have:
\begin{equation*}
\lim_{\ell\rightarrow \infty }  \norm{ f^x_{\ell}}_{L^2(\mu)}^2 = \lim_{\ell\rightarrow \infty }  \langle \mu, (f^x_{\ell})^2  \rangle = \mu(x)\norm{K}_2^2. 
\end{equation*}
We deduce  that the sequences of functions $(f_{\ell,n}^\text{shift},\, n\geq  \ell\geq 0)$, $(f_{\ell,n}^{\text{id}},\, n \geq \ell \geq 0)$ and $(f^0_{\ell,n},\, n\geq  \ell\geq 0)$ satisfy (iv) of Assumption \ref{hyp:f} with $\sigma^2$ defined by \reff{eq:limf-ln} respectively given by:
\begin{equation}\label{eq:val-s2}
(\sigma^\text{shift})^2= 2 \mu(x)\norm{K}_2^2, \, \quad (\sigma^{\text{id}})^2= 2
\mu(x)\norm{K}_2^2 \,
\quad \text{and} \quad
(\sigma^0)^2=\mu(x)\norm{K}_2^2.  
\end{equation}

\subsection{Proof of Lemma \ref{lem:cge-ps-muhat}}\label{sec:lem-Cv-mu} 
We begin the proof with $\A_{n} = \TT_{n}.$ We have the following decomposition:
\begin{equation}\label{eq:D-muhat}
\widehat{\mu}_{\TT_{n}}(x) - \mu(x) = \frac{\sqrt{|\GG_{n}|}}{|\TT_{n}|
  h_{n}^{d/2}} N_{n,\emptyset}(\bF_n) + B_{h_n}(x),  
\end{equation}
where $\bF_n=(f_{\ell, n}, \, \ell\in \N)$ with the functions
$f_{\ell,n}=f_{\ell,n}^{\text{id}}$ defined in \eqref{eq:def-f-kernel}
for $n\geq  \ell\geq 0$ and $f_{\ell,n}=0$ otherwise;
$N_{n,\emptyset}$ is defined in \eqref{eq:def-NOf} with $\bF$ replaced
by $\bF_n$; and  the bias term: 
\begin{equation*}
B_{h_n}(x) = \frac{1}{|\TT_{n}| h_{n}^{d/2}}  \sum_{\ell = 0}^{n} 2^{n -
  \ell} \langle \mu, f_{\ell,n} \rangle   - \mu(x) = \langle
\mu,h_{n}^{-d}K(h_n^{-1} (x-\cdot)) \rangle - \mu(x).  
\end{equation*}
Thanks to Section \ref{sec:hyp=hyp}, we have under the assumption of
Lemma  \ref{lem:cge-ps-muhat} that   Assumptions \ref{hyp:DenMu}, \ref{hyp:f},
  \ref{hyp:QRjfln-crit} and \ref{hyp:fln-Scrit} hold. Since $\lim_{n \rightarrow \infty} |\GG_{n}|h_{n}^{d} = \infty$ as
$\gamma<1$, we get that   $\lim_{n \rightarrow \infty} |\GG_{n}|^{1/2}/
|\TT_{n}| h_{n}^{d/2} = 0$. Thus, we get, as a direct consequence of Theorem
\ref{thm:flx} the following convergence in probability:
\begin{equation*}
\lim_{n \rightarrow \infty} \frac{\sqrt{|\GG_{n}|}}{|\TT_{n}| h_{n}^{d/2}} N_{n,\emptyset}(\bF_{n})  = 0.
\end{equation*}
Next, it follows from Lemma \ref{lem:bochner} 
that $\lim_{n \rightarrow \infty} B_{h_n}(x) =0$. By considering the
functions $f_{\ell,n}=f_{\ell,n}^{0}$ defined in \eqref{eq:def-f-kernel}, we similarly get the result for the case $\A_{n} = \GG_{n}$. 

\subsection{Proof of Theorem \ref{thm:MK}}
\label{sec:p-thm-MK}
\noindent  \textbf{The sub-critical  case and  $\A_{n} =  \TT_{n}$}.  We
keep  notations  from the  proof  of  Lemma
\ref{lem:cge-ps-muhat}. Recall that   $\bF_n=(f_{\ell, n},  \, \ell\in  \N)$  with the functions
$f_{\ell,n}=f_{\ell,n}^{\text{id}}$ defined in \eqref{eq:def-f-kernel}. Using the value of
$\sigma= \sigma^{\text{id}}$ in \eqref{eq:val-s2}, thanks to Theorem \ref{thm:flx} and the
decomposition \eqref{eq:D-muhat}, we see that to get the asymptotic normality of the
estimator \eqref{eq:MKsub2} it suffices to prove that:
\begin{equation}\label{eq:RTCLmusub}
  \lim_{n \rightarrow \infty}  |\TT_{n}|^{1/2} h_{n}^{d/2} B_{h_n}(x)
  = 0. 
\end{equation}
Using that 
\begin{multline*}
\mu(x - h_{n}y) - \mu(x) = \sum_{j = 1}^{d} (\mu(x_{1} - h_{n}y_{1},
  \ldots,
x_{j} - h_{n}y_{j}, x_{j+1}, \ldots, x_{d}) \\
- \mu(x_{1}-h_{n}y_{1}, \ldots, x_{j-1} - h_{n}y_{j-1}, x_{j},
    x_{j+1}, \ldots, x_{d})),  
\end{multline*}
the Taylor expansion and  Assumption \ref{hyp:estim-tcl}, 
we get that, for some finite constant $C > 0$,
\begin{align*}
  |\TT_{n}|^{1/2}h_{n}^{d/2} B_{h_n}(x)
  &= \sqrt{|\TT_{n}|h_{n}^{d}} \,\, \Big|\int_{\RR^{d}} h_{n}^{-d}
    K(h_n^{-1}(x-y)) \mu(y)dy - \mu(x)\Big| \\
  &= \sqrt{|\TT_{n}|h_{n}^{d}} \,\,  \Big|\int_{\RR^{d}}
    K(y)(\mu(x - h_{n}y) - \mu(x)) \, dy\Big| \\
  &\leq C \sqrt{|\TT_{n}|h_{n}^{d}}\,\, \sum_{j = 1}^{d} \,\,
    \int_{\RR^{d}} K(y)\frac{(h_{n}|y_{j}|)^{s}}{\lfloor s \rfloor!} dy  
\\
&\leq C \sqrt{|\TT_{n}|h_{n}^{2s + d}}. 
\end{align*}
Then Equation  \eqref{eq:RTCLmusub} follows,  since $\lim_{n \rightarrow \infty} |\GG_{n}|s_{n}^{2s + d} = 0$. This ends the proof  for $\A_{n} = \TT_{n}$.
\medskip

\noindent   \textbf{The sub-critical case and  $\A_{n} = \G_{n}$}. 
The proof is similar, using instead  the
functions  $f_{\ell,n}=f_{\ell,n}^{0}$ defined in
\eqref{eq:def-f-kernel}. 
\medskip
 
\noindent   \textbf{The critical and super-critical cases}. 
The proof follows the same lines,  using Theorem \ref{thm:flx} in the
critical and super-critical cases
and the decomposition \eqref{eq:D-muhat}. 

\section{Proof  of   Theorem  \ref{thm:flx} in the sub-critical case ($2\alpha^2<1$)}\label{sec:proof-sub-L2}
Recall the definition of $M_\A$ given in \eqref{eq:def-MA} and of
$\tilde f= f - \langle \mu, f \rangle$ in \eqref{eq:def-tf}. 
In order  to study the  asymptotics of $M_{\GG_{n-\ell }}(\tilde  f)$ as
$n$ goes to infinity and $\ell$  is fixed, it is convenient to consider the
contribution of the descendants of the individual $i\in \T_{n-\ell}$ for
$n\geq \ell\geq 0$:
\begin{equation}\label{eq:def-Nnil}
N^\ell_{n,i}(f)=|\G_n|^{-1/2} M_{i\G_{n-|i|-\ell}}(\tilde f), 
\end{equation}
where  $i\G_{n-|i|-\ell}=\{ij, \, j\in \G_{n-|i|-\ell}\}\subset \G_{n-\ell}$. For  all $k\in \N$ such that $n\geq k+\ell$, we have:
\begin{equation*}
M_{\G_{n-\ell}}(\tilde f)=\sqrt{|\G_n|}\,\, \sum_{i\in \G_k} N^\ell_{n,i}(f) = \sqrt{|\G_n|}\, \,N_{n, \emptyset}^\ell(f).
\end{equation*}
Let $\bF=(f_\ell, \ell\in \N)$ be a sequence of elements of
$L^1(\mu)$. 
We set for $n\in \N$ and $i\in \T_n$:
\begin{equation}\label{eq:def-NiF}
N_{n,i}(\bF)=\sum_{\ell=0}^{n-|i|} N_{n,i}^\ell(f_\ell)  = |\G_n|^{-1/2 }\sum_{\ell=0}^{n-|i|} M_{i\G_{n-|i|-\ell}}(\tilde f_\ell).
\end{equation}
We deduce that $ \sum_{i\in \G_k} N_{n,i}(\bF)=|\G_n|^{-1/2
}\sum_{\ell=0}^{n-k} M_{\G_{n-\ell}}(\tilde f_\ell)$. For $k=0$, we
recover Equation \eqref{eq:def-NOf}.

\medskip

We  consider  the  notations  of  Theorem  \ref{thm:flx}.   Recall  that
$\bF_n=(f_{\ell,  n},   \,  \ell\in   \N)$  with  the   convention  that
$f_{\ell,n}=0 $ for $\ell>n$.
In the following proofs, we will  denote by  $C$ any unimportant finite  constant which may  vary from  line to line  (in particular $C$ does not  depend on $n$ nor on $\bF_{n}$). 

\begin{rem}\label{rem:Nn=Nnk0}
Recall $k_{0}$ given in Assumption \ref{hyp:DenMu} (iii). Recall that from Assumption \ref{hyp:f} (ii), the sequence $\bF_{n}$ is bounded in $L^{2}(\mu)$.  We have
\begin{equation}\label{eq:Nn0=Nn0k0+R}
N_{n,\emptyset}(\bF_{n}) = N_{n,\emptyset}^{[k_{0}]}(\bF_{n}) \, + \, |\GG_{n}|^{-1/2} \sum_{\ell = 0}^{k_{0} - 1} M_{\GG_{\ell}}(\tilde{f}_{n-\ell,n}),
\end{equation}
where we set:
\begin{equation*}\label{eq:def-NOf2}
\boxed{ N_{n, \emptyset}^{[k_0]}(\bF_{n}) = |\G_n|^{-1/2}\sum_{\ell=0}^{n-k_0} M_{\G_{n-\ell}}(\tilde f_{\ell,n}) .} 
\end{equation*}
Using the Cauchy-Schwartz inequality, we get
\begin{equation}\label{eq:reste}
|\GG_{n}|^{-1/2} |\sum_{\ell = 0}^{k_{0} - 1}
M_{\GG_{\ell}}(\tilde{f}_{n-\ell,n})| \leq C c_{2}(\bF) |\GG_{n}|^{-1/2}
\, + \, |\GG_{n}|^{-1/2} \sum_{\ell = 0}^{k_{0}-1}
M_{\GG_{\ell}}(|f_{n-\ell,n}|).
\end{equation}   
Since the sequence $\bF_{n}$ is bounded in $L^{2}(\mu)$ and since $k_{0}$ is finite, we have, for all $\ell \in \{0, \ldots, k_{0}-1\}$, $\lim_{n \rightarrow \infty} |\GG_{n}|^{-1/2} M_{\GG_{\ell}}(|f_{n-\ell,n}|) = 0$ a.s. and then that (used \eqref{eq:reste})
\begin{equation*}
\lim_{n \rightarrow \infty} |\GG_{n}|^{-1/2} |\sum_{\ell = 0}^{k_{0} - 1} M_{\GG_{\ell}}(\tilde{f}_{n-\ell})| = 0 \quad \text{a.s.}
\end{equation*}
Therefore, from \eqref{eq:Nn0=Nn0k0+R}, the study of $N_{n,\emptyset}(\bF_{n})$ is reduced to that of $N^{[k_{0}]}_{n,\emptyset}(\bF_{n})$.   
\end{rem}

Let $(p_n, n\in \N)$ be a non-decreasing sequence of elements of $\N^*$ such that, for all $\lambda>0$:
\begin{equation}\label{eq:def-pn}
p_n< n, \quad \lim_{n\rightarrow \infty } p_n/n=1 \quad\text{and}\quad \lim_{n\rightarrow \infty } n-p_n - \lambda \log(n)=+\infty .
\end{equation}
When there is no ambiguity, we write $p$ for $p_n$. 
\medskip

Let $i,j\in \T$. We write $i\preccurlyeq  j$ if $j\in i\T$. We denote by $i\wedge j$  the most recent  common ancestor of  $i$ and $j$,  which is defined  as   the  only   $u\in  \T$   such  that   if  $v\in   \T$  and $v\preccurlyeq i$, $v \preccurlyeq j$  then $v \preccurlyeq u$. We also define the lexicographic order $i\leq j$ if either $i \preccurlyeq j$ or $v0  \preccurlyeq i$  and $v1  \preccurlyeq j$  for $v=i\wedge  j$.  Let $X=(X_i, i\in  \T)$ be  a $BMC$  with kernel  $\cp$ and  initial measure $\nu$. For $i\in \T$, we define the $\sigma$-field:
\begin{equation*}\label{eq:field-Fi}
\cf_{i}=\{X_u; u\in \T \text{ such that  $u\leq i$}\}.
\end{equation*}
By construction,  the $\sigma$-fields $(\cf_{i}; \, i\in \T)$ are nested as $\cf_{i}\subset \cf_{j} $ for $i\leq  j$.
\medskip

We define for $n\in \N$, $i\in \G_{n-p_n}$ and $\bF_{n}$ the martingale increments:
\begin{equation}\label{eq:def-DiF}
\Delta_{n,i}(\bF_{n})= N_{n,i}(\bF_{n}) - \E\left[N_{n,i}(\bF_{n}) |\, \cf_i\right] \quad \text{and} \quad \Delta_n(\bF_{n}) = \sum_{i\in \G_{n-p_n}} \Delta_{n,i}(\bF_{n}).
\end{equation}
Thanks to \reff{eq:def-NiF}, we have:
\[
\sum_{i\in \G_{n-p_n}} N_{n, i}(\bF_{n}) = |\G_n|^{-1/2} \sum_{\ell=0}^{p_n}  M_{\G_{n-\ell}} (\tilde f_{\ell,n}) = |\G_n|^{-1/2} \sum_{k=n-p_n}^{n}  M_{\G_{k}} (\tilde f_{n-k,n}).
\]
Using the branching Markov property, and \eqref{eq:def-NiF}, we get for $i\in \G_{n-p_n}$:
\[
\E\left[N_{n,i}(\bF_{n}) |\, \cf_i\right] =\E\left[N_{n,i}(\bF_{n}) |\, X_i\right] = |\G_n|^{-1/2} \sum_{\ell=0}^{p_n} \E_{X_i}\left[M_{\G_{p_n-\ell}}(\tilde f_{\ell,n})\right].
\] 
Assume  that  $n$  is  large  enough  so  that $n-p_n -1\geq k_0$. We have:
\[
N^{[k_0]}_{n, \emptyset}(\bF) = \Delta_n(\bF) + R_0^{k_0} (n)+R_1(n),
\]
where $\Delta_n(\bF)$ and $R_1(n)$ are defined in \reff{eq:def-DiF} and:
\begin{equation*}
R_0^{k_0} (n)= |\G_n|^{-1/2}\, \sum_{k=k_0}^{n-p_n-1} M_{\G_k}(\tilde f_{n-k}) \quad \text{and} \quad R_1(n)= \sum_{i\in \G_{n-p_n}}\E\left[N_{n,i}(\bF_{n}) |\, \cf_i\right].
\end{equation*}



We have the following result:
\begin{lem}\label{lem:Nk0f=Dnf}
Under the assumptions of Theorem \ref{thm:flx} ($2\alpha^{2} < 1$), we have that
\begin{equation*}
\lim_{n \rightarrow \infty} \E\left[\left(N^{[k_0]}_{n, \emptyset}(\bF_n) - \Delta_n(\bF_n)\right)^2\right] = 0.
\end{equation*} 
\end{lem} 
\begin{proof}
We deduce from Remark 5.5 in \cite{BD2} that $\E\left[\left(N^{[k_0]}_{n, \emptyset}(\bF_n) - \Delta_n(\bF_n)\right)^2\right] \leq  a_{0, n} \, \bC_2^2$ for a  sequence $(a_{0,n}, n\in \N)$ which converges to $0$ and does not depend on the sequences $\bF_n$.
\end{proof}
 
We consider  the bracket of the martingale $\Delta_n(\bF_n) $ given by $V(n) = \sum_{i\in \G_{n-p_n}} \E\left[\Delta_{n, i}(\bF_n)^2|\cf_i\right]$. Using \reff{eq:def-NiF} and \reff{eq:def-DiF}, we write:
\begin{equation}\label{eq:def-V}
V(n) = |\G_n|^{-1} \sum_{i\in \G_{n-p_n}} \E_{X_i}\left[\left(\sum_{\ell=0}^{p_n} M_{\G_{p_n-\ell}}(\tilde f_{\ell,n}) \right)^2\right] - R_2(n) = V_1(n) + 2V_2(n) - R_2( n),
\end{equation}
with:
\begin{align*}
V_1(n) & =   |\G_n|^{-1} \sum_{i\in \G_{n-p_n}}  \sum_{\ell=0}^{p_n} \E_{X_i}\left[M_{\G_{p_n-\ell}}(\tilde f_{\ell,n}) ^2 \right] ,
\\
V_2(n)& =  |\G_n|^{-1} \sum_{i\in \G_{n-p_n}} \sum_{0\leq \ell<k\leq p_n} \E_{X_i}\left[M_{\G_{p_n-\ell}}(\tilde f_{\ell,n})  M_{\G_{p_n-k}}(\tilde f_{k,n}) \right], 
\\
R_2( n) & = \sum_{i\in \G_{n-p_n}} \E\left[ N_{n,i} (\bF_{n})|X_i \right] ^2.
\end{align*}
 
\begin{lem}\label{lem:cv-R2-s-E}
Under the assumptions of Theorem \ref{thm:flx} ($2\alpha^2<1$), we have that $R_{2}(n)$ converges in probability towards 0. 
\end{lem}
\begin{proof}
We deduce from Remark 5.7 in \cite{BD2} that $\EE[|R_{2}(n)|] \leq C \bC_{2}^{2} a_{n}$ for a sequence $(a_{n}, n\in \NN)$ which converges to $0$ and does not depend on the sequence $\bF_{n}.$
\end{proof}

\begin{lem}\label{lem:cv-V2-s-E}
Under the assumptions of Theorem \ref{thm:flx} ($2\alpha^2<1$), we have that $V_{2}(n)$ converges in probability towards 0. 
\end{lem}

\begin{proof}
First, we have the following preliminary results. Let $f\in L^2(\mu)$ and recall that $\tilde f=f - \langle \mu, f \rangle$. We deduce from  $\langle \mu, f \rangle=\langle \mu, \cq f \rangle \leq  \norm{\cq f}_\infty \leq  \norm{\cq (f^2)}_\infty ^{1/2}$ that:
\begin{equation}\label{eq:q2-mu-f}
\norm{\cq \tilde f}_\infty \leq  2 \norm{\cq (f^2)}_\infty ^{1/2} \quad\text{and}\quad \norm{\cq (\tilde f^2)}_\infty \leq  4  \norm{\cq (f^2)}_\infty.
\end{equation}
Note that thanks to Assumption \ref{hyp:f} we have, for all $k, \ell,r\in
\N$, and $j>0$:
\begin{equation}\label{eq:majoqtf}
\lim_{n \rightarrow \infty} |\langle \mu, \tilde
f_{k,n}\Qq^{j}\tilde f_{\ell,n} \rangle| = 0 \quad \text{and} \quad
\lim_{n \rightarrow \infty} |\langle \mu,
\cp\left(\Qq^{r}\tilde{f}_{k,n}  \sot  \Qq^{j}\tilde{f}_{\ell,n}
\right)\rangle| = 0. 
\end{equation}
Indeed, we have thanks to Assumption \ref{hyp:f} (iii):
\[
 |\langle \mu, \tilde
f_{k,n}\Qq^{j}\tilde f_{\ell,n} \rangle|
\leq  \norm{\cq \tilde f_{\ell,n} }_\infty \,  \langle \mu, |\tilde
f_{k,n}|  \rangle 
\leq  4 \norm{\cq  f_{\ell,n}^2 }_\infty^{1/2} \langle \mu, |
f_{k,n}|  \rangle 
\leq  4\bQ_2\,   \delta_{k, n}. 
\]
We also have thanks to Assumption \ref{hyp:f} (iii), for 
$g=\cq^{j-1}|\tilde f_{\ell,n}|$ and $r=0$:
\begin{align}
\nonumber
 |\langle \mu,
\cp\left(\Qq^{r}\tilde{f}_{k,n}  \sot  \Qq^{j} \tilde{f}_{\ell,n}
  \right)\rangle| 
&\leq  \langle \mu,
\cp\left(|\tilde{f}_{k,n}|  \sot \cq g \right)\rangle\\
\nonumber
&\leq  \langle \mu, \cp\left(\ind \sot  \cq g  \right) \rangle\langle \mu, |
f_{k,n}|\rangle + \normm{\cp(| f_{k,n}| \sot \cq g)}\\ \nonumber
&\leq 2  \normm{g}\,  \delta_{k,n} \\
\nonumber
&\leq 2\bC_2\,  \delta_{k, n}, 
\end{align}
and for $r\geq 1$ using \reff{eq:q2-mu-f} and that $\langle \mu, \cp
(\ind \sot h) \rangle= \langle \mu, h \rangle$ : 
\[
 |\langle \mu,
\cp\left(\Qq^{r}\tilde{f}_{k,n}  \sot  \Qq^{j}\tilde{f}_{\ell,n} \right)\rangle|
\leq  \langle \mu,
\cp\left(\ind \sot  \cq g  \right)\rangle\,  \norm{\cq^r \tilde f_{k,
  n}}_\infty   
\leq 2  \,\bQ_2\, \delta_{\ell,n}.
\]
Then use that for all $k\in \N$ fixed, we have $\lim_{n\rightarrow \infty } \delta_{k,n}=0$ to conclude that \reff{eq:majoqtf} holds.

\medskip
 
Using \reff{eq:Q2-bis}, we get:
\begin{equation}\label{eq:decom-V2}
V_2(n)= V_5(n)+ V_6(n),
\end{equation}
with
\begin{align*}
V_5(n) &=  |\G_n|^{-1} \sum_{i\in \G_{n-p}} \sum_{0\leq \ell<k\leq  p } 2^{p-\ell} \cq^{p-k} \left( \tilde{f}_{k,n} \cq^{k-\ell} \tilde{f}_{\ell,n}\right)(X_i),\\
V_6(n) &=     |\G_n|^{-1} \sum_{i\in \G_{n-p}} \sum_{0\leq \ell<k<  p }\sum_{r=0}^{p-k-1}  2^{p-\ell+r} \, \cq^{p-1-(r+k)}\left(\cp\left(\cq^r \tilde{f}_{k,n} \sot \cq  ^{k-\ell+r} \tilde{f}_{\ell,n} \right)\right)(X_i).
\end{align*}

First, we consider the term $V_6(n)$. We have:
\begin{equation*}\label{eq:V6}
V_6(n)=|\G_{n-p}|^{-1} M_{\G_{n-p}} (H_{6,n}),
\end{equation*}
with
\[
H_{6,n} = \sum_{\substack{0\leq \ell< k \\ r\geq 0}} h_{k,\ell,r}^{(n)}\,  \ind_{\{r+k<  p\}}  \quad\text{and} \quad h_{k, \ell,r}^{(n)} =  2^{r-\ell} \, \cq^{p-1-(r+k)}\left(\cp\left(\cq^r \tilde{f}_{k,n} \sot \cq  ^{k-\ell+r} \tilde{f}_{\ell,n} \right)\right).
\]

Define
\begin{equation}\label{eq:def-B6nsubL2}
H_6^{[n]}(\bF_{n}) = \sum_{\substack{0\leq \ell< k \\ r\geq 0}}
h_{k,\ell,r} \,  \ind_{\{r+k<  p\}}, 
\end{equation}
with $h_{k,\ell,r} = 2^{r-\ell} \langle \mu, \cp\left(\cq^r \tilde{f}_{k,n} \sot \cq^{k-\ell+r} \tilde{f}_{\ell,n}\right)\rangle = \langle \mu, h_{k, \ell,r}^{(n)}\rangle$.

 \medskip
  
We set $A_{6,n}(\bF_{n})=H_{6,n}-H_6^{[n]}(\bF_{n}) = \sum_{\substack{0\leq \ell< k \\ r\geq 0}} (h_{k,\ell,r}^{(n)}-h_{k,\ell,r}) \,  \ind_{\{r+k<  p\}}$, so that from the    definition     of     $V_{6}(n)$,     we    get     that:
\begin{equation*}\label{eq:V6-H6f}
V_6(n)-  H_6^{[n]}(\bF_{n})=  |\G_{n-p}|^{-1}\, M_{\G_{n-p}}  (A_{6,n}(\bF_{n})).
\end{equation*}
We now study the second moment of $ |\G_{n-p}|^{-1}\, M_{\G_{n-p}}(A_{6,n}(\bF_{n}))$.  Using \reff{eq:EGn-nu0}, we get for $n-p\geq k_0$:
\[
|\G_{n-p}|^{-2}\, \E\left[M_{\G_{n-p}} (A_{6,n}(\bF_{n}))^2 \right] \leq  C\, |\G_{n-p}|^{-1}\, \sum_{j=0}^{n-p} 2^j \normm{\cq^j (A_{6,n}(\bF_{n}))}^2. 
\]
We deduce that
\begin{align}
\nonumber \normm {\cq^j (A_{6,n}(\bF_{n}))} &\leq \sum_{\substack{0\leq \ell< k\\ r\geq 0}}  \normm {\cq^j h_{k,\ell,r}^{(n)}-h_{k,\ell,r}} \ind_{\{r+k<  p\}} \\ 
\nonumber &\leq C \sum_{\substack{0\leq \ell< k\\ r\geq 0}} 2^{r-\ell } \, \alpha^{p-1-(r+k)+j} \normm{\cp\left(\cq^r \tilde f_{k,n} \sot \cq  ^{k-\ell+r} \tilde f_{\ell,n} \right)} \ind_{\{r+k<  p\}}\\
\label{eq:majo-A6_bis-c4q2} &\leq C \bC_2^2\, \alpha^j \, \sum_{\substack{0\leq \ell< k\\ 
r \geq k_{1}}} 2^{r-\ell } \,\alpha^{p-(r+k)} \alpha^{k-\ell+2r} \, \ind_{\{r+k<  p\}}\\
\label{eq:majo-A6-c4q2} &\hspace{0.5cm} + C  \alpha^j  \, \sum_{\substack{0\leq \ell< k \\ 0 \leq r \leq k_{1} - 1}}  2^{-\ell } \,\alpha^{p-k} \normm{\cp\left(\Qq^{r}\tilde{f}_{k,n} \sot \cq  ^{k - \ell + r} \tilde {f}_{\ell,n} \right)} \, \ind_{\{r + k <  p\}},
\end{align}
where  we  used the  triangular  inequality  for the  first  inequality; \reff{eq:L2-erg} for the second; \reff{eq:Hil-SchP-1} for $r\geq k_{1}$  and \reff{eq:L2-erg} again for the third.
The term \reff{eq:majo-A6-c4q2} can be bounded from above using \reff{eq:q2-mu-f} and $\normm{\cp(\Qq^{r}\tilde f_{k,n} \sot \cq^{k-\ell+r} \tilde f_{\ell,n})} \leq \norm{\cq \tilde f_{\ell,n}}_\infty \, \normm{\cp (\Qq^{r} \tilde f_{k,n}\sot \ind)} \leq  2 \bQ_{2}\, \bC_{2}$ as $k > \ell$, and thus \reff{eq:majo-A6_bis-c4q2} and \eqref{eq:majo-A6-c4q2} imply that
\begin{align}
\nonumber \normm {\cq^j (A_{6,n}(\bF_{n}))} &\leq C \bC_{2} (\bC_{2} + \bQ_{2}) \, \alpha^j \, \sum_{\substack{0\leq \ell< k\\ r\geq 0}} 2^{r-\ell } \,\alpha^{p-(r+k)} \alpha^{k-\ell+2r} \, \ind_{\{r+k<  p\}}\\
\label{eq:majo-A6-c4} &\leq C \bC_{2} (\bC_{2} + \bQ_{2}) \, \alpha^j, 
\end{align}
where we used that $\sum_{{0\leq \ell< k,\,  r\geq 0}}  2^{r-\ell} \alpha^{k - \ell + 2r} $ is finite for the last inequality. As $\sum_{j=0}^{\infty } (2\alpha^2)^j$ is finite, we deduce that:
\begin{equation}\label{eq:majo-L2A6}
\E\left[\left( V_6(n)-  H_6^{[n]}(\bF_{n})\right)^2\right] = |\G_{n-p}|^{-2}\, \E\left[M_{\G_{n-p}} (A_{6,n}(\bF_{n}))^2 \right] \leq C \bC^{2}_{2} (\bC_{2} + \bQ_{2})^{2} \, 2^{-(n-p)}.
\end{equation}
\medskip

We now consider the term $V_5(n)$ defined just after \reff{eq:decom-V2}:
\[
V_5(n)=|\G_{n-p}|^{-1} M_{\G_{n-p}} (H_{5,n}),
\]
with
\[
H_{5,n}=\sum_{0\leq \ell< k } h_{k,\ell}^{(n)}\,  \ind_{\{k\leq   p\}} \quad\text{and} \quad h_{k, \ell}^{(n)} =  2^{-\ell} \, \cq^{p-k}\left(\tilde f_{k,n} \cq  ^{k-\ell} \tilde f_{\ell,n} \right). 
\]
We consider the constant
\begin{equation}\label{eq:def-B5nsubL2}
H_5^{[n]}(\bF_{n}) = \!\!\sum_{0\leq \ell< k} \!\! h_{k,\ell}\ind_{\{k\leq p\}} \quad \text{with} \quad h_{k,\ell}=2^{-\ell} \langle \mu, \tilde f_{k,n} \cq^{k-\ell} \tilde f_{\ell,n}\rangle.
\end{equation} 
We set $A_{5,n}(\bF_{n}) = H_{5,n} - H_5^{[n]}(\bF_{n}) = \sum_{0\leq \ell< k } (h_{k,\ell}^{(n)}-h_{k,\ell}) \,  \ind_{\{k\leq   p\}}$, so that from the definition of $V_{5}(n)$,  we  get that:
\begin{equation*}\label{eq:V5-H5f}
V_5(n) -  H_5^{[n]}(\bF_{n}) =  |\G_{n-p}|^{-1}\, M_{\G_{n-p}}  (A_{5,n}(\bF_{n})).
\end{equation*}
We now  study the second  moment of  $|\G_{n-p}|^{-1}\,  M_{\G_{n-p}} (A_{5,n}(\bF_{n}))$. Using \reff{eq:EGn-nu0}, we get for $n-p\geq k_0$:
\[
|\G_{n-p}|^{-2}\, \E\left[M_{\G_{n-p}} (A_{5,n}(\bF_{n}))^2 \right] \leq  C\, |\G_{n-p}|^{-1}\, \sum_{j=0}^{n-p} 2^j \normm{\cq^j (A_{5,n}(\bF_{n}))}^2. 
\]
We also have  that:
\begin{align}
\nonumber \normm {\cq^j(A_{5,n}(\bF_{n}))} & \leq \sum_{0\leq \ell< k}  \normm {\cq^j h_{k,\ell}^{(n)}-h_{k,\ell}} \ind_{\{k\leq   p\}} \\ 
\label{eq:L2A5Hf-q} &\leq C \sum_{0\leq \ell< k} 2^{-\ell } \, \alpha^{p-k+j} \normm{ \tilde f_{k,n} \cq  ^{k-\ell} \tilde f_{\ell,n}} \ind_{\{k\leq   p\}},
\end{align}
where we  used the  triangular inequality for  the first  inequality and \reff{eq:L2-erg} for the last.
The term \reff{eq:L2A5Hf-q} can be bounded from above using $\normm{ \tilde f_{k,n} \cq  ^{k-\ell} \tilde f_{\ell,n}} \leq   \normm{\tilde f_{k,n}}\, \norm{ \cq  ^{k-\ell} \tilde f_{\ell,n}}_\infty \leq \bC_{2} \, \bQ_{2}$ as $k>\ell$. This implies that
\begin{equation*}\label{eq:L2A5Hf} 
\normm {\cq^j(A_{5,n}(\bF_{n}))} \leq C \bC_{2} \bQ_{2} \, \alpha^j. 
\end{equation*}
As $\sum_{j=0}^{\infty} (2\alpha^2)^j$ is finite, we deduce that: 
\begin{equation}\label{eq:majo-L2A5}
\E\left[\left(V_5(n)-  H_5^{[n]}(\bF_{n})\right)^2\right] = |\G_{n-p}|^{-2}\, \E\left[M_{\G_{n-p}} (A_{5,n}(\bF_{n}))^2 \right] \leq  C\, \bC_{2} \bQ_{2}\, 2^{-(n-p)} . 
\end{equation}
\medskip

We deduce from \reff{eq:majo-L2A6}  and \reff{eq:majo-L2A5}, as $V_2(n)=V_5(n)+V_6(n)$ (see \reff{eq:decom-V2}), that: 
\begin{equation}\label{eq:majo-L2A2}
\E\left[\left( V_2(n)-  H_2^{[n]}(\bF_n)\right)^2\right] \leq C \, \left(\bC_2^4+ \bC_2^2\, \bQ_2^2\right)  \, 2^{-(n-p)}, \quad\text{with} \quad H_2^{[n]}(\bF_{n})= H_6^{[n]}(\bF_{n}) + H_5^{[n]}(\bF_{n}).
\end{equation}
Since according  to (ii) in Assumption \ref{hyp:f} $\bC_2$ and $\bQ_2$ are finite, we deduce  that $\lim_{n\rightarrow \infty } V_2(n) - H_2^{[n]}(\bF_n)=0$ in probability.

\medskip

We now check that $\lim_{n\rightarrow \infty } H_2^{[n]}(\bF_n)=0$. Using \reff{eq:def-B6nsubL2} and \reff{eq:def-B5nsubL2}, we get  that:
\[
|H_2^{[n]}(\bF_n)| \leq  \sum_{k>\ell\geq 0} 2^{-\ell} |\langle \mu,   \tilde f_{k,n} \cq^{k-\ell} \tilde f_{\ell,n}\rangle | + \sum_{\substack{k>\ell\geq 0\\ r\geq 0}} 2^{r-\ell} |\langle \mu,\cp\left( \cq^r \tilde f_{k,n} \sot \cq^{k-\ell+r} \tilde f_{\ell,n}  \right)\rangle|.
\]
Recall   the   definition  of   $\Delta$   in     Assumption \ref{hyp:f} (iii). Thanks to   \reff{eq:L2-erg} and   \reff{eq:mp} we have:
\begin{align}\label{eq:majL2qtf1}
|\langle \mu,   \tilde f_{k,n} \cq^{k-\ell} \tilde f_{\ell,n}\rangle| & \leq   \bC_2^2\,  \alpha^{k - \ell},\\ 
\nonumber |\langle \mu,\cp\left( \cq^r \tilde f_{k,n} \sot \cq^{k-\ell+r} \tilde f_{\ell,n}  \right)\rangle| &\leq  C \, \bC_{2}^2 \,  \alpha^{k - \ell+ 2r}. 
\end{align}
Since $\sum_{0\leq \ell< k} 2^{-\ell} \alpha^{k-\ell} + \sum_{\substack{0\leq \ell< k\\ r\geq 0}} 2^{r-\ell} \alpha^{k-\ell+2r} $ is finite, we deduce from \reff{eq:def-B6nsubL2}, \reff{eq:def-B5nsubL2},  \reff{eq:majoqtf} and dominated convergence that $\lim_{n\rightarrow \infty} H_2^{[n]}(\bF_n)=0$. This implies that  $\lim_{n\rightarrow \infty } V_2(n) =0$ in probability. 
\end{proof}


\begin{lem}\label{lem:cv-V1-s-E}
Under the assumptions of Theorem \ref{thm:flx} ($2\alpha^2<1$), we have that $V(n)$ converges in probability towards  $\sigma^2$ defined by \reff{eq:limf-ln}. 
\end{lem}

\begin{proof}
Using \reff{eq:Q2}, we get:
\begin{equation}\label{eq:DV4nsub}
V_1(n) = V_3(n)+ V_4(n),
\end{equation}
with
\begin{align*}
V_3(n) &=  |\G_n|^{-1} \sum_{i\in \G_{n-p}} \sum_{\ell=0}^p 2^{p-\ell}\, \cq^{p-\ell} (\tilde f_{\ell,n}^2)(X_i),\\ 
V_4(n) &=     |\G_n|^{-1} \sum_{i\in \G_{n-p}} \sum_{\ell=0}^{p-1}\,\sum_{k=0}^{p-\ell -1} 2^{p-\ell+k} \,\cq^{p-1-(\ell+k)}\left(\cp\left(\cq^k \tilde f_{\ell,n} \otimes^2\right)\right)(X_i).  
\end{align*}
We first consider the term $V_4(n)$. We have:
\[
V_4(n)=|\G_{n-p}|^{-1} M_{\G_{n-p}} (H_{4,n}),
\]
with:
\[
H_{4,n} = \sum_{\ell\geq 0, \, k\geq 0} h_{\ell,
  k}^{(n)}\,  \ind_{\{\ell+k<  p\}} \quad \text{and} \quad h_{\ell,
  k}^{(n)} = 2^{k-\ell}\,  \cq^{p-1-(\ell+k)} \left(\cp\left(\cq^k
    \tilde f_\ell \otimes^2 \right)\right). 
\]
Define the constant
\begin{equation}\label{eq:def-B4nsubL2}
H_4^{[n]}(\bF_{n}) = \sum_{\ell\geq 0, \, k\geq 0} h_{\ell,k} \,  \ind_{\{\ell+k<  p\}} \quad \text{with} \quad h_{\ell, k} = 2^{k-\ell}\,  \langle \mu,\cp\left(\cq^k \tilde f_{\ell,n} \otimes^2 \right) \rangle.
\end{equation}

\medskip

We set $
A_{4,n}(\bF_{n}) = H_{4,n} - H_4^{[n]}(\bF_{n}) = \sum_{\ell\geq 0, \, k\geq 0} (h_{\ell,k}^{(n)}-h_{\ell,k}) \,  \ind_{\{\ell+k<  p\}} $, so that from the definition of  $V_{4}(n)$, we get that: 
\begin{equation*}\label{eq:V6-H4f}
V_4(n) -  H_4^{[n]}(\bF_{n}) =  |\G_{n-p}|^{-1}\, M_{\G_{n-p}}  (A_{4,n}(\bF_{n})).
\end{equation*}
We now study the second moment of $ |\G_{n-p}|^{-1}\, M_{\G_{n-p}}(A_{4,n}(\bF_{n}))$.  Using \reff{eq:EGn-nu0}, we get for $n-p\geq k_0$:
\[
|\G_{n-p}|^{-2}\, \E\left[M_{\G_{n-p}} (A_{4,n}(\bF_{n}))^2 \right] \leq  C\, |\G_{n-p}|^{-1}\, \sum_{j=0}^{n-p} 2^j \normm{\cq^j (A_{4,n}(\bF_{n}))}^2. 
\]
Using \reff{eq:majo-pfxf} and \eqref{eq:q2-mu-f}, we obtain  that for all $0 \leq k < k_{1},$ $\normm{\cp(\Qq^{k}\tilde f_{\ell,n} \otimes^{2}} \leq \|\Qq \tilde{f}^{2}_{\ell,n}\|_{L^{2}(\mu)} \leq 4 \bQ_{2}^{2}$. 
We deduce that:
\begin{align}
\nonumber \normm {\cq^j (A_{4,n}(\bF_{n}))} &\leq \sum_{\ell\geq 0,\, k\geq 0}  \normm {\cq^j h_{\ell,k}^{(n)}-h_{\ell,k}} \ind_{\{\ell+k<  p\}} \\ 
\nonumber &\leq C \sum_{\ell\geq 0,\, k\geq 0}  2^{k-\ell} \, \alpha^{p-1-(\ell+k)+j} \normm{\cp\left(\cq^k \tilde f_{\ell,n} \otimes^2 \right)} \ind_{\{\ell+k<  p\}}\\
\nonumber &\leq C \, \bC_{2}^{2} \, \alpha^j\!\! \!  \sum_{\ell\geq 0,\, k \geq k_{1}} \! \!\! 2^{k-\ell } \,\alpha^{p-(\ell+k)} \alpha^{2k} \, \ind_{\{\ell+k<  p\}}\\
\nonumber&\hspace{2cm} +  C \, \alpha^j \sum_{\substack{\ell\geq 0 \\ 0 \leq k < k_{1}}} 2^{k - \ell} \, \alpha^{p -(\ell + k)} \normm{\cp\left(\Qq^{k}\tilde f_{\ell,n} \otimes ^2\right)} \ind_{\{\ell<  p\}}\\
\nonumber&\leq C \, (\bC_{2}^{2} +  4\bQ_{2}^{2}) \, \alpha^j,
\end{align}
where  we  used the  triangular  inequality  for the  first  inequality;\reff{eq:L2-erg}   for   the   second;   \reff{eq:Hil-SchP-1}   for   $k\geq k_{1}$ and \reff{eq:L2-erg} again for the third; \reff{eq:majo-pfxf} and \eqref{eq:q2-mu-f} for the last. 
As $\sum_{j=0}^{\infty } (2\alpha^2)^j$ is finite, we deduce that:
\begin{equation}\label{eq:majo-L2A4}
\E\left[\left(V_4(n)-  H_4^{[n]}(\bF_{n})\right)^2\right] = |\G_{n-p}|^{-2}\, \E\left[M_{\G_{n-p}} (A_{4,n}(\bF_{n}))^2 \right]\leq C \, (\bC_{2}^{2} + \bQ_{2}^{2}) \, 2^{-(n-p)}.
\end{equation}
\medskip

We now consider the term $V_3(n)$ defined just after \reff{eq:DV4nsub}:
\[
V_3(n)=|\G_{n-p}|^{-1} M_{\G_{n-p}} (H_{3,n}),
\]
with 
\[
H_{3,n} = \sum_{\ell\geq 0} h_{\ell}^{(n)}\,  \ind_{\{\ell\leq   p\}} \quad \text{and} \quad h_{\ell}^{(n)} = 2^{-\ell}\,  \cq^{p-\ell} \left(\tilde f_{\ell,n}^2\right). 
\]
We consider the constant
\begin{equation}\label{eq:def-B3nsubL2}
H_3^{[n]}(\bF_{n}) = \sum_{\ell\geq 0} h_{\ell} \,  \ind_{\{\ell\leq p\}} \quad \text{with} \quad h_\ell= 2^{-\ell} \langle \mu,  \tilde f_{\ell,n}^2 \rangle=\langle \mu, h^{(n)}_\ell \rangle.  
\end{equation} 

\medskip

We set $A_{3,n}(\bF_{n}) = H_{3,n} - H_3^{[n]}(\bF_{n}) = \sum_{\ell\geq 0}(h_{\ell}^{(n)} - h_{\ell}) \,  \ind_{\{\ell\leq   p\}} $, so that from the definition of $V_{3}(n)$,  we get that:
\begin{equation}\label{eq:V3-H3f}
V_3(n) -  H_3^{[n]}(\bF_{n}) =  |\G_{n-p}|^{-1}\, M_{\G_{n-p}}  (A_{3,n}(\bF_{n})).
\end{equation}
We now study the second moment of $ |\G_{n-p}|^{-1}\, M_{\G_{n-p}}(A_{3,n}(\bF_{n}))$.  Using \reff{eq:EGn-nu0}, we get for $n-p\geq k_0$: 
\begin{equation}\label{eq:majoA3}
|\G_{n-p}|^{-2}\, \E\left[M_{\G_{n-p}} (A_{3,n}(\bF_{n}))^2 \right] \leq  C\, |\G_{n-p}|^{-1}\, \sum_{j=0}^{n-p} 2^j \normm{\cq^j (A_{3,n}(\bF_{n}))}^2. 
\end{equation}
Recall $c_{k}(\bF_{n})$ and $q_{k}(\bF_{n})$ defined in \eqref{eq:def-ck}. We have  that
\begin{align}
\nonumber \normm {\cq^j (A_{3,n}(\bF_{n}))} &\leq \sum_{\ell\geq 0}  \normm {\cq^j h_{\ell}^{(n)}-h_{\ell}} \ind_{\{\ell\leq   p\}} \\ 
\label{eq:L2A3Hf-q} &\leq C \sum_{\ell\geq 0}  2^{-\ell} \,\normm{ \cq^{j+p-\ell} \tilde g} \ind_{\{\ell\leq  p\}}\quad\text{with}\quad g=\tilde f^2_{\ell,n}\\
\nonumber&=  2^{-p}  \normm{  \tilde g}\, \ind_{\{j=0\}} +\sum_{\ell= 0}^p  2^{-\ell}  \normm{ \cq^{j+p-\ell-1} \cq\tilde g} \ind_{\{j+p-\ell>0\}} \\
\nonumber&\leq C \, c_4^2(\bF_{n})\,  2^{-p} \ind_{\{j=0\}} + C \sum_{\ell\geq 0}  2^{-\ell} \, \alpha^{j+p-\ell} \normm{\cq \tilde g } \\
\nonumber&\leq C \, c_4^2(\bF_{n})\,  2^{-p} \ind_{\{j=0\}} +  C \, q^2_2(\bF_{n})\, \alpha^j,
\end{align}
where  we  used the  triangular  inequality  for the  first  inequality; \eqref{eq:L2-erg} for the third and \reff{eq:q2-mu-f} for the last inequality. As $\sum_{j=0}^{\infty } (2\alpha^2)^j$ is finite, we deduce that:
\begin{equation}\label{eq:majo-L2A3}
\E\left[\left( V_3(n)-  H_3^{[n]}(\bF_{n})\right)^2\right] = |\G_{n-p}|^{-2}\, \E\left[M_{\G_{n-p}} (A_{3,n}(\bF_{n}))^2 \right] \leq C \, c_4^4(\bF_{n})\,  2^{-n}  +  C \, q^4_2(\bF_{n}) \,2^{-(n-p)}.
\end{equation}

As $V_1 = V_4+V_3$, we deduce from \reff{eq:majo-L2A4} and \reff{eq:majo-L2A3} that:
\begin{equation*}\label{eq:I-V1n-H1nS}
\E\left[\left( V_1(n)-  H_1^{[n]}(\bF_n)\right)^2\right] \leq C \, \left((c_2^4(\bF_n)+ q_2^4(\bF_n))  \, 2^{-(n-p)}  + c_4^4(\bF_n)\,  2^{-n}\right),
\end{equation*}
with  $H_1^{[n]}(\bF_{n}) = H_3^{[n]} (\bF_{n})+ H_4^{[n]}(\bF_{n})$.  Since $c_4^4(\bF_n)  \leq c_2^2(\bF_n) \, c_\infty^2  (\bF_n) \leq  C_\rho \, c_2^2(\bF_n) \, 2^{2n\rho}$ with $\rho\in (0, 1/2)$ and some finite constant $C_\rho$ according to (i) in Assumption \ref{hyp:f}, and since $\lim_{n\rightarrow \infty }p/n=1$ so that $2^{-n(1-2\rho)}\leq  2^{-(n-p)}$ (at least for $n$ large enough), we deduce from  (ii) in Assumption \ref{hyp:f} that:
\begin{equation}\label{eq:majo-L2A1}
\E\left[\left( V_1(n)-  H_1^{[n]}(\bF_n)\right)^2\right] \leq C \, \left(\bC_2^4+ \bQ_2^4+ C_\rho \bC_2^2\right)\, 2^{-(n-p)}
\end{equation}
and thus $\lim_{n\rightarrow \infty } V_1(n) - H_1^{[n]}(\bF_n)=0$ in probability.

We check  that $\lim_{n\rightarrow \infty } H_1^{[n]}(\bF_n)=\sigma^2$. Recall  (see \reff{eq:def-B3nsubL2} and \reff{eq:def-B4nsubL2}) that:
\[
H_3^{[n]} (\bF_n) =  \sum_{\ell\geq 0} 2^{-\ell} \langle \mu,  \tilde f_{\ell,n}^2 \rangle\,  \ind_{\{\ell\leq p\}}  \quad \text{and} \quad |H_4^{[n]}(\bF_n)| \leq  \sum_{\ell\geq 0, \, k\geq 0} 2^{k-\ell}\,  |\langle \mu,\cp\left(\cq^k \tilde f_{\ell,n} \otimes^2\right) \rangle|.
\]
Thanks to \reff{eq:majo-pfxf} and \reff{eq:L2-erg}, we have:
\[
|\langle \mu,\cp\left(\cq^k \tilde f_{\ell,n} \otimes^2\right) \rangle| \leq \normm{\cq^k \tilde f_{\ell,n}}^2 \leq  C\, \alpha^{2k}\normm{f_{\ell,n}}^2 \leq  C\, \alpha^{2k} \, \bC_{2}^{2}.
\]
Using Assumption \ref{hyp:f}  (iii), we get that
\begin{equation}\label{eq:majo-mPf22}
|\langle \mu, \cp(\tilde f_{\ell, n} \otimes^2 )\rangle|\leq |\langle \mu, \cp(f_{\ell, n} \otimes^2) \rangle|+ \langle \mu, f_{\ell, n} \rangle^2\leq(1+\Delta)\delta_{\ell, n}.
\end{equation}
We deduce from \reff{eq:majoqtf} (for $k\geq 1$) and the previous upper-bound   (for $k=0$) and dominated convergence that $\lim_{n\rightarrow \infty } H_4^{[n]}(\bF_n)=0$. 

We now prove that  $\lim_{n\rightarrow \infty } H_3^{[n]}(\bF_n)=\sigma^2$. We define $\sigma_n^2=\sum_{\ell=0}^n 2^{-\ell} \norm{f_{\ell,n} }_{L^2(\mu)} ^2$, so that by  Assumption \ref{hyp:f}  (iv), $\lim_{n\rightarrow\infty } \sigma^2_n=\sigma^2$. 
We have:
\[
|H_3^{[n]}(\bF_n)- \sigma_n^2| \leq \sum_{\ell=p+1}^n 2^{-\ell} \langle \mu, f_{\ell,n}^2 \rangle  + \sum_{\ell=0}^p 2^{-\ell} \langle \mu, f_{\ell,n} \rangle ^2 \leq  \bC_2^2 2^{-p} + \Delta \sum_{\ell=0}^p 2^{-\ell} \delta_{\ell, n}. 
\]
Then use dominated convergence to deduce that $\lim_{n\rightarrow \infty} |H_3^{[n]}(\bF_n)- \sigma_n^2|=0$. This implies that $\lim_{n\rightarrow \infty } V_1(n)=\sigma^2$ in probability.


\end{proof}
Using \eqref{eq:def-V}, we have the following result as a direct consequence of Lemmas \ref{lem:cv-R2-s-E}, \ref{lem:cv-V2-s-E} and \ref{lem:cv-V1-s-E}.
\begin{lem}\label{lem:cv-var-s-L2}
Under the assumptions of Theorem \ref{thm:flx} ($2\alpha^2<1$), we have that $V(n)$
converges in probability towards  $\sigma^2$ defined by
\reff{eq:limf-ln}. 
\end{lem}
\medskip

We now check the Lindeberg's condition using a fourth moment
condition. We set
\begin{equation}\label{eq:def-R3}
R_3(n) = \sum_{i\in \G_{n-p_n}} \E\left[\Delta_{n,i}(\bF_{n})^4\right].
\end{equation} 

\begin{lem}\label{lem:cv-R3-L2-n}
Under the assumptions of Theorem \ref{thm:flx} ($2\alpha^2<1$), we get $\lim_{n\rightarrow \infty } R_3(n)=0$. 
\end{lem}


\begin{proof}[Proof of Lemma \ref{lem:cv-R3-L2-n}]
We have:
\begin{align}
\nonumber R_3(n) &\leq   16 \sum_{i\in \G_{n-p}}\E\left[N_{n,i}(\bF_{n})^4\right]\\ 
\nonumber &\leq   16 (p+1)^3 \sum_{\ell=0}^p \sum_{i\in \G_{n-p}} \E\left[N_{n,i}^\ell(\tilde f_{\ell,n})^4\right],
\end{align}
where we used that $(\sum_{k=0}^r a_k)^4 \leq  (r+1)^3 \sum_{k=0}^r a_k^4$ for the two inequalities (resp. with $r=1$ and $r=p$) and also Jensen inequality and \reff{eq:def-DiF} for the first and \reff{eq:def-NiF} for the last.  Using \reff{eq:def-Nnil}, we get:
\begin{equation*}\label{eq:def-R3hnl}
\E\left[N_{n,i}^\ell(\tilde f_{\ell,n})^4\right] = |\G_n|^{-2} \E\left[h_{n,\ell}(X_i)\right], \quad\text{with}\quad h_{n,\ell}(x)=\E_x\left[M_{\G_{p-\ell}}(\tilde f_{\ell,n}) ^4\right], 
\end{equation*}
so  that:
\[
R_3(n) \leq   C n^3 \sum_{\ell=0}^p \sum_{i\in \G_{n-p}} |\G_n|^{-2} \E\left[h_{n,\ell}(X_i)\right]. 
\]
Using \reff{eq:EGn-nu0} (with $f$ and $n$ replaced by $h_{n, \ell}$ and $n-p$), we get that:
\begin{equation}\label{eq:R3-L2-9}
R_3(n) \leq   C\,  n^3\, 2^{-n-p}\,   \sum_{\ell=0}^p \E_\mu\left[M_{\G_{p-\ell}}(\tilde f_{\ell,n}) ^4\right].
\end{equation}
Now we give the main steps to get  an upper bound of
$\E_\mu\left[M_{\G_{p-\ell}}(\tilde f_{\ell,n}) ^4\right]$. Recall that:
\[
\norm{\tilde f_{\ell,n}}_{L^4(\mu)}\leq  C\, c_4(\bF_{n}).
\]
We have:
\begin{equation}\label{eq:L2R3-l=p}
\E_\mu\left[M_{\G_{p-\ell}}(\tilde f_{\ell,n}) ^4\right] \leq C\,   c_4^4(\bF_{n}) \quad \text{for $\ell\in \{p - k_{1} - 1, \ldots, p\}$.}
\end{equation}
Now  we consider  the case  $0\leq \ell  \leq p - k_{1} - 2$. Let the  functions $\psi_{j,p-\ell}$, with $1\leq  j\leq 9$, from Lemma  \ref{lem:M4}, with $f$  replaced by  $\tilde  f_{\ell,n}$ so that for $\ell \in \{0, \ldots,
 p - k_{1} - 2\} $
\begin{equation}\label{eq:R3L2-Em4}
\E_\mu\left[M_{\G_{p-\ell}}(\tilde f_{\ell,n}) ^4\right] = \sum_{j=1}^9 \langle \mu, \psi_{j,p-\ell}\rangle.
\end{equation}

\medskip 

We now look precisely at the terms in
 \reff{eq:R3L2-Em4}. 
We set $h_k = \cq^{k-1} \tilde f_{\ell,n}$ so that for $k\in \N^{*}$:
\begin{equation}\label{eq:hk}
\normm{h_k}\leq  C \, \alpha^k \bC_{2} \quad \text{and} \quad \norm{h_k}_{L^4(\mu)} \leq  C \,  c_4(\bF_{n}).
\end{equation}
We recall the notation $f \otimes f = f\otimes^2$. We deduce for $k\geq k_{1} + 1$ from \reff{eq:Hil-SchP-1} applied with $h_k = \cq^{k_{1}} h_{k - k_{1}}$ and for $1 \leq k \leq k_{1}$ from \eqref{eq:majo-pfxf} and \reff{eq:q2-mu-f} that: 
\begin{equation}\label{eq:phk}
\normm{\cp(h_k\otimes ^2)} \leq 
\begin{cases}
C\, \alpha^{2k} \bC_{2}^{2}& \text{for $k\geq k_{1} + 1$,}   \\
C\,  \bQ_{2}^{2}& \text{for $k \in \{1, \ldots, k_{1}\}$.}   
\end{cases}
\end{equation}


\medskip
\textbf{Upper bound of $\langle \mu, |\psi_{1,p-\ell}| \rangle$}. We have:
\begin{equation}\label{eq:L2R3p1}
\langle \mu, |\psi_{1,p-\ell}| \rangle\leq C\,  2^{p-\ell}\, \langle \mu,\cq^{p-\ell} (\tilde f_{\ell,n}^4) \rangle \leq C\,  2^{p-\ell}\, c_4^4(\bF_{n}).
\end{equation} 

\medskip
\textbf{Upper bound of $|\langle \mu, \psi_{2,p-\ell} \rangle|$}.
We set $g=(\tilde f_ {\ell,n})^3$. Then we have
\begin{align}
|\langle \mu, \psi_{2, p-\ell} \rangle| &\leq  C    2^{2(p-\ell)}\,\sum_{k=0}^{p-\ell-1} 2^{-k} |\langle \mu, \cq^k \cp  \left( \cq^{p-\ell-k - 1}( (\tilde f_{\ell, n})^3) \sot h_{p-\ell- k}\right) \rangle| \nonumber \\
&= C    2^{2(p-\ell)}\, \sum_{k=0}^{p-\ell-1} 2^{-k} |\langle \mu, \cp  \left( \cq^{p-\ell-k - 1}( \tilde g) \sot h_{p-\ell- k}\right) \rangle| \nonumber \\
\label{eq:L2R3p2-astuce} &\leq  C    2^{2(p-\ell)}\, \sum_{k=0}^{p-\ell-1} 2^{-k} \,  \normm{ \cq^{p-\ell-k - 1} \tilde g } \normm{ h_{p-\ell
  -k}}  \\ 
&\leq  C    2^{2(p-\ell)}\,\sum_{k=0}^{p-\ell-1} 2^{-k} \, \alpha^{2(p-k-\ell)}   \normm{  g }\normm{ \tilde f_{\ell, n}} \nonumber \\
&\leq   C\, 2^{p-\ell}\, c_6^3(\bF_n)\, \bC_2, \nonumber
\end{align}
where we used that $\langle  \mu, \cp  (\ind\sot h_{p-\ell-k}  )\rangle=2 \langle  \mu, \cq h_{p-\ell -k}\rangle=0$ for   the   equality,   \reff{eq:mp}    for   the   second   inequality, \reff{eq:L2-erg} and \reff{eq:hk}  for the third. 

\medskip
\textbf{Upper bound of $\langle \mu, |\psi_{3,p-\ell} \rangle|.$}
Using \reff{eq:mp}, we easily get:
\begin{equation}\label{eq:L2R3p3}
\langle \mu,| \psi_{3, p-\ell}| \rangle \leq  C \, 2^{2(p-\ell)} \sum_{k=0}^{p-\ell-1} 2^{-k}\, \langle \mu, \cq^k\cp \left( \cq^{p-\ell-k - 1} (\tilde f_\ell^2) \otimes^2\right) \rangle.  
\end{equation}
We deduce from \reff{eq:L2R3p3}, distinguishing  according to  $k=p-\ell -1$  (then use \reff{eq:mp}) and $k \leq p-\ell -2$    (then     use $|\cq(\tilde f_{\ell, n}^2)|\leq 4 \bQ_2^2$, see \reff{eq:q2-mu-f}) that:
\begin{equation}\label{eq:L2R3p3-fln}
\langle \mu,| \psi_{3, p-\ell}| \rangle
\leq  C\, 2^{p-\ell} \, c_4^4(\bF_n) + C\, 2^{2(p-\ell)} \, \bQ_2^2 \bC_2^2.
\end{equation}

\medskip
\textbf{Upper bound of $\langle \mu, |\psi_{4,p-\ell} \rangle|$}. Using \reff{eq:mp} and then \reff{eq:phk} with  $p-\ell -1\geq k_{1} + 1$, we get:
\begin{align}
\langle \mu,| \psi_{4, p-\ell}| \rangle & \leq C \, 2^{4(p-\ell)} \, \langle \mu, \cp \left( |\cp(h_{p-\ell-1}\otimes^2)\otimes^2|\right) \rangle \nonumber\\ 
& \leq C \, 2^{4(p-\ell)} \, \normm{\cp(h_{p-\ell-1}\otimes^2)}^2  \label{eq:L2R3p4-sub}\\
& \leq C \, 2^{4(p-\ell)} \, \alpha^{4(p-\ell)} \, \bC_{2}^{4} \nonumber\\
& \leq C \, 2^{2(p-\ell)} \, \bC_{2}^{4}. \nonumber
\end{align}
 
\medskip
\textbf{Upper bound of $\langle \mu, |\psi_{5,p-\ell} \rangle|$}. 
We have:
\[
\langle \mu,| \psi_{5, p-\ell}| \rangle \leq C \, 2^{4(p-\ell)} \, \sum_{k = 2}^{p - \ell - k_{1} - 1} \sum_{r = 0}^{k -1} 2^{-r}\,  \Gamma_{k,r}^{[5]}  \,+ \, C \, 2^{4(p-\ell)} \,  \sum_{k = p - \ell - k_{1} - 2}^{p - \ell  - 1} \sum_{r = 0}^{k -1} 2^{-r}\,  \Gamma_{k,r}^{[5]},
\]
with 
\[
\Gamma_{k,r}^{[5]}= 2^{-2k }  \langle \mu, \cp \left( \cq^{k -r- 1} |\cp (h_{p-\ell- k} \otimes^2)| \otimes^2 \right) \rangle.
\]
Using \reff{eq:mp}  and then  \reff{eq:phk}, we get:
\begin{equation}\label{eq:L2R3p5-Scrit}
\Gamma_{k,r}^{[5]} \leq C \,  2^{-2k}  \normm{\cp (h_{p-\ell- k} \otimes^2)}^2
\end{equation} 
Using \eqref{eq:majo-pfxf} and (ii) of Assumption \ref{hyp:f}, we get, for $m \in \{0, \ldots, k_{1} - 1\}$, $\cp(\Qq^{m}\tilde f_{\ell, n} \otimes^2) \leq  \Qq^{m}\cq (\tilde f_{\ell, n}^2)
\leq 4 \bQ_2^2$ and then, for all $k \in \{p-\ell-k_{1}-2, \ldots, p-\ell-1\},$ we deduce that
\begin{equation}\label{eq:IPflnpsi5} 
\normm{\cp \left(\Qq^{p-\ell-k-1}\tilde f_{\ell, n} \otimes^2 \right)} \leq C\, \bQ_2\, \bC_2. 
\end{equation}
Using \eqref{eq:Hil-SchP-1} and \eqref{eq:L2-erg}, we get, for all $k \in \{2, \ldots, k - \ell - k_{1} - 1\}$, 
\begin{equation}\label{eq:G5-eq2}
\Gamma_{k,r}^{[5]} \leq C \, 2^{-2k} \alpha^{4(p-\ell -k)} \bC_{2}^{4}.
\end{equation} 
From \eqref{eq:IPflnpsi5} and \eqref{eq:G5-eq2} we deduce that 
\begin{equation*}
\langle \mu, |\psi_{5,p-\ell}| \rangle\leq C\,  2^{2(p-\ell)}\, \bC_2^2(\bQ_2^2 + \bC_2^2).
\end{equation*}

\medskip
\textbf{Upper bound of $\langle \mu, |\psi_{6,p-\ell} |\rangle$}.
We have:
\begin{equation*}
\langle \mu,| \psi_{6, p-\ell}| \rangle \leq C\, 2^{3(p-\ell)} \, \sum_{k=1}^{p-\ell-1} \sum_{r=0}^{k -1} 2^{-r}\, \Gamma_{k,r}^{[6]},
\end{equation*}
with 
\begin{equation*}
\Gamma_{k,r}^{[6]}=  2^{-k } \, \langle \mu,  \cq^r \cp \left(\cq^{k -r-1}|\cp \left( h_{p-\ell-k} \otimes^2 \right)|\sot \cq^{p-\ell-r-1}(\tilde f^2_{\ell,n}) \right) \rangle.
\end{equation*}
Using \reff{eq:mp} and then  \reff{eq:phk}, we get:
\begin{equation}\label{eq:L2R3p6-Scrit} 
\Gamma_{k,r}^{[6]}  \leq  C\, 2^{-k } \, \normm{\cp \left( h_{p-\ell-k} \otimes^2 \right)} \normm{\cq^{p-\ell-r-1}(\tilde f^2_{\ell,n})}. 
\end{equation}
Distinguishing the cases $k > p - \ell - k_{1} - 1$ and $k \leq p - \ell - k_{1} - 1$ and using that $\normm{\cq^{j}(\tilde f^2_{\ell,n})}\leq  4\bQ_2\, \min(\bQ_2,\bC_2)$ for all $j\in \N^*$  (see \reff{eq:q2-mu-f}), \eqref{eq:IPflnpsi5} and \eqref{eq:phk}, we get:
\begin{equation*}
\Gamma_{k,r}^{[6]} \leq C\, 2^{-k} \, \bQ_2\, \min(\bQ_2,\bC_2)\,  \ind_{\{k > p-\ell - k_{1} - 1\}} +   C\, 2^{-k } \, \alpha^{2(p-\ell-k)}\, \bC_{2}^{2} \, \bQ_{2}^{2}\,\ind_{\{k \leq  p - \ell - k_{1} - 1\}} 
\end{equation*}
From the previous inequality, we conclude that
\begin{equation*}
\langle \mu, |\psi_{6,p-\ell}| \rangle\leq C\,  2^{2(p-\ell)}\, \bC_2^2\bQ_2^2.
\end{equation*}

\medskip

\textbf{Upper bound of $|\langle \mu, \psi_{7,p-\ell} \rangle|$}. 
We have:
\begin{equation}\label{eq:y7}
|\langle\mu, \psi_{7, p-\ell} \rangle| \leq C\, 2^{3(p-\ell)} \, \sum_{k=1}^{p-\ell-1} \sum_{r=0}^{k -1} 2^{-r}\, \Gamma_{k,r}^{[7]},
\end{equation}
with 
\begin{equation}\label{eq:L2R3p7-G}
\Gamma_{k,r}^{[7]}= 2^{-k }  |\langle \mu , \cq^r \cp \left( \cq^{k -r-1}\cp \left( h_{p-\ell-k} \sot  \cq^{p-\ell-k -1} (\tilde f_{\ell,n}^2) \right)\sot h_{p-\ell -r} \right)\rangle|.
\end{equation}
When $k \geq p - \ell - k_{1}$, setting $g = \cp  \left( \Qq^{p - \ell - k - 1}\tilde f_{\ell,n} \sot  \Qq^{p - \ell - k - 1}\tilde f_{\ell,n}^2\right)$,  we get that:
\begin{align}
\nonumber  \Gamma_{k,r}^{[7]} &= 2^{- k}  |\langle \mu ,  \cp \left(\cq^{k - r - 1} g \sot h_{p-\ell -r} \right)\rangle|\\ 
\nonumber &= 2^{- k}  |\langle \mu ,  \cp \left(\cq^{k - r - 1} \tilde g \sot h_{p-\ell -r} \right)\rangle|\\ 
\label{eq:Gamm7-1} &\leq  C\,  2^{- k} \normm{\cq^{k - r - 1} \tilde g }\normm{h_{p-\ell -r}}
\end{align}
where we used for $\langle  \mu, \cp(\ind \sot h_{p-\ell -r}\rangle=0$
for the second  equality; \reff{eq:mp} for the first  inequality; Using \eqref{eq:L2-erg} twice (for the first and the last inequality), \eqref{eq:Pfog2} and (ii) of Assumption \ref{hyp:f} for the second inequality, we get
\begin{equation*}
\normm{\cq^{k - r - 1} \tilde g } \leq C \alpha^{k-r} \|g\|_{L^{2}(\mu)} \leq C \bQ_{2}^{2} \alpha^{k-r} 2^{n\rho} \quad \text{and} \quad \normm{h_{p-\ell -r}} \leq C \bC_{2} \alpha^{p - \ell - r}.
\end{equation*} 
Using that $p - \ell - k_{1} \leq k \leq p - \ell - 1$ and putting the last inequalities in \eqref{eq:Gamm7-1}, we deduce that
\begin{equation*} 
\Gamma_{k,r}^{[7]} \leq  C\, \bC_{2}\, \bQ_{2}^{2}\,  2^{-(p-\ell)}  \alpha^{2(p-\ell -r)} \,  2^{n\rho}.
\end{equation*}
We now consider $k \leq p - \ell - k_{1} - 1$. We have:
\begin{align}
\Gamma_{k,r}^{[7]} & \leq C\,     2^{-k }  \normm{\cp \left( h_{p-\ell-k} \sot  \cq^{p-\ell-k -1} (\tilde f_{\ell,n}^2) \right)}\normm{h_{p-\ell -r} } \nonumber\\
\nonumber &\leq  C\, 2^{-k }  \normm{ h_{p-\ell-k-1}}\normm{\cq^{p-\ell-k -2} (\tilde f_{\ell,n}^2)} \, \normm{h_{p-\ell -r} }\\
&\leq C\,     2^{-k }  \normm{h_{p-\ell-k} } \bQ_2^2\, \normm{h_{p-\ell -r} } \nonumber \\
& \leq  C\, 2^{-k }  \, \alpha ^{2(p-\ell-k)}\, \bC_2^2 \,\bQ_2^2, \nonumber
\end{align}
where we used  \reff{eq:mp} for the first  inequality; \reff{eq:Hil-SchP-1} for the  second; and \reff{eq:hk}  for the two lasts. We deduce from \reff{eq:y7} that:
\[
|\langle\mu, \psi_{7, p-\ell} \rangle| \leq C\,  2^{p-\ell}\, 2^{n\rho}\, \bC_2 + C\, 2^{2(p-\ell)}\, \bC_2^2\bQ_2^2.
\]

\medskip
\textbf{Upper bound of $\langle \mu, |\psi_{8,p-\ell} |\rangle$}. 
We have:
\begin{equation}\label{eq:y8}
\langle \mu,| \psi_{8, p-\ell}| \rangle \leq   C\,  2^{4(p-\ell)} \,
\sum_{k=2}^{p-\ell-1} \sum_{r=1}^{k -1} \sum_{j=0}^{r-1} \,  2^{-j }\, \Gamma_{k,r,j}^{[8]},
\end{equation}
with 
\[
\Gamma_{k,r,j}^{[8]}
\leq  2^{-k-r } \langle \mu,\cq^j\cp \left(|
\cq^{r-j-1}\cp \left( h_{p-\ell-r} \otimes^2 \right)|
\sot
|\cq^{k-j-1}\cp \left( h_{p-\ell-k} \otimes^2 \right)|
\right) \rangle.
\]
When $k \geq p - \ell - k_{1}$ and $r > p - \ell - k_{1}$, we have, according to \reff{eq:majo-pfxf} and \reff{eq:q2-mu-f}:
\begin{equation}\label{eq:IPhp-l-kScri}
\cp \left( h_{p-\ell-k} \otimes^2 \right) \leq  \cq (h_{p-\ell-k}^2) \leq 4 \bQ_2^2 \quad \text{and} \quad  \cp \left( h_{p-\ell-r} \otimes^2 \right) \leq \cq (h_{p - \ell - r}^2) \leq 4 \bQ_2^2.
\end{equation}
Distinguishing the three cases $k < p - \ell - k_{1}$, $k \geq p - \ell - k_{1}$ and $r > p - \ell - k_{1}$, $k \geq p - \ell - k_{1}$ and $r \leq p - \ell - k_{1}$,  using \reff{eq:mp}, \eqref{eq:phk} and \reff{eq:IPhp-l-kScri} (noticing that $p - \ell - r \geq k_{1} + 2$ if $k < p - \ell - k_{1}$), we get:
\begin{align}
\Gamma_{k,r,j}^{[8]} &\leq C \, 2^{-k-r } \normm{\cp \left( h_{p-\ell-r} \otimes^2 \right)} \normm{\cp \left( h_{p-\ell-k} \otimes^2 \right)} \nonumber\\
&\leq \begin{cases} C \, \bC_{2}^{4} \, 2^{-k-r } \, \alpha^{4(p-\ell)} \alpha^{-2(k+r)} & \text{if $k < p - \ell - k_{1}$}  \\ 
C \, \bC_{2}^{2} \, \bQ_{2}^{2} \,  2^{-k-r } \, \alpha^{2(p-\ell -r)} & \text{if $k \geq p - \ell - k_{1}$ and $r \leq p - \ell - k_{1}$} \\ 
C \, \bQ_{2}^{4} \, 2^{-k-r } & \text{if $k \geq p - \ell - k_{1}$ and $r > p - \ell - k_{1}.$} \end{cases} \label{eq:G8krj-Sc-R3}
\end{align}
We deduce from \eqref{eq:y8} that:
\[
\langle \mu,| \psi_{8, p-\ell}| \rangle \leq  C\, 2^{2(p-\ell)} \, \left( \bC_2^2+ \bQ_2^2\right)^{2}.
\]

\medskip
\textbf{Upper bound of $\langle \mu, |\psi_{9,p-\ell} |\rangle$}. 
We have:
\begin{equation}\label{eq:L2R3p9}
\langle \mu,| \psi_{9, p-\ell}| \rangle \leq   C\,  2^{4(p-\ell)} \,\sum_{k=2}^{p-\ell-1} \sum_{r=1}^{k -1} \sum_{j=0}^{r-1} \,  2^{-j } \, \Gamma_{k,r,j}^{[9]},
\end{equation}
with 
\[
\Gamma_{k,r,j}^{[9]} \leq 2^{-k-r } \,\langle \mu, \cq^j\cp \left(\cq^{r-j-1}|\cp \left( h_{p-\ell-r} \sot \cq^{k -r -1} \cp\left(h_{p-\ell-k}\otimes^2 \right)\right)|\sot |h_{p-\ell-j}| \right)  \rangle .
\]
For $r \leq k - k_{1} - 1$, we have $r \leq p - \ell - k_{1} - 1$ and:
\begin{align}\nonumber
\Gamma_{k,r,j}^{[9]} &\leq  C\,  2^{-k-r } \,\normm{\cp \left( h_{p-\ell-r} \sot \cq^{k -r -1}\cp\left(h_{p-\ell-k}\otimes^2 \right)\right)}\normm{h_{p-\ell-j}}\\ 
\label{eq:L2R3p9-1-1st} &\leq   C\,  2^{-k-r } \,\normm{ h_{p-\ell-r-1}}\normm{\cp\left(h_{p-\ell-k}\otimes^2 \right)}\normm{h_{p-\ell-j}}\\
&\label{eq:L2R3p9-1}\leq   \begin{cases} C\, \bC_{2}^{4} \,  2^{-k-r } \, \alpha^{4(p-\ell)} \alpha^{-2(k+r)} & \text{if $p - \ell - k > k_{1}$}\\
 C \bC_{2}^{2} \bQ_{2}^{2} \, 2^{- k - r} \alpha^{2(p - \ell - r)} & \text{if $p - \ell - k \leq k_{1},$}
\end{cases}
\end{align}
where we used \reff{eq:mp} for the first inequality; \reff{eq:Hil-SchP-1} as $p-\ell -r \geq  k_{1} + 1$ and $k-r-1\geq k_{1}$ for the  second; and \reff{eq:hk} (two times) and \eqref{eq:phk} (one time) for the last. For $r > k - k_{1} - 1$ and $k \leq p - \ell - k_{1} - 1$, we have:
\begin{align}
\nonumber \Gamma_{k,r,j}^{[9]} & \leq   C\,  2^{- k - r } \, \normm{\cp \left( h_{p - \ell - r} \sot \Qq^{k - r - 1} \cp\left(h_{p-\ell-k}\otimes^2 \right)\right)}\normm{h_{p-\ell-j}}\\ 
&\leq  C\,  2^{- k - r } \, \normm{ h_{p - \ell - r - k_{1}}}\normm{h_{p - \ell - k - k_{1}}}^2 \normm{h_{p-\ell-j}} \label{eq:L2R3p9-2-1st}\\ 
\nonumber &\leq  C\, \bC_{2}^{4} \,   2^{- k - r } \, \alpha^{4(p-\ell)} \, \alpha^{- 2 (r - k)},
\end{align}
where we used \reff{eq:mp} for the first inequality;
\reff{eq:Hil-SchP-2}\footnote{Notice this is the only place in the proof
  of Theorem \ref{thm:flx}  where we use \reff{eq:Hil-SchP-2}.}  as $p-\ell -k \geq  k_{1} + 1$  for the  second; and \reff{eq:hk} (three times) for the last. For $r > k - k_{1} - 1$ and $k > p - \ell - k_{1} - 1$, we have:
\begin{align}
\nonumber \Gamma_{k,r,j}^{[9]} &\leq  C\,  2^{-k-r } \,\normm{\cp \left( h_{p-\ell-r} \sot \cq^{k -r -1}\cp\left(h_{p-\ell-k}\otimes^2 \right)\right)}\normm{h_{p-\ell-j}}\\ 
\nonumber&\leq   C\,  2^{- k - r} \,\normm{\cp \left(|\cq^{p - \ell - r - 1} \tilde f_{\ell,n}| \sot \Qq^{p - \ell - 2} \left(\Qq \tilde{f}_{\ell,n}^2 \right)\right)} \normm{h_{p-\ell-j}}\\ 
\label{eq:L2R3p9-3-1st} &\leq   C\, \bQ_{2}^{2} \,  2^{- k - r} \, \normm{\Qq^{p - \ell - r - 1} \tilde{f}_{\ell,n}} \normm{h_{p-\ell-j}}\\ 
\nonumber &\leq   C\,  2^{- k - r} \, \bQ_{2}^{2} \, \bC_{2}^{2}\, \alpha^{2(p-\ell)} \alpha^{- r - j}, 
\end{align}
where we used  \reff{eq:mp} for the first inequality,
\reff{eq:majo-pfxf} (with $f$ replaced by $\tilde{f}_{\ell,n}$) for the
second, \eqref{eq:q2-mu-f} for the third and \reff{eq:hk} (two times)   for the last. We then deduce from \reff{eq:L2R3p9} and the
computations thereafter, that:
\[
\langle \mu,| \psi_{9, p-\ell}| \rangle \leq   C\,   2^{2(p-\ell)} \, \bC_2^2 \, \left(\bC_2^2+ \bQ_2^2\right).
\]
\medskip

In conclusion, we get that:
\begin{align*}
R_3(n) &\leq  C\, n^3 2^{-n -p} \left( c_4^4(\bF_n) + \sum_{\ell=0}^{p - k_{1} - 2} \sum_{j=1}^9 \langle \mu, \psi_{j,p-\ell}\rangle\right)\\
&\leq  C\, n^3 2^{-n -p} \left( c_4^4(\bF_n) + \sum_{\ell=0}^{p - k_{1} - 2} \left[2^{p-\ell} (c_4^4(\bF_n) + c_6^3(\bF_n) \bC_2)  +  2^{2(p-\ell)} ( \bC^2_2+ \bQ_2^2)^{2}\right]\right)\\
&\leq  C\, n^3 \left(  2^{-n(1- 2\rho)} + 2^{(n-p)}  ( \bC^2_2 + \bQ_2^2)\right) (\bC_2^2 + \bQ_{2}^{2}), 
\end{align*}
where, we used  \eqref{eq:R3-L2-9}, \eqref{eq:L2R3-l=p} and
\eqref{eq:R3L2-Em4} for the first inequality; and $c_4^4(\bF_n)\leq  C\,
\bC_2^2 \, 2^{2n \rho}$ and  $c_6^3(\bF_n)\leq  C\,
\bC_2 \, 2^{2n \rho}$ with $\rho\in (0,1/2)$ thanks to Remark
\ref{rem:fln} and (i) from Assumption \ref{hyp:f} for the last one. 
As $\rho\in (0,1/2)$ by  Assumption \ref{hyp:f} (i), we deduce  that
$\lim_{n\rightarrow \infty } R_3(n)=0$.  
\end{proof}

\medskip

We can  now use  Theorem 3.2 and  Corollary 3.1, p.~58, and  the Remark p.~59  from  \cite{hh:ml}  to  deduce  from  Lemmas \ref{lem:cv-var-s-L2} and \ref{lem:cv-R3-L2-n} that  $\Delta_n(\bF_{n})$ converges in distribution towards a Gaussian  real-valued   random  variable  with   deterministic variance $\sigma^{2}$ defined  by \reff{eq:limf-ln}.  Using Remark \ref{rem:Nn=Nnk0} and Lemma  \ref{lem:Nk0f=Dnf}, we then  deduce Theorem \ref{thm:flx}.

\section{Proof of Theorem \ref{thm:flx} in the critical case ($2\alpha^2=1$)}\label{sec:proof-crit-L2}

We keep notations from Section  \ref{sec:proof-sub-L2}. We   assume that    Assumption \ref{hyp:Q1}  holds with $\alpha=1/\sqrt{2}$.   Let $(f_{ \ell,n}, \,  n\geq \ell\geq 0)$ be a  sequence of  function satisfying  Assumptions \ref{hyp:f} and \ref{hyp:QRjfln-crit}. We set $f_{\ell, n}=0$ for $\ell>n\geq 0$       and $\bF_n=(f_{\ell, n}, \ell\in \N)$. Recall the definition of $c_k(\bF)$ and $q_k(\bF)$ in \reff{eq:def-ck}.  Assumption \ref{hyp:f} (ii) gives that $\bC_{2}    =     \sup_{n\in     \N}    c_{2}(\bF_n)$ and $\bQ_{2} = \sup_{n\in \N} q_{2}(\bF_n)$ are finite. Recall from Remark \ref{rem:Nn=Nnk0} that the study of $N_{n,\emptyset}(\bF_{n})$ is reduced to that of $N_{n,\emptyset}^{[k_{0}]}(\bF_{n}).$ \medskip
 
\begin{lem}\label{lem:cvRk0L2crit-density}
Under the assumptions of Theorem \ref{thm:flx} ($2\alpha^2=1$), we get  $\lim_{n \rightarrow \infty} \EE[R_{0}^{k_{0}}(n)^{2}] = 0.$
\end{lem}

\begin{proof}
Assume $n-p\geq k_0$. We write:
\[
R_0^{k_0} (n)= |\G_n|^{-1/2}\, \sum_{k=k_0}^{n-p-1}\sum_{i\in \G_{k_0}}  M_{i\G_{k-k_0}}(\tilde f_{n-k,n}).
\]
We have that $
\sum_{i\in \G_{k_0}} \E[ M_{i\G_{k-k_0}}(\tilde f_{n-k,n})^2]  =
\E[M_{\G_{k_0}}(h_{k,n})]$,
where:
\[
h_{k,n}(x)=\E_x[ M_{\G_{k-k_0}}(\tilde f_{n-k,n})^2].
\] 
We deduce from \reff{eq:EGn-nu0}, that $\E[M_{\G_{k_0}}(h_{k,n})]\leq  C \langle \mu,  h_{k, n}\rangle$.
We have also that:
\begin{align*}
\langle \mu,  h_{k, n}\rangle = \E_\mu[ M_{\G_{k-k_0}}(\tilde
  f_{n-k,n})^2]
  & \leq  C \, 2^k \sum_{\ell=0}^k 2^\ell \normm{\cq^\ell \tilde
    f_{n-k,n}}^2\\
  & \leq  C \, 2^k\, c_{2}^2(\bF_{n}) \sum_{\ell=0}^{k}
    2^{\ell}\alpha^{2\ell}\\
  &\leq  C k \, 2^{k}\,c_{2}^2(\bF_{n}) ,
\end{align*}
where we used \reff{eq:EGn-nu0} for the first inequality (notice one can take $k_0=0$ in this case as we consider the expectation $\E_\mu$), \reff{eq:L2-erg}  in the second, and $2\alpha^2<1$ in the last. We deduce that:
\begin{equation}\label{eq:IRk0nL2sub}
\E[R_0^{k_0}(n)^2]^{1/2} \leq  |\G_n|^{-1/2}\, \sum_{k=k_0}^{n-p-1} \left(2^{k_0} \E\left[M_{\G_{k_0}}(h_{k,n})\right]\right)^{1/2} \leq  C\, \bC_{2} \, n^{1/2} \, 2^{- p/2}.
\end{equation}
As $\lim_{n\rightarrow \infty } p/n=1$, we get $\lim_{n\rightarrow \infty} n\, 2^{-p} = 0$ and this ends the proof using \eqref{eq:IRk0nL2sub}. 
\end{proof}

\begin{lem}\label{lem:cvR1L2crit-density}
Under the assumptions of Theorem \ref{thm:flx} ($2\alpha^2=1$), we get $\lim_{n \rightarrow \infty} \EE[R_{1}(n)^{2}] = 0.$
\end{lem}

\begin{proof}
Notice that \reff{eq:bQfln-crit} implies that:
\begin{equation}\label{eq:bQfln-crit-cons}
\normm{\cq \tilde f_{\ell,n}} \leq  \normm{\cq  f_{\ell,n}} \leq  \normm{\cq (| f_{\ell,n}|)} \leq \delta_{\ell, n}. 
\end{equation} 
We deduce that for $k\in \N$:
\begin{equation} \label{eq:bQtfln-crit}
\normm{\cq^k \tilde f_{\ell,n}} \leq  \alpha^{k-1} \delta_{\ell, n} \ind_{\{ k \geq 1\}} + \bC_2 \ind_{\{k=0\}}.
\end{equation}
We  set  for $p\geq  \ell \geq 0$, $n-p\geq k_0$ and $j\in \G_{k_0}$:
\begin{equation*}\label{eq:def-R1ln-L2}
R_{1,j} (\ell,n ) = \sum_{i\in j\G_{n-p-k_0}}\E\left[N_{n,i}^\ell(f_{\ell,n}) |\,\cf_i\right],
\end{equation*}
so that $R_1(n)=\sum_{\ell=0}^{p} \sum_{j\in \G_{k_0}} R_{1,j} (\ell,n )
$. We have for $i\in \G_{n-p}$:
\begin{equation}\label{eq:ENlni}
|\G_n|^{1/2} \E\left[N_{n,i}^\ell(f_{\ell,n}) |\, \cf_i\right] = \E\left[M_{i\G_{p-\ell}}(\tilde f_{\ell,n}) |X_i\right] = \E_{X_i} \left[M_{\G_{p-\ell}}(\tilde f_{\ell,n})\right] = |\G_{p-\ell}|\, \cq^{p-\ell} \tilde f_{\ell,n}(X_i),
\end{equation}
where we used definition \reff{eq:def-Nnil} of $N_{n,i}^\ell$ for the first equality, the Markov property of $X$ for the second and  \reff{eq:Q1} for the third. Using \reff{eq:ENlni}, we get for $j\in \G_{k_0}$: 
\[
R_{1,j}(\ell,n)= |\G_n|^{-1/2} \,  |\G_{p-\ell}|\, M_{j\G_{n-p-k_0}}(\cq^{p-\ell} \tilde f_\ell).
\]
We deduce from the Markov property of $X$ that $\E[R_{1,j}(\ell,n)^2|\, \cf_j]=2^{-n+2(p-\ell)} \, h_{\ell, n} (X_j)$ with $h_{\ell,n}(x)=\E_x\left[M_{\G_{n-p-k_0}} (\cq^{p-\ell} \tilde f_\ell)^2\right]$. 
Using \reff{eq:EGn-nu0}, we get: 
\[
\sum_{j\in \G_{k_0}} \E[R_{1,j}(\ell,n)^2] = 2^{-n+2(p-\ell)} \,
\E\left[M_{\G_{k_0}} (h_{\ell, n})\right]  
\leq C 2^{-n+2(p-\ell)} \,  \langle \mu,  h_{\ell,n}\rangle . 
\]
Using \reff{eq:EGn-nu0}, we have:
\begin{equation*}\label{eq:Imu-hln}
\langle \mu,  h_{\ell, n} \rangle = \E_\mu \left[M_{\G_{n-p-k_0}}(\cq^{p-\ell} \tilde f_{\ell,n})^2\right]  \leq   C  \, 2^{n-p} \sum_{k=0}^{n-p-k_0} 2^{k} \,  \normm{\cq^k \cq^{p-\ell} \tilde f_{\ell,n}} ^2. 
\end{equation*}
Using \eqref{eq:L2-erg} and \reff{eq:bQtfln-crit}, the latter inequality implies that:
\begin{align*}
\langle \mu,  h_{\ell, n} \rangle
& \leq   C  \, 2^{n-p} \sum_{k=0}^{n-p-k_0} 2^{k} \,  \normm{\cq^k
  \cq^{p-\ell} \tilde f_{\ell,n}} ^2 \\
&\leq  C \,2^{n-p} \alpha^{2(p-\ell)} \sum_{k=0}^{n-p-k_0} 2^{k}
  \alpha^{2k-2} \normm{\cq \tilde f_{\ell,n}} ^2\ind_{\{k+p-\ell \geq
  1\}}
+ C\, 2^{n-p} \normm{f_{\ell, n}}^2 \ind_{\{\ell=p\}}\\
&\leq  C \, (n-p) \, 2^{n-2p + \ell} \, \delta_{\ell, n}^2
+ C\, 2^{n-p} \bC_2^2 \ind_{\{\ell=p\}}.
\end{align*}
Using the following inequality,
\begin{equation*}\label{eq:IR1nL2sub}
\E\left[R_1(n)^2\right]^{1/2} \leq  \sum_{\ell=0}^{p}  \left(2^{k_0} \sum_{j \in \G_{k_0}} \E\left[R_{1,j}(\ell, n)^2\right] \right)^{1/2},
\end{equation*}
we have
\begin{align*}
   \EE[R_{1}(n)^{2}]^{1/2} 
& \leq C\, \sum_{\ell=0}^p \left( (n-p) 2^{-\ell}\,  \delta_{\ell, n}^2 +
  \ind_{\{\ell=p\}} 2^{-p} \bC_2^2 \right)^{1/2} \\
&\leq  C\, \left(   2^{-p/2} \bC_2 + \sqrt{n}\sum_{\ell=0}^ n 2^{-\ell/2}\,
  \delta_{\ell, n}  \right). 
\end{align*}
Then use \reff{eq:cv-nldln} to conclude. 
\end{proof}

\begin{lem}\label{lem:cvR2L2crit-density}
Under the assumptions of  Theorem \ref{thm:flx} ($2\alpha^2=1$), we get $\lim_{n \rightarrow \infty} \EE[ R_{2}(n)] = 0.$
\end{lem}

\begin{proof}
Using  \eqref{eq:Q1}, we have:
\begin{align}
\nonumber \E\left[R_2(n)\right] & = |\G_n|^{-1} \, \sum_{i\in \G_{n-p}} \E\left[ \E\left[\sum_{\ell=0}^{p} M_{i\G_{p-\ell}}(\tilde f_{\ell,n})
  |X_i\right]^2\right]\\
\nonumber &= |\G_n|^{-1} \, \sum_{i\in \G_{n-p}} \E\left[\left(\sum_{\ell=0}^{p} \E_{X_i} \left[M_{\G_{p-\ell}}(\tilde f_{\ell,n})\right]\right)^2\right]\\
\nonumber &= |\G_n|^{-1} \,  |\G_{n-p}|\, \cq^{n-p} \langle \nu, \left(\Big(\sum_{\ell=0}^{p}  |\G_{p-\ell}|\, \cq^{p-\ell} \tilde f_{\ell,n}\Big)^2\right)\rangle.
\end{align}

Next, using (iii) from Assumption \ref{hyp:DenMu} and \eqref{eq:bQtfln-crit}, we deduce that:
\begin{equation*}
\E\left[R_2(n)\right] \leq C 2^{-p} \left(\sum_{\ell=0}^{p}  |\G_{p-\ell}|\, \norm{\cq^{p-\ell} \tilde f_{\ell,n}}_{L^2(\mu)}\right)^2 \leq  C\, \bC_{2}^{2} \, 2^{-p} \,  + \, C \, \left(\sum_{\ell = 0}^{p-1} 2^{-\ell/2} \delta_{\ell,n}\right)^{2}.
\end{equation*}
Now, the result follows using the fact that $\lim_{n\rightarrow}p = \infty$ and the dominated convergence theorem.
\end{proof}

We now consider the limit of $V_2(n)$. 

\begin{lem}\label{lem:cvV2-L2-density}
Under   the  assumptions   of  Theorem   \ref{thm:flx} ($2\alpha^2=1$),  we  get $\lim_{n\rightarrow \infty } V_2(n) =0$ in
probability. 
\end{lem}

\begin{proof}
To  prove that $\lim_{n\rightarrow \infty } V_2(n)=0$ in
probability, we give a closer look at the proof of
\reff{eq:majo-L2A2}. Using $2\alpha^{2} = 1$, we get that
the upper
bound in \reff{eq:majo-L2A6}  can be replaced by 
$C c_2^2(\bF)\left(c_2(\bF)+q_2(\bF)\right) ^2 \, (n-p)2^{-(n-p)}$
and  the upper
bound in \reff{eq:majo-L2A5}  can be replaced by 
$C c_2^2(\bF)q_2^2(\bF) \, (n-p) 2^{-(n-p)}$. 
As $V_2=V_6+V_5$, we deduce that (compare with \reff{eq:majo-L2A2}):
\begin{align*}
\E\left[\left( V_2(n)-  H_2^{[n]}(\bF_n)\right)^2\right]
& \leq C \, \left(c_2^4(\bF_n)+ c_2^2(\bF_n)\, q_2^2(\bF_n)\right)  \, (n-p)\,
2^{-(n-p)} \\
& \leq C \, \left(\bC_2^4+ \bC_2^2\, \bQ_2^2\right)  \,
(n-p)\, 2^{-(n-p)},  
\end{align*}
with  $H_2^{[n]}(\bF_n)=H_5^{[n]} (\bF_n)+ H_6^{[n]}
(\bF_n)$. Since according  to (ii) in
Assumption, \ref{hyp:f} $\bC_2$ and $\bQ_2$ are finite, 
 we deduce  that 
$\lim_{n\rightarrow \infty } V_2(n) - H_2^{[n]}(\bF_n)=0$ in
probability. We now check that $\lim_{n\rightarrow \infty }
H_2^{[n]}(\bF_n)=0$. 
From 
\reff{eq:def-B5nsubL2}, we get  that $| H_5^{[n]}(\bF_n) |\leq  \sum_{k>\ell\geq 0}
2^{-\ell} |\langle \mu,   \tilde f_{k,n} \cq^{k-\ell} \tilde
f_{\ell,n}\rangle |$, and using \reff{eq:majoqtf} and \reff{eq:majL2qtf1}
which are a consequence
of Assumption \ref{hyp:f},  and the fact that
$\sum_{k>\ell\geq 0} 2^{-\ell} \alpha^{k-\ell}$ is finite, we get by
dominated convergence that $\lim_{n\rightarrow \infty }  H_5^{[n]}(\bF_n)=0$. 

Using \reff{eq:def-B6nsubL2}, we get  that:
\begin{equation}
   \label{eq:H6n-sec9}
|H_6^{[n]}(\bF_n)|
\leq \sum_{k>\ell\geq 0} \sum_{r=0}^{p-k-1} 2^{r-\ell} |\langle \mu,
\cp\left( \cq^r \tilde
    f_{k,n} \sot
\cq^{k-\ell+r} \tilde f_{\ell,n}  \right)\rangle|.
\end{equation}
Using \reff{eq:mp} and \reff{eq:bQfln-crit} (or more precisely 
\reff{eq:bQfln-crit-cons}) , we obtain:
\begin{align*}
|H_6^{[n]}(\bF_n)|
& \leq  \sum_{k>\ell\geq 0}\sum_{r=0}^{p-k-1}  2^{r-\ell} \normm{\cq^r 
\tilde     f_{k,n}}\normm{ \cq^{k-\ell+r} \tilde f_{\ell,n} }\\
& \leq  \sum_{k>\ell\geq 0}\sum_{r=0}^{p-k-1}  2^{r-\ell}
  \alpha^{k-\ell +2 r} \bC_2 
\normm{ \cq \tilde f_{\ell,n} }\\
& \leq  n \sum_{\ell=0}^n 2^{-\ell} \delta_{\ell, n}. 
\end{align*}
Then, use  \reff{eq:cv-nldln} to conclude. 
\end{proof}

\begin{lem}\label{lem:cvV1-L2-density}
Under   the  assumptions   of  Theorem   \ref{thm:flx} ($2\alpha^2=1$),  we   get $ \lim_{n\rightarrow \infty } V_1(n) =\sigma^2$ in probability.
\end{lem}

\begin{proof}
To  prove that $\lim_{n\rightarrow \infty } V_1(n)=\sigma^2$ in
probability, we give a closer look at the proof of
\reff{eq:majo-L2A1}.  
Using $2\alpha^{2} = 1,$ we get that the upper bound in \reff{eq:majo-L2A4}  can be replaced by $C \left(\bC_{2}^{4} + \bQ_{2}^{4}\right)  \, (n-p)2^{-(n-p)}$ and the upper bound in \reff{eq:majo-L2A3} can then be replaced by $C  \, c_4^4(\bF_{n})\, n2^{-n}$. As $V_1=V_4+V_3$, using (i) from Assumption \ref{hyp:f}, we deduce that (compare with \reff{eq:majo-L2A1}):
\begin{align*}
\E\left[\left( V_1(n)-  H_1^{[n]}(\bF_n)\right)^2\right] & \leq C \, \left(\bC_{2}^{4} + \bQ_{2}^{4}\right)  \, (n-p)\, 2^{-(n-p)}  + C\,  c_4^4(\bF_{n})\, n^22^{-n}\\
& \leq C \, \left(\bC_2^4+ \bQ_2^4\right)  \, (n-p)\, 2^{-(n-p)} + C\, n2^{-n(1-2\rho)},  
\end{align*}
with  $H_1^{[n]}(\bF_n)=H_4^{[n]} (\bF_n)+ H_3^{[n]} (\bF_n)$. This implies that  $\lim_{n\rightarrow \infty } V_1(n) - H_1^{[n]}(\bF_n)=0$ in probability. See the proof of Lemma \ref{lem:cv-var-s-L2} to get that $\lim_{n\rightarrow \infty } H_3^{[n]}(\bF_n)=\sigma^2$. Recall \reff{eq:def-B4nsubL2} for the definition of $H_4^{[n]}(\bF)$. We have:
\begin{align*}
|H_4^{[n]}(\bF_n)| &\leq  \sum_{\ell\geq 0; \, k\geq 0}  2^{k-\ell} |\langle \mu, \cp \left(\cq^k \tilde f_{\ell, n} \otimes ^2\right) \rangle|\,\ind_{\{\ell+k <p\}} \\
&\leq  \sum_{\ell\geq 0; \, k\geq 1}  2^{k-\ell} \alpha^{2k} \norm{Q |\tilde f_{\ell, n}|}_\infty  \,\ind_{\{\ell+k <p\}}  +   \sum_{\ell=0}^{p-1}  2^{-\ell}  |\langle \mu, \cp( \tilde f_{\ell, n} \otimes ^2) \rangle|\\
&\leq C n \sum_{\ell=0}^n 2^{-\ell} \delta_{\ell, n} +  C(1+\Delta) \sum_{\ell=0}^n 2^{-\ell}\delta_{\ell,n}.
\end{align*}
Thanks to \reff{eq:cv-nldln} from Assumption \ref{hyp:QRjfln-crit}, we
get
$\lim_{n\rightarrow \infty }
H_4^{[n]}(\bF_n)=0$, and thus $\lim_{n\rightarrow \infty }
H_1^{[n]}(\bF_n)=\sigma^2$. This finishes the proof. 
\end{proof}

As a conclusion of Lemmas \ref{lem:cvR2L2crit-density},
\ref{lem:cvV2-L2-density} and \ref{lem:cvV1-L2-density} and since
$V(n)=V_1(n)+2V_2(n)- R_2(n)$ (see \reff{eq:def-V}), we deduce that
  $ \lim_{n\rightarrow \infty } V(n) =\sigma^2$ in probability.
\medskip

We now check the Lindeberg condition using a fourth moment
condition. Recall $R_3(n)=\sum_{i\in \G_{n-p_n}}
\E\left[\Delta_{n,i}(\bF_{n})^4\right]$ defined in \reff{eq:def-R3}.

\begin{lem}\label{lem:cvG-L2-b-dens}
Under the assumptions of Theorem \ref{thm:flx} ($2\alpha^2=1$),  we get $\lim_{n\rightarrow\infty } R_3(n)\!=0.$
\end{lem}
\begin{proof}
Following line by line the proof of Lemma \ref{lem:cv-R3-L2-n} with the same
notations and taking $\alpha
= 1/ \sqrt{2}$, we get that concerning $|\langle \mu, \psi_{i,p-\ell}
\rangle|$ or $\langle \mu, | \psi_{i,p-\ell}| \rangle$,  the bounds
for $i\in \{1, 3, 4\}$ are the same; the bounds  for $i\in \{2, 5,
6\}$ have an extra $(p-\ell)$ term, the bounds   for $i\in \{7, 8, 9\}$  have an
extra $(p-\ell)^2$ term. This leads to:
\[
R_3(n) \leq   C\,  n^5\, \left( 2^{-n}   c_4^4(\bF_n) + 2^{-n}  c_6^3(\bF_n)\, \bC_2 +  2^{-(n-p) } \,\bC_2^2( \bC_2^2 +  \bQ_2^2 )\right). 
\]
Then conclude as in the proof of Lemma \ref{lem:cv-R3-L2-n}.
\end{proof}

Then, we end the proof of Theorem \ref{thm:flx} with $2\alpha^2=1$ by arguing as in the (end of the) proof of Theorem \ref{thm:flx} with $2\alpha^{2} < 1$.

\section{Proof of Theorem \ref{thm:flx} in the super-critical case ($2\alpha^2>1$)}\label{sec:proof-Scrit-L2} 
We assume $\alpha\in (1/\sqrt{2}, 1)$. We  follow line  by  line  the proof  of  Theorem \ref{thm:flx}  in Section \ref{sec:proof-crit-L2}  with $\alpha  > 1/\sqrt{2}$  instead of $\alpha=1/\sqrt{2}$, and use notations from Sections  \ref{sec:proof-sub-L2}. We  recall                that $\bC_{2} = \sup\{     c_{2}(\bF_n),\,  n\in  \N\}$   and $\bQ_{2} = \sup_{n\in \N} q_{2}(\bF_n)$  are finite thanks to Assumption \ref{hyp:f} (ii). We  will denote $C$  any unimportant finite  constant which may  vary  line   to  line,  independent  on  $n$   and  $\bF_{n}$. Let $(p_{n}, n \in \NN)$ be an increasing sequence of elements of $\NN$ such that \eqref{eq:def-pn} holds.  When there is no ambiguity,  we write $p$ for $p_{n}.$

\begin{lem}\label{lem:cvRk0L2Scri-density}
Under the assumptions of Theorem \ref{thm:flx} ($2\alpha^2>1$), we get $\lim_{n \rightarrow \infty} \EE[R_{0}^{k_{0}}(n)^{2}] = 0.$ 
\end{lem}

\begin{proof}
Mimicking  the proof of Lemma \ref{lem:cvRk0L2crit-density}, we get, as
$\lim_{n\rightarrow \infty } p/n=1$:
\begin{equation*}
\lim_{n \rightarrow \infty} \EE[R_{0}^{k_{0}}(n)^{2}] \leq  C \lim_{n
  \rightarrow \infty} c_{2}^2(\bF_n) (2 \, \alpha^{2})^{n-p}\, 2^{-p} \leq C \bC_2^2 \lim_{n
  \rightarrow \infty} (2 \, \alpha^{2})^{n-p}\, 2^{-p} = 0.
\end{equation*}
\end{proof}

\begin{lem}\label{lem:cvR1L2Scri-density}
Under the assumptions of Theorem \ref{thm:flx} ($2\alpha^2>1$), we get
$\lim_{n \rightarrow \infty} \EE[R_{1}(n)^{2}] = 0.$ 
\end{lem}

\begin{proof}
Following the proof of Lemma \ref{lem:cvR1L2crit-density} with $\alpha^{2} > 1/2$, we get:
\begin{align*}
   \EE[R_{1}(n)^{2}]^{1/2} 
& \leq C\, \sum_{\ell=0}^p \left( 2^{-\ell}\, (2\alpha^{2})^{n-\ell}  \delta_{\ell, n}^2 +
  \ind_{\{\ell=p\}} 2^{-p} \bC_2^2 \right)^{1/2} \\
&\leq  C\, \left(   2^{-p/2} \bC_2  + \sum_{\ell=0}^ n 2^{-\ell/2}\,
(2\alpha^{2})^{(n-\ell)/2}  \delta_{\ell, n}  \right). 
\end{align*}
Then use \reff{eq:Cfln-Scrit} and dominated convergence theorem to conclude. 
\end{proof}

From Lemmas \ref{lem:cvRk0L2Scri-density} and \ref{lem:cvR1L2Scri-density}, it follows that
\begin{equation*}
\lim_{n \rightarrow \infty} \E[(N^{[k_0]}_{n, \emptyset}(\bF_n) - \Delta_n(\bF_n) )^2] = 0.
\end{equation*}

\begin{lem}\label{lem:cvR2L2Scri-density}
Under the assumptions of  Theorem \ref{thm:flx} ($2\alpha^2>1$), we get $\lim_{n \rightarrow \infty} \EE[ R_{2}(n)] = 0.$
\end{lem}
\begin{proof}
Following the proof of Lemma \ref{lem:cvR2L2crit-density}, we get
\begin{equation*}
\EE[R_{2}(n)] \leq C \bC_{2}^{2} 2^{-p} + (\sum_{\ell = 0}^{n} 2^{-\ell/2} (2 \alpha^{2})^{(p - \ell)/2} \delta_{\ell,n})^{2}.
\end{equation*}
Then use \reff{eq:Cfln-Scrit}, $2\alpha^2>1$ and  dominated convergence
theorem to conclude.
\end{proof}

We now consider the limit of $V_2(n)$. 

\begin{lem}\label{lem:cvV2-L2-Scri-density}
Under   the  assumptions   of  Theorem   \ref{thm:flx} ($2\alpha^2>1$),
we  get $\lim_{n\rightarrow \infty } V_2(n) =0$ in 
probability. 
\end{lem}
\begin{proof}
Using with $\alpha > 1/\sqrt{2}$,  we get that the upper-bound in \reff{eq:majo-L2A5} can be replaced by  $C \bC_{2}^{2} \bQ_{2}^{2} \, \alpha^{2(n-p)}$. We get that for $r \geq k_{1}$:
\[
\normm{\cp (\cq^r \tilde f_{k,n} \sot \cq^{k-\ell+r} \tilde f_{\ell,
    n})}
\leq   C \alpha^{2r+k-\ell} \delta_{k, n}\delta_{\ell, n}
\leq  C (2 \alpha^2)^{-n}\, 2^{(\ell+k)/2} \alpha^{2(r+ k)} ,
\]
where  we  used  Assumption  \ref{hyp:QRjfln-crit} (vi)  for  the  first
inequality and Assumption \ref{hyp:fln-Scrit}  for the second.  Thus the
bound     \reff{eq:majo-A6_bis-c4q2}     can      be     replaced     by
$ C\, (2\alpha^2)^{-(n-p)} \alpha^j$. The term \reff{eq:majo-A6-c4q2} is
handled as in  the proof of Lemma  \ref{lem:cvV2-L2-density}. This gives
that \reff{eq:majo-A6-c4} can be replaced by $C \, \alpha^j$.  Therefore
the   upper   bound   in   \reff{eq:majo-L2A6}  can   be   replaced   by
$C   \,   \alpha^{2(n-p)}$.    As    $V_2=V_6+V_5$,   we   deduce   that
$\EE[(V_{2}(n)  -  H_{2}^{[n]}(\bF_{n}))^{2}] \leq  C  \alpha^{2(n-p)}$.
(Compare  with \reff{eq:majo-L2A2}  and replace  $\bF$ by  $\bF_n$.)  It
follows                                                             that
$\lim_{n  \rightarrow \infty}  V_{2}(n) -  H_{2}^{[n]}(\bF_{n}) =  0$ in
probability.  As in the proof of Lemma \ref{lem:cvV2-L2-density} we also
have  $\lim_{n  \rightarrow   \infty}H_{5}^{[n]}(\bF_{n})  =  0$.  Using
\reff{eq:mp} and \reff{eq:bQfln-crit-cons} , we deduce from
\reff{eq:H6n-sec9} and Assumption \ref{hyp:fln-Scrit}  that:
\begin{align*}
|H_{6}^{[n]}(\bF_{n})| 
&\leq C \sum_{0 \leq \ell < k \leq p} 2^{-\ell} \alpha^{k-\ell}
  \delta_{\ell,n} + C \sum_{0 \leq \ell < k \leq p} \sum_{r = 1}^{p - k
  - 1} 2^{r - \ell} \delta_{k,n} \delta_{\ell,n} \alpha^{2r + k - \ell}
  \\ 
&\leq C \sum_{0 \leq \ell < k \leq p} 2^{-\ell} \alpha^{k-\ell}
  (2\alpha^2)^{-(n-\ell)/2}  + C \sum_{0 \leq \ell < k \leq n} (2\alpha^2)^{-(n-p)}
  2^{-(\ell+k)/2}\\
& \leq  C\, (2\alpha^2)^{-(n-p)}. 
\end{align*}
 Since $H_{2}^{[n]}(\bF_{n}) = H_{5}^{[n]}(\bF_{n}) +
 H_{6}^{[n]}(\bF_{n})$, it follows that $\lim_{n \rightarrow \infty}
 |H_{2}^{[n]}(\bF_{n})| = 0$. We deduce  that $\lim_{n \rightarrow \infty}
 V_{2}(n) = 0$ in probability. 
\end{proof}

\begin{lem}\label{lem:cvV1-L2-Scri-density}
Under   the  assumptions   of  Theorem   \ref{thm:flx} ($2\alpha^2>1$),
we   get $ \lim_{n\rightarrow \infty } V_1(n) =\sigma^2$ in
probability. 
\end{lem}
\begin{proof}
We follow the proof of Lemma \ref{lem:cvV1-L2-density} with $\alpha >
1/\sqrt{2}$ and use the same trick as in the proof of Lemma
\ref{lem:cvV2-L2-Scri-density} based on Assumption
\ref{hyp:fln-Scrit}. We get, with the details left to the reader:
\[
\EE[(V_{4}(n) - H_{4}^{[n]}(\bF_{n}))^{2}] \leq C \alpha^{2(n-p)}.
\]

We set $g_{\ell,n} = \tilde{f}_{\ell,n}^{2}.$ From \eqref{eq:L2A3Hf-q}, we have for $j\in \{0, \ldots, n-p\}$:
\begin{align*}
\|\cq^j (A_{3,n}(\bF_{n}))\|_{L^{2}(\mu)} &\leq C \,  \sum_{\ell = 0}^{p}  2^{-\ell} \,  \normm{ \cq^{p-\ell} \tilde g_{\ell,n}} \ind_{\{j = 0\}} + C \, \sum_{\ell = 0}^{p}  2^{-\ell} \,\normm{ \cq^{j+p-\ell} \tilde g _{\ell, n}} \ind_{\{j \geq  1\}} \\
&= C \, 2^{-p} \, \|\tilde{g}_{p,n}\|_{L^{2}(\mu)} \ind_{\{j = 0\}} \, + \,  C \,  \sum_{\ell = 0}^{p - 1}  2^{-\ell} \,  \normm{ \cq^{p - \ell - 1} (\Qq\tilde g_{\ell,n})} \ind_{\{j = 0\}} \\
& \hspace{3cm} + \, C \,  \sum_{\ell = 0}^{p}  2^{-\ell} \,  \normm{ \cq^{j + p - \ell - 1} (\Qq\tilde g_{\ell,n})} \ind_{\{j \geq 1\}} \\
& \leq C \, 2^{-p} \, \|f_{\ell,n}^{2}\|_{L^{2}(\mu)} \, \ind_{\{j = 0\}} \, + \, C \, \sum_{\ell = 0}^{p-1} 2^{-\ell} \alpha^{p - \ell - 1} \|\Qq \tilde{g}_{\ell,n}\|_{L^{2}(\mu)} \ind_{\{j = 0\}} \\
& \hspace{3cm} + \, C \, \sum_{\ell = 0}^{p} 2^{-\ell} \alpha^{j + p - \ell - 1} \|\Qq \tilde{g}_{\ell,n}\|_{L^{2}(\mu)} \ind_{\{j \geq 1\}} \\
&\leq C \, 2^{-p} \, \|f_{\ell,n}\|_{L^{4}(\mu)}^{2} \, \ind_{\{j = 0\}} \, + \, C \, (\bQ_{2}^{2}+ \bC_2^2) \, \sum_{\ell = 0}^{p} 2^{-\ell} \alpha^{j + p - \ell} \\
&\leq C \bC_2^2) \, 2^{-p} \,  2^{n\rho} \, \ind_{\{j = 0\}}  \, + \, C \,  (\bQ_{2}^{2}+ \bC_2^2) \, \alpha^{j + p}, 
\end{align*}
where we used Remark \ref{rem:fln}, (ii) of Assumption \ref{hyp:f},  \eqref{eq:L2-erg} and \eqref{eq:q2-mu-f}. From the latter inequality, we get using \eqref{eq:V3-H3f} and \eqref{eq:majoA3}:
\[
\EE[(V_{3}(n) - H_{3}^{[n]}(\bF_{n}))^{2}] \leq C (2^{-(1-2\rho)n} +
\alpha^{2n}).
\]
The latter inequalities imply that $\lim_{n \rightarrow \infty} \EE[(V_{1}(n) - H_{1}^{[n]}(\bF_{n}))^{2}] = 0$, with $H_{1}^{[n]}(\bF_{n})=H_{4}^{[n]}(\bF_{n})+H_{3}^{[n]}(\bF_{n})$. From the proof of Lemma \ref{lem:cv-var-s-L2} we have
$\lim_{n \rightarrow \infty} H_{3}^{[n]}(\bF_{n}) = \sigma^{2}$. 
Next, we have
\begin{align*}
|H_4^{[n]}(\bF_n)| 
&\leq  \sum_{\ell\geq 0, \, k\geq 0}  2^{k-\ell} |\langle \mu, \cp
  \left(\cq^k \tilde f_{\ell, n} \otimes ^2\right) \rangle|\,
  \ind_{\{\ell+k <p\}} \\ 
& \leq C \sum_{\ell\geq 0, \, k\geq 1}  2^{k-\ell} \alpha^{2k} \norm{Q
  |\tilde f_{\ell, n}|}_\infty^{2}  \,\ind_{\{\ell+k <p\}}  +
  C\sum_{\ell=0}^{p-1}  2^{-\ell}  |\langle \mu, \cp ( \tilde f_{\ell, n}
  \otimes ^2) \rangle|\\ 
&\leq C \sum_{\ell\geq 0, \, k\geq 1}  2^{k-\ell} \alpha^{2k}
 (2\alpha^2)^{-(n-\ell)} \,\ind_{\{\ell+k <p\}}   
  + C(1+\Delta)\sum_{\ell=0}^{p-1} 2^{-\ell} (2 \alpha^2)^{-(n-\ell)/2} \\ 
&\leq C \, (2 \alpha^2)^{-(n-p)}   
 + C(1+\Delta)\, (2\alpha^2)^{-n/2},
\end{align*}
where we used \reff{eq:def-B4nsubL2}  and the
definition of $h_{\ell,k}$ therein for the first inequality;
\eqref{eq:mp} for the second; Assumption \ref{hyp:f} (iii), Assumption \ref{hyp:QRjfln-crit} (vi),
\reff{eq:majo-mPf22} and Assumption \ref{hyp:fln-Scrit} (twice) for the third. 
We deduce that   $\lim_{n \rightarrow \infty} |H_{4}^{[n]}(\bF_{n})| =
0$.  This
ends the proof. 
\end{proof}

We now check the Lindeberg condition. For that purpose, we have the following result. 
\begin{lem}\label{lem:cvG-L2-b-dens-2}
Under the assumptions of Theorem \ref{thm:flx} ($2\alpha^2>1$),  we have  $\lim_{n\rightarrow\infty } R_3(n)=0.$
\end{lem}

\begin{proof}
From \eqref{eq:R3-L2-9}, \eqref{eq:L2R3-l=p} and \eqref{eq:R3L2-Em4}, we have
\begin{equation*}\label{eq:IR3ScritL2}
R_{3}(n) \leq C n^{3} 2^{-n-p} c_{4}^{4}(\bF_{n}) + n^{3} 2^{-n-p} \sum_{\ell = 0}^{p} \sum_{j = 1}^{9} \langle \mu,\psi_{j,p-\ell} \rangle.
\end{equation*}

Now, we will bound above each term in the latter sum. For that purpose, we will follow line by line the proof of Lemma \ref{lem:cv-R3-L2-n} and we will intensively use \eqref{eq:bQfln-crit} and \eqref{eq:bQfln-crit-cons}. We will also use the fact that for all nonnegative sequence $(a_{\ell},\ell \in \NN)$ such that $\sum_{\ell \geq 0} a_{\ell} < \infty$, the sequence $(\sum_{\ell = 0}^{n} a_{\ell}(2\alpha^{2})^{n-\ell}\delta_{\ell,n}^{2}, n \in \NN)$ is bounded  as a consequence of the first part of  \eqref{eq:Cfln-Scrit} from Assumption \ref{hyp:fln-Scrit}. (Notice that by the second part of \eqref{eq:Cfln-Scrit} and the dominated convergence theorem, the latter sequence converges towards 0; but we shall not need this.). Recall from Assumption \ref{hyp:f} that $\rho \in (0,1/2).$
 
\medskip

\textbf{The term $n^{3} 2^{-n-p} c_{4}^{4}(\bF_{n}).$} From the first inequality in Remark \ref{rem:fln}, we have
\begin{equation*}
n^{3} 2^{-n-p} c_{4}^{4}(\bF_{n}) \leq C\, n^{3} 2^{-(1-2\rho )n - p}.
\end{equation*}

\medskip
\textbf{The term $n^{3} 2^{-n-p} \sum_{\ell = 0}^{p-3}  \langle \mu,|\psi_{1,p-\ell}| \rangle.$} Using \eqref{eq:L2R3p1} and Remark \ref{rem:fln}, we get:
\begin{equation*}
n^{3} 2^{-n-p} \sum_{\ell = 0}^{p}  \langle \mu,|\psi_{1,p-\ell}| \rangle \leq C \bC_{2}^{2} n^{3} 2^{-(1-2\rho)n}.
\end{equation*}

\medskip
\textbf{The term $n^{3} 2^{-n-p} \sum_{\ell = 0}^{p-3}  |\langle \mu,\psi_{2,p-\ell} \rangle|.$} Distinguishing  the case $k = p - \ell - 1$ and $k \leq p - \ell - 2$ in \eqref{eq:L2R3p2-astuce} and using Remark \ref{rem:fln} and \eqref{eq:bQfln-crit}, we get:
\begin{equation*}
|\langle \mu,\psi_{2,p-\ell} \rangle| \leq C\, 2^{p - \ell} 2^{2n\rho}  + C\, (2\alpha)^{2(p - \ell)} 2^{2n\rho} \delta_{\ell,n}^{2}. 
\end{equation*}
This implies that
\begin{align*}
  n^{3} 2^{-n-p} \sum_{\ell = 0}^{p-3}  \langle \mu,|\psi_{2,p-\ell}| \rangle
  &\leq C \,n^{3} 2^{-n(1-2\rho)} \sum_{\ell = 0}^{p - 3} 2^{-\ell}  + C\, n^{3} 2^{-n(1-2\rho )} \sum_{\ell = 0}^{p - 3} 2^{-\ell}\,   (2\alpha^{2})^{p - \ell} \delta_{\ell,n}^{2} \\
  & \leq C\, n^{3} 2^{-n(1-2\rho)}.
\end{align*}  

\medskip
\textbf{The term $n^{3} 2^{-n-p} \sum_{\ell = 0}^{p-3}  \langle \mu,|\psi_{3,p-\ell}| \rangle.$}
From \eqref{eq:L2R3p3-fln} we have
\begin{equation*}
n^{3} 2^{-n-p} \sum_{\ell = 0}^{p}  \langle \mu,|\psi_{3,p-\ell}| \rangle \leq C\, n^{3} 2^{-(1 - 2\rho)n} + C\, n^{3} 2^{- n + p}
\end{equation*}

\medskip
\textbf{The term $n^{3} 2^{-n-p} \sum_{\ell = 0}^{p-3}  \langle \mu,|\psi_{4,p-\ell}| \rangle.$} From \eqref{eq:L2R3p4-sub} we have
\begin{equation*}
\langle \mu,|\psi_{4,p-\ell}| \rangle \leq C\, 2^{2(p - \ell)} ((2\alpha^{2})^{p-\ell}\delta_{\ell,n}^{2})^{2} \leq C\, 2^{2(p - \ell)}
\end{equation*}
and thus
\begin{equation*}
n^{3} 2^{-n-p} \sum_{\ell = 0}^{p-3}  \langle \mu,|\psi_{4,p-\ell}| \rangle \leq C\, n^{3} 2^{- n + p}.
\end{equation*}

\medskip
\textbf{The term $n^{3} 2^{-n-p} \sum_{\ell = 0}^{p-3}  \langle \mu,|\psi_{5,p-\ell}| \rangle.$} From \eqref{eq:L2R3p5-Scrit} and distinguishing the case $k > p - \ell - k_{1} - 1$ (and then using (ii) of Assumption \ref{hyp:f}) and $k \leq p - \ell - 1$ (and then using  \eqref{eq:Hil-SchP-1}, \eqref{eq:L2-erg} with $\Qq\tilde{f}_{\ell,n}$ instead of $f$, \eqref{eq:bQfln-crit} and \eqref{eq:Cfln-Scrit} of Assumption \ref{hyp:fln-Scrit}),
we get $\langle \mu,|\psi_{5,p-\ell}| \rangle \leq C 2^{2(p - \ell)}$ and thus $$n^{3} 2^{-n-p} \sum_{\ell = 0}^{p-3}  \langle \mu,|\psi_{5,p-\ell}| \rangle \leq C n^{3} 2^{- n + p}.$$ 

\medskip
\textbf{The term $n^{3} 2^{-n-p} \sum_{\ell = 0}^{p-3}  \langle
  \mu,|\psi_{6,p-\ell}| \rangle.$} Very similarly, from
\eqref{eq:L2R3p6-Scrit}, we have $$n^{3} 2^{-n-p} \sum_{\ell = 0}^{p-3}
\langle \mu,|\psi_{6,p-\ell}| \rangle \leq C n^{3} 2^{- n + p}.$$   

\medskip
\textbf{The term $n^{3} 2^{-n-p} \sum_{\ell = 0}^{p-3}  |\langle
  \mu,\psi_{7,p-\ell} \rangle|.$} We set $g_{k,n} = \cp \left(
  h_{p-\ell-k} \sot  \cq^{p-\ell-k -1} (\tilde f_{\ell,n}^2) \right).$ Using
that $\langle  \mu, \cp(\ind \sot h_{p-\ell -r}\rangle=0$,
\eqref{eq:L2R3p7-G} and \eqref{eq:mp}, we obtain  
\begin{align}
\Gamma_{k,r}^{[7]} &= 2^{-k} |\langle \mu,\Pp(\Qq^{k-r-1}(\tilde{g}_{k,n}) \sot \Qq^{p - \ell - r - 1}(\tilde{f}_{\ell,n})) \rangle| \nonumber \\
&\leq C 2^{-k} \normm{\Qq^{k - r - 1}(\tilde{g}_{k,n})} \normm{\Qq^{p - \ell - 2 - r}(\Qq\tilde{f}_{\ell,n})}. \label{eq:G7krSc-R3}
\end{align}
For $k \geq p - \ell - k_{1} - 1$, we have
\begin{equation*}
\Gamma_{k,r}^{[7]} \leq  C\, 2^{- k} \normm{\Qq^{k - r - 1}(\tilde{g}_{k,n})} \normm{\Qq^{p - \ell - 2 - r}(\Qq\tilde{f}_{\ell,n})} \leq C\, 2^{- (p - \ell)} \alpha^{2(p - \ell -
  r)} \delta_{\ell,n}^2 2^{2\rho n}, 
\end{equation*}
where we used \eqref{eq:L2-erg}, \eqref{eq:bQfln-crit-cons} and the following inequalities:
\begin{equation*}
\normm{\cp (\tilde f_{\ell,n} \sot  (\tilde f_{\ell,n}^2))} \leq C  \,
\delta_{\ell,n}2^{2n\rho} 
\end{equation*}
which is a consequence of (i) and (iii) of Assumption \ref{hyp:f},
\eqref{eq:bQfln-crit} from Assumption \ref{hyp:QRjfln-crit} and
\begin{equation*}
\normm{\cp (|f_{\ell,n}| \sot  f_{\ell,n}^2)} \leq C
\normm{\Qq(|f_{\ell,n}| ^{3})^{1/3} \Qq(|f_{\ell,n}|^{3})^{2/3}} = \!C
\normm{\Qq(|f_{\ell,n}|^{3})}\leq  C \delta_{\ell, n} 2^{2n\rho},
\end{equation*}
and
\begin{equation*}
\normm{\cp (\Qq^{p-\ell-k-1}\tilde f_{\ell,n} \sot  (\Qq^{p-\ell-k-1}\tilde f_{\ell,n}^2))} \leq C  \, \bQ_{2}^{2} \, \delta_{\ell,n} \quad \text{for all $k \leq p - \ell - k - 2.$}
\end{equation*}
which is a consequence of (ii) of Assumption \ref{hyp:f}, \eqref{eq:bQfln-crit} from Assumption \ref{hyp:QRjfln-crit}. 
Next, for $k \leq p - \ell - k_{1} - 2,$ using \eqref{eq:G7krSc-R3} for the first inequality, \eqref{eq:L2-erg} and \eqref{eq:bQfln-crit-cons} twice (for the second and the last inequality) and (ii) of Assumption \ref{hyp:f} for the third inequality, we obtain:
\begin{align*}
\Gamma_{k,r}^{[7]} &\leq  C\, 2^{-k} \normm{\Qq^{k  -  r  - 1}(\tilde{g}_{k,n})} \normm{\Qq^{p - \ell - 2 - r}(\Qq\tilde{f}_{\ell,n})}\\
& \leq C \, 2^{-k} \alpha^{k - r} \|g_{k,n}\|_{L^{2}(\mu)} \alpha^{p - \ell - r} \delta_{\ell,n} \\
& \leq C \, 2^{-k} \alpha^{k - r + (p - \ell - r)} \, \delta_{\ell,n} \, \|\Qq^{p - \ell - k - 1} (\Qq \tilde{f}_{\ell,n})\|_{L^{2}(\mu)} \, \|\Qq(\tilde{f}_{\ell,n}^{2})\|_{\infty}\\
&\leq C\, 2^{-k} \alpha^{2(p - \ell - r)} \delta_{\ell,n}^{2}.
\end{align*}
Thanks to Assumption \ref{hyp:fln-Scrit}, it follows from the foregoing that
\begin{equation*}
|\langle \mu,\psi_{7,p-\ell} \rangle| \leq C 2^{3(p - \ell)} \sum_{k = 1}^{p - \ell - 1} \sum_{r = 0}^{k - 1} 2^{-r} \Gamma^{[7]}_{k,r} \leq C 2^{2(p - \ell)} (\alpha^{2(p-\ell)} \delta_{\ell,n}^2 2^{2n\rho} + (2\alpha^{2})^{p-\ell}\delta_{\ell,n}^{2}),
\end{equation*}
and thus, we obtain
\begin{equation*}
n^{3} 2^{-n-p} \sum_{\ell = 0}^{p-3}  |\langle \mu,\psi_{7,p-\ell} \rangle| \leq C n^{3} (2^{-(1-2\rho )n} + 2^{- n + p}). 
\end{equation*} 
  
\medskip
\textbf{The term $n^{3} 2^{-n-p} \sum_{\ell = 0}^{p-3}  \langle \mu,|\psi_{8,p-\ell}| \rangle.$} From \eqref{eq:G8krj-Sc-R3} we have, 
\begin{align*}
\Gamma^{[8]}_{k,r,j} & \leq C \, 2^{- k - r} \alpha^{4(p - \ell)} \alpha^{- 2(k + r)}\delta_{\ell,n}^{4} \, \ind_{\{k \leq p - \ell - k_{1} - 2\}}  \\ 
& + C \, 2^{- k - r} \alpha^{2(p - \ell - r)} \delta_{\ell,n}^{2} \, \ind_{\{k \geq p - \ell - k_{1} - 1; \, r \leq p - \ell - k_{1} - 2\}}  \, + \, C \, 2^{2(p-\ell)} \, \ind_{\{k \geq p - \ell - k_{1} - 1; \, r \geq p - \ell - k_{1} - 1\}},
\end{align*}
where we use \eqref{eq:Hil-SchP-1}, \eqref{eq:L2-erg} and \eqref{eq:bQfln-crit-cons} for the cases $k \leq p - \ell - k_{1} - 2$ and $\{k \geq p - \ell - k_{1} - 1; \, r \leq p - \ell - k_{1} - 2\}$, and we used in addition \eqref{eq:IPflnpsi5} for the case $\{k \geq p - \ell - k_{1} - 1; \, r \geq p - \ell - k_{1} - 1\}.$ From \eqref{eq:y8}, the latter inequality implies that
\begin{equation*}
n^{3} 2^{-n-p} \sum_{\ell = 0}^{p-3}  \langle \mu,|\psi_{8,p-\ell}|\rangle  \leq \, C \, n^{3} 2^{- n + p}. 
\end{equation*}   

\medskip
\textbf{The term $n^{3} 2^{-n-p} \sum_{\ell = 0}^{p-3}  \langle \mu,|\psi_{9,p-\ell}| \rangle.$} From \eqref{eq:L2R3p9-1-1st}, \eqref{eq:L2R3p9-2-1st} and \eqref{eq:L2R3p9-3-1st}, using \eqref{eq:L2-erg}, \eqref{eq:bQfln-crit-cons} and \eqref{eq:IPflnpsi5}, we obtain
\begin{align*}
\Gamma^{[9]}_{k,r,j} \leq C \, 2^{- k - r} \alpha^{4(p - \ell)} \alpha^{- 2(k + r)} \delta_{\ell,n}^{4} \, \ind_{\{k \leq p - \ell - k_{1} - 1\}} \, + \, C \, 2^{- k - r} \, \alpha^{2(p - \ell)} \, \delta_{\ell,n}^{2} \, \alpha^{- r - j} \, \ind_{\{k \geq p - \ell - k_{1}\}}. 
\end{align*}
Using \eqref{eq:L2R3p9}, it then follows that
\begin{equation*}
n^{3} 2^{-n-p} \sum_{\ell = 0}^{p-3}  \langle \mu,|\psi_{9,p-\ell}| \rangle \leq C n^{3} 2^{- n + p}.
\end{equation*}
From the previous bounds, we deduce that $\lim_{n \rightarrow \infty} R_{3}(n) = 0.$ 
\end{proof}
Finally, arguing as in the (end of the) proof of Theorem \ref{thm:flx} (sub-critical case), we end the proof of Theorem \ref{thm:flx} in the
super-critical case. 

\section{Appendix}

In this section, we recall useful results on BMC which are recalled in  \cite{BD2}.

\begin{lem}\label{lem:Qi}
Let $f,g\in \cb(S)$, $x\in S$ and $n\geq m\geq 0$. Assuming that all the quantities below are well defined,  we have:
\begin{align} 
\label{eq:Q1} \E_x\left[M_{\G_n}(f)\right] &=|\G_n|\, \cq^n f(x)= 2^n\, \cq^n f(x) ,\\ \label{eq:Q2} \E_x\left[M_{\G_n}(f)^2\right] &=2^n\, \cq^n (f^2) (x) + \sum_{k=0}^{n-1} 2^{n+k}\,   \cq^{n-k-1}\left( \cp\left(\cq^{k}f\otimes \cq^k f \right)\right) (x),\\ \label{eq:Q2-bis} \E_x\left[M_{\G_n}(f)M_{\G_m}(g)\right] &=2^{n} \cq^{m} \left(g \cq^{n-m} f\right)(x)\\ \nonumber &\hspace{2cm} + \sum_{k=0}^{m-1} 2^{n+k}\, \cq^{m-k-1} \left(\cp\left(\cq^k g \sot \cq^{n-m+k} f\right) \right)(x). 
\end{align}
\end{lem}

\begin{lem}\label{lem:L2MG}
Let $X$  be a BMC  with kernel  $\cp$ and initial  distribution $\nu$ such that (iii)  from Assumption \ref{hyp:DenMu} (with $k_0\in  \N$) is in force. There exists a finite constant $C$, such that for all $f\in \cb_+(S)$ all $n\geq k_0$, we have: 
\begin{equation}\label{eq:EGn-nu0}
|\G_n|^{-1} \E[M_{\G_n} (f)]\leq  C  \norm{f}_{L^1(\mu)}\quad\text{and}\quad |\G_n|^{-1} \E\left[M_{\G_n} (f)^2\right]\leq  C  \sum_{k=0}^{n} 2^{k} \normm{\cq^kf}^2. 
\end{equation}
 \end{lem}

We also give some bounds on $\E_x\left[M_{\G_{n}}(f) ^4\right]$,  see the proof of Theorem 2.1 in \cite{BDG14}.  We will use the notation:
\[
g\otimes^2=g\otimes g.
\]
\begin{lem}\label{lem:M4}
There exists a  finite constant $C$ such that for  all $f\in \cb(S)$,$n\in \N$ and  $\nu$ a probability measure on $S$,  assuming that all the quantities below are well defined, there exist functions $\psi_{j,n}$ for $1\leq j\leq 9$ such that:
\[
\E_\nu\left[M_{\G_n}(f)^4\right]= \sum_{j=1}^9 \langle \nu, \psi_{j, n}\rangle,
\]
and, with $h_{k}= \cq^{k - 1} (f) $ and (notice that either $|\psi_j|$ or $|\langle \nu, \psi_j \rangle|$ is bounded), writing  $\nu   g = \langle \nu  ,   g   \rangle$:
\begin{align*}
| \psi_{1, n}| &\leq C \,2^n \cq^n(f^4),\\
| \nu \psi_{2, n}| &\leq C\,  2^{2n}\, \sum_{k=0}^{n-1} 2^{-k} |\nu \cq^k \cp  \left( \cq^{n-k - 1}( f^3) \sot h_{n- k} \right)|,\\
|\psi_{3, n}| &\leq C 2^{2n} \sum_{k=0}^{n-1} 2^{-k}\,  \cq^k \cp \left( \cq^{n-k - 1} (f^2) \otimes^2\right),\\
|\psi_{4, n}| &\leq C \, 2^{4n} \, \cp \left( |\cp(h_{n-1}\otimes^2)\otimes^2|\right), \\
|\psi_{5, n}| &\leq C\,  2^{4n} \, \sum_{k=2}^{n-1} \sum_{r=0}^{k -1}  2^{-2k-r }  \cq^r \cp \left( \cq^{k -r- 1} |\cp (h_{n- k} \otimes^2)|\otimes^2 \right),\\
|\psi_{6, n}| &\leq C\, 2^{3n} \, \sum_{k=1}^{n-1} \sum_{r=0}^{k -1} 2^{-k-r }  \cq^r| \cp \left(\cq^{k -r-1}\cp \left( h_{n-k} \otimes^2 \right)\sot\cq^{n-r-1}(f^2) \right)|,\\
|\nu \psi_{7, n}| &\leq  C\, 2^{3n} \,   \sum_{k=1}^{n-1} \sum_{r=0}^{k -1} 2^{-k-r } |\nu \cq^r \cp \left(\cq^{k -r-1}\cp \left( h_{n-k} \sot  \cq^{n-k -1} (f^2) \right)\sot h_{n-r} \right)|,\\
|\psi_{8, n}| &\leq C\,  2^{4n} \, \sum_{k=2}^{n-1} \sum_{r=1}^{k -1} \sum_{j=0}^{r-1} 2^{-k-r-j } \cq^j \cp \left(|\cq^{r-j-1}\cp \left( h_{n-r} \otimes^2 \right)|\sot |\cq^{k-j-1}\cp \left( h_{n-k} \otimes^2 \right)|\right),\\
|\psi_{9, n}| & \leq C\,  2^{4n} \, \sum_{k=2}^{n-1} \sum_{r=1}^{k -1} \sum_{j=0}^{r-1}2^{-k-r-j } \cq^j |\cp \left(\cq^{r-j-1}|\cp \left( h_{n-r} \sot \cq^{k -r -1}\cp\left(h_{n-k}\otimes^2 \right)\right) \sot h_{n-j} \right)|.
\end{align*}
\end{lem}


\bibliographystyle{abbrv}
\bibliography{biblio}

\end{document}